\documentclass[ijoc,sglanonrev]{class_paper}
\RequirePackage{tgtermes}
\RequirePackage{newtxtext}
\RequirePackage{newtxmath}
\RequirePackage{bm}
\RequirePackage{endnotes}

\OneAndAHalfSpacedXII 

\usepackage{tikz}
\usepackage{verbatim}
\usepackage{adjustbox}
\usepackage{amsmath}
\allowdisplaybreaks[4]
\newcommand{\SMIPD}{{\texttt{2-SMIPD}}}
\newcommand{\NDFPP}{{\texttt{NDFPP}}}
\newcommand{\DEOP}{{\texttt{DE-OP}}}
\newcommand{\SAA}{{\texttt{SAA}}}
\newcommand{\MILP}{{\texttt{MILP}}}
\newcommand{\EV}{{\texttt{EV}}}
\newcommand{\MINLP}{{\texttt{MINLP}}}
\newcommand{\MICP}{{\texttt{MICP}}}
\newcommand{\GBD}{{\texttt{GBD}}}
\newcommand{\APS}{{\texttt{APS}}}
\newcommand{\VSS}{{\texttt{VSS}}}
\newcommand{\EEV}{{\texttt{EEV}}}
\newcommand{\SAATP}{{\texttt{SAA-TP}}}
\usepackage[inline]{enumitem}
\usepackage{algorithmic}
\usepackage[ruled,linesnumbered]{algorithm2e}
\SetKw{Return}{Return:}
\SetKw{Returns}{returns}
\SetKw{Output}{Output:}
\SetKw{Input}{Input:}
\SetKw{Stepthree}{\emph{Step 3:}}
\SetKw{Steptwo}{\emph{Step 2:}}
\SetKw{Stepone}{\emph{Step 1:}}
\SetKwFor{While}{while}{do}{endw}
\SetKwIF{If}{ElseIf}{Else}{if}{then}{else if}{else}{endif}
\usepackage{xcolor}
\usepackage{multirow}
\newcolumntype{H}{>{\setbox0=\hbox\bgroup}c<{\egroup}@{}}
\newtheorem{result}{Result}
\usepackage{pgfplots}
\pgfplotsset{compat=1.18}

\usepackage{natbib}
 \bibpunct[, ]{(}{)}{,}{a}{}{,}%
 \def\bibfont{\small}%
\EquationsNumberedThrough    

\TheoremsNumberedThrough     
\ECRepeatTheorems  %

\begin{document}

\RUNAUTHOR{Maria Bazotte, Margarida Carvalho, and Thibaut Vidal}

\RUNTITLE{Solving Programs with Endogenous Uncertainty via Random Variable Transformation}

\TITLE{Solving Two-Stage Programs with Endogenous Uncertainty via Random Variable Transformation}

\ARTICLEAUTHORS{%

\AUTHOR{Maria Bazotte}
\AFF{MAGI \& CIRRELT, Polytechnique Montr\'{e}al, Montr\'{e}al, QC H3T 1J4, Canada, \EMAIL{maria-carolina.bazotte-corgozinho@polymtl.ca}}

\AUTHOR{Margarida Carvalho}
\AFF{DIRO \& CIRRELT, Universit\'{e} de Montr\'{e}al, Montr\'{e}al, QC H3T 1J4, Canada, \EMAIL{carvalho@iro.umontreal.ca}}

\AUTHOR{Thibaut Vidal}
\AFF{MAGI \& CIRRELT, Polytechnique Montr\'{e}al, Montr\'{e}al, QC H3T 1J4, Canada, \EMAIL{thibaut.vidal@polymtl.ca}}

} 

\ABSTRACT{%
Real-world decision-making problems often involve decision-dependent uncertainty, where the probability distribution of the random vector depends on the model's decisions. Few studies focus on two-stage stochastic programs with this type of endogenous uncertainty, and those that do lack general methodologies. We propose a general method for solving a class of these programs based on random variable transformation, a technique widely employed in probability and statistics. The random variable transformation converts a stochastic program with endogenous uncertainty (original program) into an equivalent stochastic program with decision-independent uncertainty (transformed program), for which solution procedures are well-studied. Additionally, endogenous uncertainty usually leads to nonlinear nonconvex programs, which are theoretically intractable. Nonetheless, we show that, for some classical endogenous distributions, the proposed method yields mixed-integer linear or convex programs with exogenous uncertainty. We validate this method by applying it to a network design and facility-protection problem, considering distinct decision-dependent distributions for the random variables. While the original formulation of this problem is nonlinear nonconvex for most endogenous distributions, the proposed method transforms it into mixed-integer linear programs with exogenous uncertainty. We solve these transformed programs with the sample average approximation method. We highlight the superior performance of our approach compared to solving the original program in the case a mixed-integer linear formulation of this program exists.
}


\maketitle


\section{Introduction}

Uncertainty is an inherent characteristic of most real-world decision-making problems~\citep{kuccukyavuz2017introduction}. Stochastic programming is a well-known framework in mathematical programming to model optimization problems with uncertain parameters; see, e.g.,~\citep{shapiro2021lectures} for a summary of the theoretical aspects and the most popular algorithms. Research in this area has focused mainly on \emph{exogenous or decision-independent uncertainty}~\citep{li2021review}, where the random vector is independent of the model decisions. In contrast, \emph{endogenous} or \emph{decision-dependent uncertainty}, where the random vector depends on the decisions, has received less attention, although it is prevalent in critical applications such as pre-disaster supply chain management~\citep{peeta2010pre,bhuiyan2020stochastic}, and production planning~\citep{ekin2018integrated}.

In the context of endogenous uncertainty, different types of decision dependency may occur~\citep{jonsbraaten1998class}. \citet{goel2006class} propose classifying decision dependency into two types: For problems with Type~1 endogenous uncertainty or decision-dependent probability, decisions impact the underlying probability distribution of the random parameters. Although little research focuses on this type of uncertainty, as noted by~\citet{li2021review}, it deals with various problems of practical relevance, such as the design of communication networks, power plants, and warehouses and resource allocation in queuing systems~\citep{pflug1990line}. With Type~2 endogenous uncertainty, decisions influence the time at which the realization of the random variables is totally or partially revealed. Thus, with Type~2 endogenous uncertainty, the scenario tree depends on the decisions~\citep{jonsbraaten1998class}. Problems with Type~2 uncertainty have received more attention than those with Type~1, being mostly multistage programs. \citet{li2021review} provide a concise review of methodologies and algorithms for multistage programs with Type~2 endogenous uncertainty. In the present work, we focus on Type~1 endogenous uncertainty, in particular, on \textit{two-stage stochastic mixed-integer linear programs with decision-dependent probability}~({\SMIPD}), in which the probability distribution depends on the first-stage decisions. In what follows, we refer to Type~1 endogenous uncertainty as endogenous uncertainty. 

The mathematical properties of a stochastic program, such as convexity, are strongly linked to whether its uncertainty is modeled as endogenous or exogenous. Consequently, solution procedures typically employed for stochastic programs with exogenous uncertainty become complex when dealing with endogenous uncertainty, necessitating \textit{a priori} a comprehensive analysis of its underlying properties~\citep{varaiya1988stochastic}. In fact, endogenous uncertainty usually leads to nonlinear nonconvex programs~\citep{dupacova2006optimization}, which are theoretically intractable (NP-hard). Moreover, few studies focus on general-purpose methodologies to solve programs with endogenous uncertainty sharing common properties (i.e., methodologies are typically problem-specific).

\paragraph{Contributions.} The main contribution of this work is a general method based on random variable transformation to address {\SMIPD} with discrete or continuous endogenous random variables. This method resolves two main issues in the literature: the lack of generalist methods and the difficulty of solving the usual nonconvex nonlinear programs with endogenous uncertainty. More specifically, the contributions of this work are as follows:
\begin{enumerate*}[label=(\roman*)]
    \item We define a transformation function to model the endogenous uncertainty by combining first-stage decisions with decision-independent random variables, thereby converting stochastic programs with endogenous uncertainty (i.e., original programs) into stochastic programs with exogenous uncertainty (i.e., transformed programs).
    \item We identify different techniques to define appropriate transformation functions for endogenous random variables following classical decision-dependent probability distributions, including a commonly studied case known as the ``discrete selection of distributions''.
    \item We show that, for these distributions, the obtained transformed programs have a convex continuous relaxation.
    \item We investigate the implications of endogenous uncertainty when applying sampling methods, particularly the sample average approximation ({\SAA}) method, to solve the transformed program.
    \item We apply the proposed approach to the network design and facility protection problem ({\NDFPP}), considering different decision-dependent distributions.
    \item Finally, for all endogenous distributions, we analyze the statistical estimators of the {\SAA} applied to the transformed program ({\SAATP}) of the {\NDFPP}.
\end{enumerate*}
In less than one hour on average, the {\SAATP} finds feasible solutions with significantly small over-estimated optimality gaps (between $-1\%$ and $1\%$) and a high confidence interval for all decision-dependent distributions. For the {\NDFPP} with most endogenous distributions, the deterministic equivalent of the original program~({\DEOP}) is nonlinear nonconvex. The exception is the discrete selection of distributions case, for which the {\DEOP} is a mixed-integer linear program~({\MILP}) obtained by auxiliary linearizations~\citep{bhuiyan2020stochastic}. In a considerably shorter time, the {\SAATP} obtains a better over-estimator of the optimality gap than the optimality gap obtained by solving this mixed-integer {\DEOP} (a difference of around $24\%$ on average). Solving the {\SAATP} also obtains better solutions than when using an average scenario, i.e., the expected value ({\EV}) problem.

\paragraph{Organization.} The remainder of this work is structured as follows: Section~\ref{sec:problem-defition} formalizes the class {\SMIPD} and discusses related works on endogenous uncertainty. Section~\ref{sec:methodology} presents our general methodology to obtain the transformed program and discusses the application of methods such as the {\SAA} to solve this resulting program. Section~\ref{sec:def-transf-function} catalogs techniques to define transformations to various endogenous distributions, and Section~\ref{sec:case-study} explains how to apply this approach to the {\NDFPP}. Finally, Section~\ref{sec:experimental-results} reports the experimental results, comparing the proposed method with other solution procedures, and Section~\ref{conclusion} concludes our analyses.

\section{Problem Definition and Related Works}
\label{sec:problem-defition}

We focus on two-stage stochastic mixed-integer linear programs with decision-dependent probability ({\SMIPD}). Let $x \in \mathbb{R}^{n_1}_{+}$ be the vector of first-stage variables (i.e., the ``here-and-now'' decisions) and $\xi(\omega) \in \mathbb{R}^{U}$ a random vector that depends on the random events $\omega \in \Omega$. For the sake of simplicity, we also denote $\boldsymbol{\xi}=\xi(\omega)$ and designate $\xi$ as a realization of this random vector. We use $\mathbb{E}_{\boldsymbol{\xi} \mid x}[\cdot]$ to denote the expected value operator with respect to $\boldsymbol{\xi}$ given the first-stage variables $x$. Finally, $y(\omega) \in \mathbb{R}^{n_2}_{+}$ is the vector of second-stage variables given a realization of the random events. The standard mathematical formulation for the {\SMIPD} is:
\begin{subequations}\label{Prog:basic}
\allowdisplaybreaks
\begin{align} 
\label{initial_1}\min \ \  & c^{\top}x + \mathbb{E}_{\boldsymbol{\xi} \mid x} \left[ Q\boldsymbol{(} x,\xi(\omega) \boldsymbol{)} \right]&\\
    \label{initial_2}\mbox{s.t.} \ \ &Ax \geq b \ , \ x \in \mathcal{X} & \\
    \label{prog:second-stage}& \textrm{where } Q\boldsymbol{(}x,\xi(\omega)\boldsymbol{)} = \min\{ q(\omega)^{\top}y(\omega): T(\omega)x + W(\omega)y(\omega) \geq h(\omega), y(\omega) \in \mathcal{Y} \}. &
\end{align}
\end{subequations}
In Program~\eqref{Prog:basic}, $\mathcal{X}\subseteq \mathbb{R}^{n_1}_{+}$ and $\mathcal{Y}\subseteq \mathbb{R}^{n_2}_{+}$ are two compact sets imposing non-negativity and continuous, integer, or binary restrictions. The first stage has the matrices $A$, $b$, and $c$ with appropriate dimensions. The second-stage Program~\eqref{prog:second-stage} is established after the realization of the random events $\omega \in \Omega$ and has the matrices $q(\omega)$, $h(\omega)$, $T(\omega)$, and~$W(\omega)$ with possible random elements and appropriate dimensions. Thus, the random vector $\xi^{T}(\omega)$ combines these random elements, i.e., $\xi^{T}(\omega) = \boldsymbol{(} q(\omega)^{\top},h(\omega)^{\top}, T(\omega), W(\omega) \boldsymbol{)}$, and its support (according to its unconditional distribution on $x$) is denoted by $\Xi \subseteq \mathbb{R}^U$. For the sake of clarity, we define the random vectors ${\boldsymbol{\xi}_x^{\top} = \boldsymbol{(} q_x(\omega)^{\top}, h_x(\omega)^{\top}, T_x (\omega), W_x(\omega) \boldsymbol{)}}$ for $x \in \mathcal{X}$, which follow the conditional probability distribution of the random vector $\boldsymbol{\xi}$ given $x$ and have support $\Xi_x \subseteq \Xi$. We also set ${g(x) = \mathbb{E}_{\boldsymbol{\xi} \mid x}[Q\boldsymbol{(}x, \xi(\omega)\boldsymbol{)}]}$.

The {\SMIPD} assumes that the probability distribution of $\boldsymbol{\xi}$ depends on the first-stage decisions $x \in \mathcal{X}$ (i.e., it computes the conditional expectation $\mathbb{E}_{\boldsymbol{\xi}|x}[\cdot]$ instead of the standard expectation $\mathbb{E}_{\boldsymbol{\xi}}[\cdot]$ used in the exogenous case). This endogenous probability influences the properties of the resulting formulation. For example, unlike the exogenous case, the convexity of the continuous relaxation of $Q\boldsymbol{(}x,\xi(\omega)\boldsymbol{)}$ in Program~\eqref{Prog:basic} does not necessarily yield the convexity of the continuous relaxation of $\mathbb{E}_{\boldsymbol{\xi}|x}[Q\boldsymbol{(}x,\xi(\omega)\boldsymbol{)}]$ 
\citep{varaiya1988stochastic,dupacova2006optimization}. As a result, algorithms commonly used for solving two-stage programs with exogenous uncertainty may not be suitable for {\SMIPD}.

This work assumes that the conditional probability distribution is known \emph{a priori} and defined by a distinct probability mass or density function for $x \in \mathcal{X}$ (i.e., $x$ is a parameter of the distribution). We justify this assumption based on its broad application in the literature, as reviewed in Section~\ref{subsec:literature-review}. We consider it a reasonable extension of the usual stochastic optimization assumption for programs with exogenous uncertainty (i.e., random vectors with known probability distributions). 

For the sake of simplicity, we next discuss our reasoning assuming that the random vector $\boldsymbol{\xi}$ is discrete and has finite support $\Xi$; the process is analogous for the continuous case. We write the {\DEOP} in its extensive form (i.e., the deterministic equivalent of Program~\eqref{Prog:basic}):
\begin{equation}\label{Prog:discrete}
    \min  \left\{ c^{\top}x  + \sum_{k \in K} p_k(x)q_k^{\top} y_k: \eqref{initial_2} \ , \ T_kx + W_k y_k \geq h_k \ \forall\ k \in K \ , \ y_k \in \mathcal{Y} \ \forall\ k \in K \right\},
\end{equation}
where $K$ is the set of realizations of $\boldsymbol{\xi}$ (i.e., the set of scenarios defined by the support $\Xi$). For each scenario $k \in K$, $q_k$, $h_k$, $T_k$, and~$W_k$ have appropriate dimensions, and $\xi_k^{\top} = (q_k^{\top},h_k^{\top},T_k,W_k)$. 
In addition, $p_k(x): \mathcal{X} \rightarrow [0,1]$ is a function returning the conditional probability of scenario $k \in K$ given the first-stage variables~$x$, in other words, $\boldsymbol{\xi}_x$ follows the probability distribution defined by $p_k(x)$ for every $k \in K$. We remark that scenarios in $\Xi \backslash \Xi_x$ have zero probability according to $p_k(x)$.

\subsection{Background of Decision-Dependent Probability Distribution} \label{subsec:background}
In Program~\eqref{Prog:discrete}, the probability function $p_{k}(x)$ multiplies the second-stage objective term $q_k^{\top}y_k$, leading to nonlinear terms. Handling this nonlinearity requires describing $p_{k}(x)$, which depends on the assumptions of $\boldsymbol{\xi}$. When the $U$ elements of the random vector $\boldsymbol{\xi}$ are independent, we have 
\begin{equation} \label{probability}
    p_k(x) = \prod_{u = 1}^{U}\hat{p}_{\hat{r}(u,k)}^u(x) \quad \forall\ k \in K,
\end{equation} 
where $\hat{r}(u,k)$ returns the realization in $\mathcal{R}_u=\{1,\ldots, R_u\}$ of a random element $u \in \mathcal{U}=\{1,\ldots, U\}$ in scenario $k \in K$. In addition, $\hat{p}_{\hat{r}(u,k)}^u(x): \mathcal{X} \rightarrow [0,1]$ is a function on $x$ returning the conditional probability of the realization $\hat{r}(u,k)$ of the random element $u \in \mathcal{U}$ in scenario $k \in K$. The function $p_k(x)$ can assume other configurations that do not necessarily use $\hat{p}_{\hat{r}(u,k)}^u(x)$ but instead return the conditional probability of the whole random vector $\boldsymbol{\xi}$.

A commonly studied case where decisions involve selecting a probability distribution from a predefined set is named ``discrete selection of distributions"~\citep{dupacova2006optimization,hellemo2018decision}. In this case, the decision-dependent probability is represented by a probability distribution for each independent random element $u \in \mathcal{U}$ from a set $\mathcal{D}_u= \{1, \ldots, D_u\}$ of available distributions, which can be either discrete or continuous. For discrete endogenous random variables, we have ${\hat{p}_{\hat{r}(u,k)}^u(x) = \sum_{d \in \mathcal{D}_u} \overline{p}_{\hat{r}(u,k)}^{ud} x_u^d}$ for every element $u \in \mathcal{U}$. Here, $x_{u}^d \in \{0,1\}$ is a subvector of $x \in \mathcal{X}$ representing the chosen distribution, such that $\sum_{d \in \mathcal{D}_u} x_u^d = 1$ for every $u \in \mathcal{U}$. Moreover, the parameter $\overline{p}_{\hat{r}(u,k)}^{ud}$ represents the probability of realization $\hat{r}(u,k)$ of random component $u$ when distribution $d$ is selected. If a distribution from the set $\mathcal{D}=\{1, \ldots, D\}$ is chosen for the entire random vector $\boldsymbol{\xi}$ (i.e., for scenario $k \in K$), then $p_k(x) = \sum_{d \in \mathcal{D}}\overline{p}^d_k x^d$, where $\overline{p}^d_k$ is the probability of scenario $k$ under distribution $d$. Similarly, $x^d \in \{0,1\}$ are first-stage variables representing the selected distribution, such that $\sum_{d \in \mathcal{D}}x^d = 1$. Additionally,~\citet{dupacova2006optimization} also categorizes problems involving the selection of parameters for a specific probability distribution as parameter selection.~\citet{hellemo2018decision} also introduce the probability-distortion case, where the chosen probability given $x$ is derived by combining a set of distributions. A general probability mass function $p_k(x)$ cannot be defined for these cases as it relies on additional assumptions.

\subsection{Literature Review}\label{subsec:literature-review}

Consistent with the assumptions of this work, the literature on decision-dependent probability focuses essentially on problems where the probability mass or density functions are known \emph{a priori}. In this section, we review only works that use this assumption and that consider two-stage programs. We split our review into works on discrete and continuous random variables.

\paragraph{Discrete random variables.}
For discrete random vectors, the elements are usually assumed to be independent. For the discrete selection of distributions, distinct applications of two-stage programs have been studied, such as network design protection \citep{viswanath2004investing,da2010stochastic,peeta2010pre,laumanns2014distribution}, facility protection \citep{medal2016allocating}, facility location-protection \citep{li2021stochastic}, fire risk management~\citep{bhuiyan2019stochastic}, supplier risk management~\citep{zhou2022stochastic}, facility protection and network design~\citep{bhuiyan2020stochastic}, and 
airline flight network expansion problems~\citep{csafak2022two}.
All these works used decision-dependent probabilities as in Equation~\eqref{probability}, except for~\citet{csafak2022two}, who chose the probabilities of the entire random vector. These applications usually lead to large-scale instances with numerous scenarios.

\citet{laumanns2014distribution} and \citet{bhuiyan2019stochastic} respectively obtained exact {\MILP} reformulations of the {\DEOP} by using a linear scaling procedure and the probability chain technique described by~\citet{o2013probability}. To solve large-scale problems,~\citet{laumanns2014distribution} proposed a bundle method, and~\citet{medal2016allocating,bhuiyan2020stochastic}, and~\citet{zhou2022stochastic} implemented the {\sf L}-shaped decomposition for solving reformulations resulting from the probability chain technique. Usually, $T_k$ is a zero matrix in these two-stage programs (i.e., the second stage can be solved \emph{a priori}). The exception is~\citet{zhou2022stochastic}, who used an exact McCormick envelope~\citep{mccormick1976computability} to obtain an MILP. These reformulations apply specifically to the case of the discrete selection of distributions and rely on the binary variables in $\hat{p}_{\hat{r}(u,k)}^u(x)$. However, they require numerous variables and big-M constraints. Thus, the practical use of these reformulations is limited to problems with a sufficiently small number of scenarios.~\citet{da2010stochastic} proposed a different approximate {\MILP} reformulation for the {\DEOP} by using the binary variables in $\hat{p}_{\hat{r}(u,k)}^u(x)$ to reformulate $p_k(x)$ as a set of exponential and linear constraints. The authors ensured a predefined precision by approximating the exponential constraints by piecewise linear cuts included in a branch-and-cut scheme.

Various decision-dependent probability mass functions have also been studied, such as the binomial distribution~\citep{karaesmen2004overbooking,green2013nursevendor}, linear interpolation~\citep{du2014stochastic}, a linear scaling and convex combination of a set of distributions~\citep{hellemo2018decision}, and the exponential function~\citep{bhuiyan2021models}. In most studies, the endogenous uncertainty in general cases of discrete random variables resulted in {\DEOP}s, which are mixed-integer nonlinear nonconvex programs ({\MINLP}s) for which proving global optimality is harder than for {\MILP}s or mixed-integer convex programs ({\MICP}s). Even if $p_k(x)$ are convex or concave~\citep[see, e.g.,][]{hellemo2018decision} or $\hat{p}_{\hat{r}(u,k)}^u(x)$ are convex or concave~\citep[see, e.g.,][]{du2014stochastic,bhuiyan2021models,holzmann2021shortest,homem2022simulation}, the term $p_k(x)q_k^{\top}y_k$ could lead to an {\MINLP}.
Within this context,~\citet{krasko2017two} studied a two-stage preparedness hazard problem for debris-flow-hazard mitigation and solved the {\DEOP} with the BARON solver.~\citet{hellemo2018decision} solved the {\DEOP} of an energy investment problem using a generalized Benders decomposition ({\GBD}) and the BARON solver. However, these approaches are limited to solving instances with few scenarios because the {\DEOP} in extensive form explicitly enumerates all scenarios.

Approximative approaches have been proposed for facility protection~\citep{medal2016allocating}, sensors for holding water quality~\citep{karaesmen2004overbooking}, and overbooking problems~\citep{bhuiyan2021models}. The objective function in these problems is submodular, allowing the design of greedy algorithms. However, approximation guarantees are generally low (e.g., around $63\%$ for~\citet{medal2016allocating}).~\citet{li2021stochastic} and~\citet{li2022locating} used genetic algorithms to solve two-stage programs of a facility location-protection problem with a bilevel program as the first stage. 

\citet{viswanath2004investing,peeta2010pre}, and~\citet{du2014stochastic} used the {\SAA} within iterative heuristic procedures to evaluate current first-stage solutions for two-stage network protection problems. Unfortunately, these procedures do not provide performance guarantees. The {\SAA} was also used to estimate the performance of current solutions within stochastic gradient methods for overbooking~\citep{karaesmen2004overbooking} and appointment-scheduling problems~\citep{homem2022simulation}.~\citet{homem2022simulation} proved the local convergence of the algorithm, which is compatible with its nonconvexity. Finally, by considering a specific gradient estimator,~\citet{karaesmen2004overbooking} proved the local convergence of the algorithm. 

\citet{da2010stochastic} used importance sampling to apply the {\SAA}. The author sampled from a biased probability distribution and adjusted the weight of the objective function according to $p_k(x)$ and the biased distribution. However, importance sampling depends on the choice of the biased distribution~\citep[Chap. 5]{shapiro2021lectures}, and finding a unique biased distribution appropriate for all conditional distributions on $x$ could be difficult.~\citet{holzmann2021shortest} proposed an adapted {\SAA} for bilevel shortest-path problems with randomized strategies. They sampled realizations from the continuous uniform distribution $U(0,1)$ and included big-M constraints in the model to transform these samples into realizations of the decision-dependent discrete distribution. 

\paragraph{Continuous random variables.}
To the best of our knowledge, only~\citet{ernst1995optimal,ekin2017augmented,ekin2018integrated}, and~\citet{hellemo2018decision} investigated continuous endogenous random variables.~\citet{ernst1995optimal} proposed a search procedure for solving a newsvendor problem with uncertain demand following a normal distribution with mean and variance selected by first-stage decisions. The augmented probability simulation ({\APS}), in which decisions are also simulated, was used by~\citet{ekin2017augmented} and~\citet{ekin2018integrated} for solving, respectively, a newsvendor problem with uncertain demand depending on stock level decisions and a production-planning problem with random yield depending on maintenance and production decisions. In both works, the dependency was modeled with a normal distribution with mean and variance defined according to the first-stage decisions. Finally,~\citet{hellemo2018decision} used the discretized function of the Kumaraswamy distribution and approximated the discretized function of the normal distribution. The discretization trades an increase in the number of scenarios for a confident discrete approximation. Approximations of the discretized function may also be needed when the cumulative distribution has no closed form or a complex closed form. In some applications, the search procedure cannot analyze the entire decision space, implying the need for a heuristic procedure. Finally, the {\APS} remains limited to problems where the second stage is represented by analytical expressions. 

\section{Methodology}\label{sec:methodology}

Our methodology is based on using random variable transformation, which is a classical topic in probability and statistics (see \citealt[Chap. 5]{hogg1977probability}), to model endogenous uncertainty. This section outlines Assumption~\ref{assump:transform}, which is necessary for applying our methodology. Section~\ref{sec:def-transf-function} shows that this assumption is not restrictive by evidencing many cases in which it holds. 
\begin{assumption}
There exists 
\begin{enumerate*}[label=(\arabic*)]
    \item a random vector $\vartheta^{\top}(\omega) = \boldsymbol{(} \widetilde{q}(\omega)^{\top},\widetilde{h}(\omega)^{\top},\widetilde{T}(\omega),\widetilde{W}(\omega)\boldsymbol{)}$ 
    with a known probability distribution that is independent of the first-stage decisions $x \in \mathcal{X}$ and
    \item a transformation function $t(x,\boldsymbol{\vartheta})$, such that $\boldsymbol{\xi}_x = t\left(x,\boldsymbol{\vartheta}\right)$.
\end{enumerate*}
\label{assump:transform}
\end{assumption}
Assumption~\ref{assump:transform} allows us to mimic the endogenous uncertainty by using a function $t(\cdot,\cdot)$ that combines the first-stage variables $x$ and the decision-independent random vector $\boldsymbol{\vartheta}$. For example, $t(\cdot,\cdot)$ can be the inverse transformation of $\boldsymbol{\xi}_x$~\citep[Chap. 2]{montecarlostatistical}. Here, the random vector $\boldsymbol{\vartheta}$ does not necessarily have the same dimension as $\boldsymbol{\xi}$ or $\boldsymbol{\xi}_x$. Furthermore, the random elements of $ \boldsymbol{\xi}_x$ may be dependent. For the sake of simplicity, we consider that $t( x, \boldsymbol{\vartheta}) = ( t_1 \boldsymbol{(} x,\widetilde{q}(\omega)^{\top} \boldsymbol{)}, t_2\boldsymbol{(} x,\widetilde{h}(\omega) \boldsymbol{)},t_3\boldsymbol{(} x,\widetilde{T}(\omega) \boldsymbol{)}, t_4\boldsymbol{(} x,\widetilde{W}(\omega) \boldsymbol{)} ) $, where $t_1 \boldsymbol{(} x,\widetilde{q}\left(\omega\right)^{\top} \boldsymbol{)}$, $t_2\boldsymbol{(}x,\widetilde{h}\left(\omega\right) \boldsymbol{)}$, $t_3\boldsymbol{(} x,\widetilde{T}\left(\omega\right) \boldsymbol{)}$, and $t_4\boldsymbol{(} x,\widetilde{W}\left(\omega\right) \boldsymbol{)}$ are transformation functions applied to the corresponding subvectors and matrices of $\boldsymbol{\vartheta}$. These functions return, respectively, the subvectors and matrices of $\boldsymbol{\xi}_x$, i.e., $q_{x}(\omega)^{\top}$, $h_{x}(\omega)$, $T_{x}(\omega)$, and $W_{x}(\omega)$. From Assumption~\ref{assump:transform}, $t(.)$ and thus $t_1(.)$, $t_2(.)$, $t_3(.)$, and $t_4(.)$ must account for any dependency between the random elements of $\boldsymbol{\xi}_x$. Section~\ref{sec:def-transf-function} details different methods to define appropriate transformations respecting Assumption~\ref{assump:transform}.

In the field of stochastic optimization,~\citet{varaiya1988stochastic} and~\citet{dupacova2006optimization} introduced the idea of transforming a decision-dependent probability distribution into an independent distribution in what is called ``plug-in'' or ``change of measure'' (see \citealt{pflug2012optimization} for a detailed explanation) and might employ random variable transformation. Despite these initial works, no study has directly modeled endogenous uncertainty as a function $t(.)$ nor investigated appropriate transformations for random variables following classical endogenous probability distributions. Specifically, the definition of $t(.)$ remains unexplored. Moreover, no study has focused on the nature and solution of problems resulting from this transformation nor, notably, addressed their potential nonconvexity. Finally, \citet{holzmann2021shortest} used random variable transformation for a specific endogenous discrete random variable, for which Section~\ref{subsubsec:discrete-rv} highlights the similarities and differences with the approach proposed herein. Therefore, the present work proposes a general framework using functions of random variables for addressing {\SMIPD} given Assumption~\ref{assump:transform}. We also investigate practical transformations for endogenous random variables following classical probability distributions, as well as the convexity of the resulting transformed program.

Based on Assumption~\ref{assump:transform}, we can formulate {\SMIPD} as the transformed program
\begin{equation}\label{Prog:transformed}
    \min \{ c^{\top}x + \mathbb{E}_{\boldsymbol{\vartheta}} \left[ Q_T\boldsymbol{(}x,\vartheta(\omega)\boldsymbol{)} \right]:  \eqref{initial_2} \},
\end{equation}
where the transformed second stage program $Q_T\boldsymbol{(}x,\vartheta(\omega)\boldsymbol{)}$ is defined as
\begin{equation}\label{Prog:2stagetransformed}
    \min  \{ t_1 \boldsymbol{(} x,\widetilde{q}\left(\omega\right)^{\top} \boldsymbol{)} y(\omega): t_3 \boldsymbol{(} x,\widetilde{T}\left(\omega\right) \boldsymbol{)}x + t_4  \boldsymbol{(} x,\widetilde{W}\left(\omega\right) \boldsymbol{)} y(\omega) \geq t_2 \boldsymbol{(} x,\widetilde{h}\left(\omega\right) \boldsymbol{)}, y(\omega) \in \mathcal{Y}  \},
\end{equation}
and we can outline the following lemma: 
\begin{lemma}\label{lemma1}
    Under Assumption~\ref{assump:transform}, 
    $g(x)=\mathbb{E}_{\boldsymbol{\xi} \mid x} \left[ Q\boldsymbol{(}x,\xi(\omega) \boldsymbol{)}\right] = \mathbb{E}_{\boldsymbol{\vartheta}} \left[ Q_T\boldsymbol{(}x,\vartheta(\omega)\boldsymbol{)} \right] $ 
    and the {\SMIPD} (Program~\eqref{Prog:basic}) can be reformulated into the transformed Program~\eqref{Prog:transformed}. In other words, Programs~\eqref{Prog:basic} and~\eqref{Prog:transformed} have the same optimal solutions. 
\end{lemma} 
Thus, our methodology is established on finding an appropriate transformation function to obtain Program~\eqref{Prog:2stagetransformed}. In that way, under Assumption~\ref{assump:transform}, {\SMIPD} can be transformed into a \emph{two-stage stochastic program with exogenous uncertainty}. Solution procedures for the transformed Program~\eqref{Prog:transformed} based on its properties are well-studied~\citep[see, e.g.,][]{shapiro2021lectures} and, thus, can be applied in this case. Program~\eqref{Prog:2stagetransformed} has the terms $t_1\boldsymbol{(}x,\widetilde{q}\left(\omega\right)^{\top} \boldsymbol{)} y$, $t_2\boldsymbol{(}x,\widetilde{h}\left(\omega\right) \boldsymbol{)}$, $t_3\boldsymbol{(} x,\widetilde{T}\left(\omega\right) \boldsymbol{)} x$, and $t_4\boldsymbol{(} x,\widetilde{W}\left(\omega\right) \boldsymbol{)} y$, which might be nonlinear nonconvex. We then introduce the following assumption:
\begin{assumption}\label{assump:convexity}
The continuous relaxation of the second-stage Program~\eqref{Prog:2stagetransformed} is convex.
\end{assumption}
Program~\eqref{Prog:2stagetransformed} might need to be reformulated to convexify its continuous relaxation. Section~\ref{sec:def-transf-function} presents different cases for which Assumption~\ref{assump:convexity} holds. In addition, the exogenous random vector~$\boldsymbol{\vartheta}$ might have a large or infinite dimension (see, e.g., the inversion in Section~\ref{subsec:inverse-transf}). Hence, next, we describe the application of the {\SAA}~\citep[Chap. 5]{shapiro2021lectures} to solve Program~\eqref{Prog:transformed}, which, although intuitive, shows that the transformation does not compromise desirable {\SAA} properties. Nonetheless, other procedures can also be employed to solve Program~\eqref{Prog:transformed} in this case, e.g., stochastic gradient~\citep{shapiro2021lectures}.

\subsection{Sample Average Approximation} \label{subsec:saa}

We consider $N$ independent and identically distributed samples $( \boldsymbol{\vartheta^1}, \ldots, \boldsymbol{\vartheta^N})$ of $\boldsymbol{\vartheta}$ with realizations $( \vartheta^1, \ldots, \vartheta^N )$. The transformed {\SAA} program is
\begin{equation}\label{Prog:saa}
  \min \{ c^{\top}x + \hat{g}_N(x) : \ \   \eqref{initial_2}\},
\end{equation}
where $\hat{g}_N(x)$ is the {\SAA} estimator of $\mathbb{E}_{\boldsymbol{\vartheta}} \left[ Q_T \boldsymbol{(}x, \vartheta(\omega) \boldsymbol{)} \right]$ defined as
\begin{equation}
    \label{SAA_1} \hat{g}_N(x) =  \frac{1}{N} \sum_{i = 1}^{N} Q_T(x,\vartheta^i).  
\end{equation}
This estimator converges point-wise with probability one (w.p.1) to $\mathbb{E}_{\boldsymbol{\vartheta}} [Q_T\boldsymbol{(}x,\vartheta(\omega)\boldsymbol{)}]$ when ${N\rightarrow +\infty}$. From Lemma~\ref{lemma1}, we have $\mathbb{E}[ \hat{g}_N(x) ] = g(x)$, so $\hat{g}_N(x)$ is an unbiased estimator of $g(x)$. We denote $v^{*}$, $v_{t(x,\boldsymbol{\vartheta})}^{*}$ and $v_N$ as the optimal values of Program~\eqref{Prog:basic}, Program~\eqref{Prog:transformed}, and Program~\eqref{Prog:saa}, respectively. Under some mild conditions, such as when Program~\eqref{Prog:2stagetransformed} is convex, an {\MILP} or when the first-stage set $\mathcal{X}$ is finite, $v_N$ converges to $v_{t(x,\boldsymbol{\vartheta})}^{*}$ w.p.1 as ${N \rightarrow +\infty}$ (cf. \citealt[Chap. 5]{shapiro2021lectures}, and \citealt{kleywegt2002sample}, respectively). In these cases and under Assumption~\ref{assump:transform}, $v_N$ converges to $v^{*}$ w.p.1 as $N \rightarrow +\infty$, since $v^{*}=v_{t(x,\boldsymbol{\vartheta})}^{*}$ from Lemma~\ref{lemma1}.   This result, and its rate of convergence (cf.~\citep{shapiro2021lectures}), depends on the properties of $\boldsymbol{\vartheta}$ (i.e., on the convergence of $v_N$ to $v_{t(x,\boldsymbol{\vartheta)}}$), which may have larger dimension and support than $\boldsymbol{\xi}$ and may require a larger number of samples.  

The {\SAA} finds a feasible first-stage solution $\overline{x} \in \mathcal{X}$. For that, we solve $M$ instances of the {\SAA} Program~\eqref{Prog:saa}, each with $N$ independent samples. We then evaluate the optimal solutions of every instance by using a larger sample size $N'$ and define $\overline{x}$ as the optimal solution with the best evaluation. When $Q\boldsymbol{(}x, \xi(\omega)\boldsymbol{)}$ is finite and Program~\eqref{Prog:basic} has relatively complete recourse (i.e., the second stage associated with $x \in \mathcal{X}$ is feasible for almost every realization of $\boldsymbol{\xi}$), then $Q_T\boldsymbol{(}x, \vartheta(\omega)\boldsymbol{)}$ is also finite, and Program~\eqref{Prog:transformed} also has relatively complete recourse, given that $\boldsymbol{\xi}_x = t(x,\boldsymbol{\vartheta})$ from Assumption~\ref{assump:transform}. In this case, in addition to finding $\overline{x}$, we can compute a lower-bound estimator $\overline{v}_M^N$ of $v^{*}$, an upper-bound estimator $\hat{v}_{N'}(\overline{x})$ of $v^{*}$, an over-estimator $\overline{gap}(\overline{x})$ of the optimality gap $c^{\top}\overline{x} + g(\overline{x})-v^{*}$, and unbiased estimators of the variances of $\overline{v}_N^M$, $\hat{v}_{N{'}}(\overline{x})$, and $\overline{gap}(\overline{x})$, i.e., $\hat{\sigma}_{N,M}^2$, $\hat{\sigma}_{N{'}}^2(\overline{x})$, and $\hat{\sigma}_{gap}^2(\overline{x})$, respectively. The Appendix~\ref{app:saa-estimators} details these estimators and the {\SAA}'s procedure.

By defining the $\boldsymbol{\xi}_x$ with the function $t(x,\boldsymbol{\vartheta})$, the estimators $\hat{g}_N(x)$ for all $x \in \mathcal{X}$ are based on the same samples of $\boldsymbol{\vartheta}$. Choosing a transformation $t(x,\boldsymbol{\vartheta})$ that creates a positive correlation between these estimators $\hat{g}_N(x)$ is desirable because it reduces the variance of their difference. Indeed, when comparing two first-stage decisions $x_1$ and $x_2$, the variance of the difference $\hat{g}_{N}(x_1) - \hat{g}_{N}(x_2)$ when both estimators use the same random samples of $\boldsymbol{\vartheta}$ is $Var\boldsymbol{(}\hat{g}_{N}(x_1)\boldsymbol{)} + Var\boldsymbol{(}\hat{g}_{N}(x_2)\boldsymbol{)} - 2Cov\boldsymbol{(}\hat{g}_{N}(x_1),\hat{g}_{N}(x_2)\boldsymbol{)}$, where $Var(\cdot)$ and $Cov(\cdot,\cdot)$ denote the variance and covariance. This is the reasoning behind the well-known efficiency-improvement technique of common random numbers~\citep{glasserman1992some}. This discussion extends to the application of any sampling technique for solving the transformed program. Section~\ref{sec:def-transf-function} further discusses this question.

\section{Transformation Functions: Definition and Examples}\label{sec:def-transf-function}

This section presents appropriate transformations $\boldsymbol{\xi}_x = t(x,\boldsymbol{\vartheta})$ with $x \in \mathcal{X}$ for classical distributions in which both Assumptions~\ref{assump:transform}~and~\ref{assump:convexity} hold. We classify the transformation approaches according to four main categories to organize our discussion: inversion (Section~\ref{subsec:inverse-transf}), enumeration (Section~\ref{subsec:enumeration}), standardization (Section~\ref{subsec:standardization}), and convolution (Section~\ref{subsec:convolution}). Then, Section~\ref{subsec:adv-disadv} discusses the proposed methodology's advantages and disadvantages. We use the notation introduced in Section~\ref{subsec:background}, and we denote $\boldsymbol{\xi}_x^u$ as the $u$-th element of the random vector $\boldsymbol{\xi}_x$, and $\boldsymbol{\vartheta}_u$ as the random variable or vector used to obtain $\boldsymbol{\xi}_x^u$. Similarly, $\xi_x^{u}$ and $\vartheta_u$ denote their realizations. Unless stated otherwise, we consider the random elements of $\boldsymbol{\xi}_x$ to be independent. If the random elements of $\boldsymbol{\xi}_x$ are independent, then the $\boldsymbol{\vartheta}_u$ for $u \in \mathcal{U}$ must also be independent. Thus, we can model each random element $\boldsymbol{\xi}_x^u$ separately by using its marginal probability distribution conditioned on $x$ and the $\boldsymbol{\vartheta}_u$. Nonetheless, in certain cases, when the elements of $\boldsymbol{\xi}_x$ are correlated, we can have correlated $\boldsymbol{\vartheta}_u$ and also model each $\boldsymbol{\xi}_x^u$ separately in the same way as in the independent case, as explained further in this section. As a result, we mainly focus on the individual transformation of an element $u$ in this section, i.e., $\boldsymbol{\xi}_x^u = t(x,\boldsymbol{\vartheta}_u)$. We remark that most studies have focused on the independent case (see Section~\eqref{subsec:literature-review}), while our approach is suitable for both independent and dependent cases.

\subsection{Inversion}\label{subsec:inverse-transf}

The inversion method is commonly used to sample a random variable when the probability mass or density function is known (e.g., the inverse sampling method, cf. \citealt[Chap. 2]{montecarlostatistical}). This technique transforms a random variable following the continuous uniform distribution $U(0,1)$ into the desired random variable. We define $\mathbb{F}_{\boldsymbol{\xi}^u}(x,\beta)(\cdot)$ to be the conditional marginal cumulative distribution function of the random variable $\boldsymbol{\xi}^u$ given $x$, the first-stage variables, and $\beta$, a set of remaining parameters that are not decisions of the optimization model. To apply the inverse transformation, we compute the generalized inverse of the conditional cumulative distribution:
\begin{equation}\label{eq:generalized-inverse}
    \mathbb{F}^{-1}_{\boldsymbol{\xi}^u}(x,\beta)(j) = \min \left\{ \xi \in \Xi_x \mid \mathbb{F}_{\boldsymbol{\xi}^u}(x,\beta)(\xi) \geq j \right\}.  
\end{equation}
The transformation is then defined as $\boldsymbol{\xi}_x^u = t(x,\boldsymbol{\vartheta}_u) =\mathbb{F}^{-1}_{\boldsymbol{\xi}^u}(x,\beta)(\boldsymbol{\vartheta}_u)$ (Property 2 of Assumption~\ref{assump:transform}), such that $\boldsymbol{\vartheta}_u$ follows the $U(0,1)$ and is independent of $x$ (Property 1 of Assumption~\ref{assump:transform}). The inversion method, in many cases, has the desirable property of inducing a positive correlation (see Section~\ref{subsec:saa}) between the estimators of different systems (i.e., the $\hat{g}_N(x)$ for $x \in \mathcal{X}$) when their samples are generated from the same sequence of uniform random numbers~\citep{whitt1976bivariate}.

If the elements of $\boldsymbol{\xi}_x$ are dependent and their correlation matrix is known, then we can define correlated random elements $\boldsymbol{\vartheta}_u$ for $u \in \mathcal{U}$ using an appropriate copula. A copula is a multivariate cumulative distribution function that ensures the marginal distribution of each random variable $\boldsymbol{\vartheta}_u$ is the $U(0,1)$. Copulas are well-known for modeling dependencies~\citep[see, e.g.,][]{nelsen2006introduction}. Thus, as in the independent case, each $\boldsymbol{\vartheta}_u$ can be transformed into $\boldsymbol{\xi}_x^u$ separately using the marginal distribution of this endogenous random element (i.e., the transformation is equal for both independent and dependent cases; the difference lies in the definition of $\boldsymbol{\vartheta}_u$). To keep the paper self-contained, we illustrate the application of the inverse in Appendix~\ref{app:examples} (Example~\ref{ex:inverse}) for both cases. We note that this method is limited to cases where the function $\mathbb{F}^{-1}_{\boldsymbol{\xi}^u}(x,\beta)(j)$ has a well-defined or approximated form. Next, we investigate the application of this technique for discrete random variables and continuous random variables according to the properties of their cumulative probability distribution. 

\subsubsection{Discrete Random Variables.}\label{subsubsec:discrete-rv}

In this case, the inverse transformation $\boldsymbol{\xi}_x^u = \mathbb{F}^{-1}_{\xi^u}(x,\beta)(\boldsymbol{\vartheta}_u)$ is simplified and consists of defining the realization $r \in \mathcal{R}_u$ such that
\begin{equation} \label{discrete_transformation}
   \sum_{j=1}^{r-1} \hat{p}_j^u(x) \leq \vartheta_u < \sum_{j=1}^{r} \hat{p}_j^u(x),
\end{equation}
where the sum on the left-hand-side equals zero when $r=1$, and $\hat{p}_{j}^{u}(x)$ and $\sum_{j=1}^{r} \hat{p}_{j}^{u}(x)$ are the marginal and cumulative probability functions (see Section~\ref{subsec:background}) of the realization $ r$ of random element $u \in \mathcal{U}$. Thus, we define the variables $\theta_{u}^r \in \{0,1\}$, which take the value one if the realization $r \in \mathcal{R}_u$ respects Inequalities~\eqref{discrete_transformation}. We also have the variables $\pi_{u}^r \in \{0,1\}$, which equal the sum $\sum_{j=1}^{r}\theta_{u}^{j}$ for every realization $r \in \mathcal{R}_u$. The following constraints model Inequalities~\eqref{discrete_transformation}:
\begin{align}
    \label{general_drv_1}& \pi_{u}^{r} \leq \frac{1}{\vartheta_u}\sum_{j=1}^{r} \hat{p}_{j}^u(x) - \epsilon && \forall\ r \in \{1,\ldots,R_u-1\} \\
    \label{general_drv_2}& (1-\pi_{u}^{r-1}) \leq \frac{1}{(1-\vartheta_u)}\sum_{j=r}^{R_u} \hat{p}_{j}^u(x)  && \forall\ r \in \{2,\ldots,R_u\} \\
    \label{general_drv_3}& \pi_{u}^r = \pi_{u}^{r+1} - \theta_{u}^{r+1} && \forall\ r \in \{1,\ldots,R_u-1\} \\
    \label{general_drv_4}& \sum_{r \in \mathcal{R}_u} \theta_{u}^r = 1, &&
\end{align}
where $\epsilon>0$ is a sufficiently small value used to define the strict inequality in~\eqref{discrete_transformation}. Constraints~\eqref{general_drv_1} and~\eqref{general_drv_2} ensure that the right- and left-hand side in Inequalities~\eqref{discrete_transformation} are respected, respectively. Constraints~\eqref{general_drv_3} establish the connection between $\pi^r_{u}$ and $\theta^r_{u}$ so that $\theta^r_{u}=1$ for the right realization, where $\pi_{u}^{R_u}$ is set to one. Although we could model Constraints~\eqref{general_drv_1} and~\eqref{general_drv_2} with variables $\theta_u^r$, we use the variables $\pi_u^{r}$ to reduce the number of nonzero elements in these constraints. Constraint~\eqref{general_drv_4} restricts the choice of a unique realization for element $u$. Finally, $\nu_r^u$ is the value of realization $r$ of the random element $u$, so $t(x,\vartheta_u) = \sum_{r=1}^{R_u} \nu_r^u \theta_u^r$. To ensure Assumption~\ref{assump:convexity}, the continuous relaxation of Constraints~\eqref{general_drv_1} and~\eqref{general_drv_2} must be convex. This leads to the following theorem:
\begin{theorem}\label{theorem_general_discrete}
    If all the functions $\hat{p}^u_{r}(x)$ for $u \in \mathcal{U}$ and $r \in \mathcal{R}_u$ are either convex or concave over $\left[0,1\right]^{n_1}$, then Constraints~\eqref{general_drv_1} and~\eqref{general_drv_2} have a convex continuous relaxation.
\end{theorem}
The Appendix~\ref{app:proofs} presents the proof of Theorem~\ref{theorem_general_discrete}. Different studies investigate problems with convex or concave $p_r^u(x)$ (see Section~\ref{subsec:literature-review}). Despite this desirable property, the original program is usually nonconvex. In contrast, using the inverse transformation function results in only linear or convex constraints. In addition, even if the continuous relaxation of the $\hat{p}_{r}^u(x)$ is neither concave nor convex, we can convexify Constraints~\eqref{general_drv_1} and \eqref{general_drv_2} in certain cases (this analysis is problem-specific). Nonetheless, in addition to relying on the probability function $\hat{p}^u_{r}(x)$ to ensure convexity, the number of integer variables and constraints increases polynomially with the number of random parameters $U$ and their realizations $R_u$, which impacts the problem's complexity and highlights the computational challenges of handling endogenous uncertainty.

To the best of our knowledge, only \citet{holzmann2021shortest} used inversion to model a specific discrete endogenous random variable in {\SMIPD}. In contrast, we describe a general framework adapted to various functions $\hat{p}_{r}^u(x)$. The set of Constraints~\eqref{general_drv_1}--\eqref{general_drv_4} also has a tighter continuous relaxation feasibility set than that employed by~\citet{holzmann2021shortest} (see Appendix~\ref{app:comparison-tranf-discrete}).

In certain cases, such as when $\boldsymbol{\xi}_x^u$ has only two possible realizations and $p_r^u(x)$ is a linear function, simpler models can be designed. This is the case of the Bernoulli distribution, described next.

\paragraph{Bernoulli Distribution.} We consider that $\boldsymbol{\xi}_x^u \sim Bern(\phi)$, with the probability of success $\phi \in [0,1]$ being a subvector of $x$. The inverse is such that $\xi_x^{u}=1$ if $\vartheta_u < \phi$, and $\xi_x^{u}=0$ otherwise, with $\boldsymbol{\vartheta}_u \sim U(0,1)$. We model this with the variables $\upsilon_u \in \{0,1\}$, i.e., $t(x,\vartheta_u) = \upsilon_u$, and the constraints
\begin{equation}\label{const:bernoulli}
    \{\vartheta_u \upsilon_u \leq \phi - \epsilon \; ; \; \upsilon_u \geq \phi - \vartheta_u \},
\end{equation}
where $\epsilon>0$ is again a sufficiently small value. Constraints~\eqref{const:bernoulli} guarantee that realization $0$ and $1$ are chosen if $\vartheta_u$ is superior to $\phi$ and if $\vartheta_u$ is inferior or equal to $\phi$, respectively.

For both transformations (general and specific to Bernoulli), if only $h(\omega)$ is uncertain, then Assumption~\ref{assump:convexity} is respected. Otherwise, linearization techniques, such as McCormick, can be applied to linearize the bilinear terms in $t_1\boldsymbol{(}x,\widetilde{q}(\omega)^{\top}\boldsymbol{)} y$, $t_3\boldsymbol{(}x,\widetilde{T}(\omega)\boldsymbol{)} x$, and $t_4\boldsymbol{(}x,\widetilde{W}(\omega)\boldsymbol{)} y$, so Assumption~\ref{assump:convexity} holds. In this case, the problem's complexity increases further, and the transformation is restricted by the dimensions of variables $x$ and $y$, due to the inclusion of extra variables and constraints when applying linearization techniques.

\subsubsection{Continuous Random Variables.}

In general, the transformation $\xi_x^u = \mathbb{F}^{-1}_{\xi^u}(x,\beta)(\vartheta_u)$ cannot be simplified for continuous random variables. Nevertheless, when, for fixed $x$ and $\beta$, the conditional cumulative distribution is strictly increasing and continuous, the generalized inverse is redundant because $\mathbb{F}^{-1}_{\boldsymbol{\xi}^u}(x,\beta)(\cdot)$ is the unique inverse function. This property is commonly respected for well-behaved distributions and leads to a simpler transformation. Two applications of the inverse for some well-known continuous distributions respecting this property are listed below.

\paragraph{Uniform distribution.} If $\boldsymbol{\xi}^u \sim U(a,b)$ and both parameters $a \in \mathbb{R}$ and $b \in \mathbb{R}$ are first-stage variables, then we can use the transformation function $t(a,b,\boldsymbol{\vartheta}_u) = a + (b-a) \boldsymbol{\vartheta}_u$.

\paragraph{Exponential distribution.} If $\boldsymbol{\xi}^u \sim \exp(\lambda)$ with $\lambda > 0$ a first-stage decision, then the transformation is defined as $t(\lambda,\boldsymbol{\vartheta}_u) = -{\ln(1-\boldsymbol{\vartheta}_u)}/{\lambda}$. This function is concave for $\lambda > 0$. 

Both distributions are part of the special case where the unique inverse $\mathbb{F}^{-1}_{\boldsymbol{\xi}^u}(x,\beta)(\cdot)$ is convex or concave on $x \in \mathcal{X}$, i.e., the inverse function is convex or concave on the parameters of the distribution. Other distributions that fall into this special case are the Cauchy and logistic distributions. Additionally, there are distributions for which $\mathbb{F}^{-1}_{\boldsymbol{\xi}^u}(x,\beta)(\cdot)$ has no closed-form expression, but good numerical approximations are available, e.g. student, chi-square, and normal distributions. Nevertheless, transformations based on other techniques, which lead to simpler formulations, can be used for certain distributions (see Section~\ref{subsec:standardization} for the normal distribution). The analysis of the special class of distributions introduced before leads to the following theorem:
\begin{theorem}\label{theorem:continuous-dist-conv}Assumption~\ref{assump:convexity} holds if the following criteria are met:
    \begin{enumerate*} [label=(\arabic*)]
        \item the cost vector $q$, technology matrix $T$ and recourse matrix $W$ are fixed;
        \label{item1}
        \item $\mathbb{F}_{\boldsymbol{\xi}^u}(x,\beta)(\cdot)$ for fixed $x$ and $\beta$ is continuous and strictly increasing; and \label{item2}
        \item the transformation $t(x,\boldsymbol{\vartheta}_u) = \mathbb{F}^{-1}_{\boldsymbol{\xi}^u}(x,\beta)(\cdot)$ is a convex function on $x \in \mathcal{X}$, or concave if the recourse constraint in Program~\eqref{Prog:2stagetransformed} is an inferior or equal to inequality.\label{item3}
        \end{enumerate*}
\end{theorem}
The Appendix~\ref{app:proofs} presents the proof of Theorem~\ref{theorem:continuous-dist-conv}. If Criterion~\ref{item1} does not hold, the products $t_1\boldsymbol{(}x,\widetilde{q}(\omega)^{\top}\boldsymbol{)} y$, $t_3\boldsymbol{(}x,\widetilde{T}(\omega)\boldsymbol{)} x$ and $t_4\boldsymbol{(}x,\widetilde{W}(\omega)\boldsymbol{)} y$ may be nonconvex even if Criterion~\ref{item3} holds. Thus, a problem-specific convexity analysis is required to confirm if Assumption~\ref{assump:convexity} is satisfied. This evidences the complexity added by endogenous uncertainty, as the transformation can, in some cases, produce a nonconvex nonlinear transformed program with exogenous uncertainty, which is theoretically intractable.

\subsection{Enumeration} \label{subsec:enumeration}

For a random element $\boldsymbol{\xi}_x^u$ with $u \in \mathcal{U}$, a simple transformation can be employed if the set of alternative marginal distributions $\mathcal{D}_u$ of element $u$ is sufficiently small to allow enumeration. This corresponds to the ``discrete selection of distributions" (see Section~\ref{subsec:background}). We have $ \boldsymbol{\vartheta}_u = \{ \boldsymbol{\vartheta}_{ud} \}_{d \in \mathcal{D}_u}$, where $\boldsymbol{\vartheta}_{ud}$ is a random variable following distribution $d \in \mathcal{D}_u$. 
The $\boldsymbol{\vartheta}_{u}$ is independent of $x$ (Property~1 of Assumption~\ref{assump:transform}). We thus have ${t(x,\boldsymbol{\vartheta}_u)=\sum_{d = 1}^{D_u}\boldsymbol{\vartheta}_{ud}x_u^d}$  (Property~2 of Assumption~\ref{assump:transform}), with $x_u^d$ as defined in Section~\ref{subsec:background}. This transformation can be applied for $\boldsymbol{\xi}_x$ with independent or correlated random elements by also using copulas. To keep the paper self-contained, we illustrate this in Appendix~\ref{app:examples}, Example~\ref{ex:discrete-selection}, and also present the case where the distribution is selected for the entire random vector $\boldsymbol{\xi}_x$ in Example~\ref{ex:discrete-selection-all}. If only $h(\omega)$ is uncertain, then Assumption~\ref{assump:convexity} is directly satisfied. Otherwise, bilinear terms may appear due to $t_1\boldsymbol{(}x,\widetilde{q}(\omega)^{\top}\boldsymbol{)}y$, $t_3\boldsymbol{(}x,\widetilde{T}(\omega)\boldsymbol{)}x$ and $t_4\boldsymbol{(}x,\widetilde{W}(\omega)\boldsymbol{)}y$. In such cases, auxiliary linearizations such as the McCormick envelope can be applied so Assumption~\ref{assump:convexity} is respected. However, this approach is limited not only by the size of $\mathcal{D}_u$ but also by the dimensions of $y$ or $x$, as it increases the number of variables, leading to computational challenges.

\subsection{Standardization} \label{subsec:standardization}

A location-scale family of distributions has distributions with location $\mu$ and scale $\sigma$ parameters~\citep[Chap. 3.5]{casella2024statistical}. The distributions of this family can be derived from a standard family distribution member with location and scale parameters equal to one and zero, respectively. If $\boldsymbol{\xi}^u$ follows a location-scale family distribution with parameters $\mu$ and $\sigma$ as first-stage decisions, then $t(\mu,\sigma,\boldsymbol{\vartheta}_u) = \mu +\sigma \boldsymbol{\vartheta}_u$ is valid (Property 2 of Assumption~\ref{assump:transform}). The $\boldsymbol{\vartheta}_u$ follows the standard family distribution member, which is independent of the values of $\mu$ and $\sigma$ (Property 1 of Assumption~\ref{assump:transform}). If only the vector $h(\omega)$ has random elements, then Assumption~\ref{assump:convexity} holds. Otherwise, a specific analysis of the terms $t_1(x,\widetilde{q}(\omega)^{\top}) y$, $t_3(x,\widetilde{T}(\omega)) x$ and $t_4(x,\widetilde{W}(\omega)) y$ should be conducted. Thus, similar to the inversion, the transformed program may also be nonlinear and nonconvex, resulting in computational challenges. Many distributions belong to the location-scale family, including normal, generalized extreme value, exponential, and logistic distributions. Next, we exemplify this method for two of these distributions.

\paragraph{Normal distribution.} If $\boldsymbol{\xi}^u \sim N(\mu,\sigma^2)$ and $\mu$ and $\sigma$ are first-stage decisions, then we can apply the proposed transformation with $\boldsymbol{\vartheta}_u$ following the standard normal $N(0,1)$.

\paragraph{Generalized extreme value distribution.} If $\boldsymbol{\xi}^u \sim GEV(\mu,\sigma,\beta)$ with $\mu$ and $\sigma$ first-stage decisions and $\beta$ a fixed parameter, then we can apply the proposed transformation with $\boldsymbol{\vartheta}_u \sim GEV(0,1,\beta)$.

\subsection{Convolution}\label{subsec:convolution}

Various random variables that follow well-known distributions can be expressed as the sum of other independent random variables. Notable examples include the binomial distribution, which is the sum of Bernoulli random variables, and the hypoexponential distribution, which is the sum of exponential random variables. Suppose that $\boldsymbol{\xi}_x^u$ is equal to the sum $\boldsymbol{\zeta}_x^{u1} + \boldsymbol{\zeta}_x^{u2} + \ldots + \boldsymbol{\zeta}_x^{un}$ of independent random variables $\boldsymbol{\zeta}_x^{uj}$'s following a specified endogenous distribution. In this case, the probability distribution of $\boldsymbol{\xi}_x^u$ is the convolution of the marginal distributions of the $\boldsymbol{\zeta}_x^{uj}$'s~\citep[Chap. 2]{hogg2013introduction}. We can define a transformation $\boldsymbol{\zeta}_x^{uj} = \hat{t}_j(x,\boldsymbol{\vartheta}_{uj})$ for every endogenous random variable $\boldsymbol{\zeta}_x^{uj}$, $j \in \{1,\ldots,n\}$, independently. Thus, $t(x,\boldsymbol{\vartheta}_u) = \sum_{j=1}^n \hat{t}_j(x,\boldsymbol{\vartheta}_{uj})$, with $\boldsymbol{\vartheta}_u = \{\boldsymbol{\vartheta}_{uj}\}_{j=1}^n$. The complexity of the problem depends on the definition of each transformation $\hat{t}_j(x,\boldsymbol{\vartheta}_{uj})$ and the number of components $n$, thus limiting the application of this method by the size of $n$. Next, we demonstrate this transformation for the binomial distribution. We also provide an example of dependent random variables for which this methodology can be applied in Appendix~\ref{app:examples} (Example~\ref{ex:correlated}).

\paragraph{Binomial distribution.} We consider that $\boldsymbol{\xi}_x^{u} \sim Bin(n,\phi)$, where the probability of success of one trial $\phi \in [0,1]$ is a subvector of $x$, and the number of trials $n$ is a parameter. We have $\boldsymbol{\xi}_x^u = \boldsymbol{\zeta}_x^{u1} + \boldsymbol{\zeta}_x^{u2} + \ldots + \boldsymbol{\zeta}_x^{un}$, where the $\boldsymbol{\zeta}_x^{uj}$'s are endogenous and follow the $Bern(\phi)$ independently. We define the $\boldsymbol{\zeta}_x^{uj}$'s with the inversion for the Bernoulli distribution of Section~\ref{subsubsec:discrete-rv}. We have ${\boldsymbol{\vartheta}_{uj} \sim U(0,1)}$ independently for $j \in \{1,\ldots,n\}$. We set the variables $\upsilon_{ur} \in \{0,1\}$ equal to one if realization $\vartheta_{ur} < \phi$ and zero otherwise, and include Constraints~\eqref{const:bernoulli} for all $\boldsymbol{\zeta}_x^{uj}$'s. Thus, $t(x,\vartheta_u) = \sum_{j=1}^{n} \upsilon_{uj}$. If only $h(\omega)$ is uncertain, then Assumption~\ref{assump:convexity} is respected. Otherwise, the bilinear terms in $t_1(x,\widetilde{q}(\omega)^{\top})y$, $t_3(x,\widetilde{T}(\omega))x$, and $t_4(x,\widetilde{W}(\omega))y$ can be linearized by using the McCormick envelope (Assumption~\ref{assump:convexity} holds). In this case, the application of this technique is limited to the dimensions of $y$ and $x$ due to the increase in the number of variables. Furthermore, inversion (see Section~\ref{subsubsec:discrete-rv}) may be used in this case because the binomial probability of every realization is a convex function on the probability of success $\phi$ (i.e., we can apply Theorem~\ref{theorem_general_discrete}). Both approaches have a similar number of variables and constraints and include big-M constraints. However, inversion leads to convex (potentially nonlinear) constraints, whereas convolution leads to linear constraints.

\subsection{Advantages and Disadvantages}\label{subsec:adv-disadv}

The transformed Program~\eqref{Prog:transformed}, a two-stage program with exogenous uncertainty, offers several advantages:
\begin{enumerate*}[label=(\roman*)]
    \item Solution procedures for such programs are extensively researched and can be applied based on their specific properties.
    \item Large-scale decomposition techniques, such as the L-shaped algorithm, are applicable depending on the characteristics of the transformed second-stage Program~\eqref{Prog:2stagetransformed}.
    \item Sampling techniques become straightforward, enabling the treatment of problems with continuous endogenous random variables or discrete ones with many realizations. In contrast, using sampling to solve the original program with endogenous uncertainty is challenging, as the distribution depends on first-stage decisions.
\end{enumerate*}
Besides, our method is versatile as various techniques can be employed to define the transformation function $t(x,\boldsymbol{\vartheta})$. As a result, multiple transformations may exist to model the same endogenous random vector. For instance, both the inverse and convolution methods can be applied to model the binomial random variable. The choice of a transformation should aim to achieve a transformed program with desirable properties. Furthermore, as demonstrated earlier in this section, in many cases, a transformation function leading to a mixed-integer linear or convex transformed Program~\eqref{Prog:transformed} exists, even if the original Program~\eqref{Prog:basic} is nonlinear nonconvex. This is particularly the case for problems where only the vector $h(\omega)$ is uncertain, as the transformation $t_2\boldsymbol{(}x,\widetilde{h}\left(\omega\right) \boldsymbol{)}$ does not multiply the variables $x$ or $y$ in Program~\eqref{Prog:2stagetransformed}. For the discrete selection of distributions where only $h(\omega)$ is uncertain, the original program is an {\MILP} obtained by linearizing the objective function with many variables and big-M constraints (see~\citep{bhuiyan2020stochastic}). However, we propose a transformation based on the enumeration of distributions that reduces the need for additional variables or constraints in this case.

Nevertheless, certain transformations, such as enumeration and the inverse transformation for discrete random variables, may introduce in the formulation many constraints and variables, which are often integer when dealing with discrete distributions. Hence, we explore the impact of such transformations in our experimental results. For problems where the cost vector $q(\omega)$, technology matrix $T(\omega)$, and/or recourse matrix $W(\omega)$ are uncertain, obtaining a transformed program with a convex relaxation is more challenging due to the terms $ t_1(x,\widetilde{q}(\omega)^{\top})y $, $ t_3(x,\widetilde{T}(\omega))x $, and $ t_4(x,\widetilde{W}(\omega))y$. For discrete distributions, introducing integer variables to define the transformations enables the use of McCormick linearization as a last resort. However, this approach is constrained by the dimensions of $x$, $y$, and the number of uncertain parameters. Additionally, the dimension and support of the exogenous random vector $\boldsymbol{\vartheta}$ are typically larger than those of the endogenous random vector $\boldsymbol{\xi}_x$. The convergence properties of sampling techniques depend on the support of $\boldsymbol{\vartheta}$, potentially requiring more samples. However, leveraging the same random samples of $\boldsymbol{\vartheta}$ for the {\SAA} estimators $\hat{g}_N(x)$ for all $x \in \mathcal{X}$ to induce correlation between them is consistent with the efficiency-improvement technique of common random numbers (see discussion in Section~\ref{subsec:saa}).

\section{Case Study}\label{sec:case-study}

Pre-disaster or pre-interdiction problems with endogenous uncertainty have been investigated in the context of facility protection~\citep{medal2016allocating,bhuiyan2020stochastic,li2021stochastic,li2022locating} and network design protection~\citep{viswanath2004investing,da2010stochastic,peeta2010pre,du2014stochastic,laumanns2014distribution}. To the best of our knowledge, these applications are the most recurrent in previous research. Thus, this section showcases our methodology by applying it to the {\NDFPP}, which was introduced by~\citet{bhuiyan2020stochastic}.~\citet{bhuiyan2020stochastic} defined the endogenous random variables through a discrete selection of distributions. In contrast, we explore this and three other decision-dependent probability mass or density functions for the random variables. This section introduces our transformed program for the {\NDFPP}, which has different versions according to the definition of the endogenous distribution. 

\subsection{Network Design and Facility Protection Problem }

The {\NDFPP} is a {\SMIPD} with complete recourse. First-stage decisions determine facility protection allocations and the network structure. Next, the scenario (natural disruptions) is realized, affecting facility capacity levels, and second-stage decisions establish an optimal flow to satisfy client demands. Facility capacity levels are mutually independent endogenous random variables. 

We consider a set $\mathcal{D}_0$ of natural disruption events, representing different types of natural disasters centered at specific locations. Facilities experience varying disruption intensities based on their distance from the disruption center $d \in \mathcal{D}_0$, with intensity levels in the set $\mathcal{L}=\{1,\ldots,L\}$. We define the set of possible events $\mathcal{D}=\mathcal{D}_0 \cup \{nd\}$, which includes the no-natural-disruption event $nd$. The network is represented by the undirected graph $G=(\mathcal{N}_\mathcal{T},\mathcal{E})$, where $\mathcal{N}_\mathcal{T}=\mathcal{N_F} \cup \mathcal{N_C} \cup \{m\}$ consists of facilities ($\mathcal{N_F}$), clients ($\mathcal{N_C}$), and a dummy facility ($m$). The set $\mathcal{N}_\mathcal{T}^0=\mathcal{N}_\mathcal{T}\backslash\{m\}$ excludes the dummy facility. We also specify a directed graph $G_o= (\mathcal{N}_\mathcal{T},\mathcal{A})$, including arcs $(i,j)$ and $(j,i)$ for each edge $[i,j] \in \mathcal{E}$. The edge corresponding to arc $(i,j) \in \mathcal{A}$ is denoted as $e(i,j)$. To represent unsatisfied demand, we include arcs $(m,n)$ for every client node $n \in \mathcal{N_C}$ in $\mathcal{A}$, with the auxiliary set of arcs arriving at the dummy facility defined as $\mathcal{A}_m=\{(m,n)\ \forall\ n \in \mathcal{N_C}\}$. The set $\mathcal{P}=\{1,\ldots,P\}$ contains available protection levels for each facility, where larger values indicate better protection. Finally, $\mathcal{S}$ is the set of scenarios, each composed of an event and a realization of the exogenous random vector. The function $\hat{d}(s)$ returns the event $d \in \mathcal{D}$ for scenario $s \in \mathcal{S}$. The realization of the exogenous random vector in scenario $s$ is denoted $\vartheta^s$. The definition of this random vector depends on the endogenous distribution and transformation function, which we define later. Moreover, protecting a facility $f \in \mathcal{N_F}$ with protection level $ p \in \mathcal{P} $ costs $ c_f^p $, and opening an edge $ e \in \mathcal{E} $ costs $ c_e $. Each client node $ n \in \mathcal{N_C} $ has a demand $ b_n $. The budget for edge opening and facility protection is denoted by $C$. The per-unit-flow transportation cost for arc $ (i,j) \in \mathcal{A} $ is $ q_{ij} $. Moreover, the cost for arcs $ \mathcal{A}_m $ leaving the dummy facility is $ a \cdot \max_{(i,j) \in \mathcal{A}\backslash \mathcal{A}_m } q_{ij} $, where $ a $ is the penalty per unit of unmet demand. The big-M value for edge capacity is $ B = \sum_{n \in \mathcal{N_C}} b_n $.
For the first stage, we define the binary variables $x_f^p$, which are equal to one if protection level $ p \in \mathcal{P} $ is allocated to facility $f \in \mathcal{N_F}$ and zero otherwise. We also include the binary variables $z_e$, which are equal to one if edge $e \in \mathcal{E}$ is installed and zero otherwise. For the second stage, we specify the variables $ y_{ij}^s $, which represent the flow through arc $ (i,j) \in \mathcal{A} $ in scenario $ s \in \mathcal{S} $. The endogenous capacity of facilities, $\xi_x^{fs}$ for facility $f \in \mathcal{N_F}$ in scenario $s \in \mathcal{S}$, depends on the $x$ variables. We summarize this notation in Appendix~\ref{app:detailed-notation}. The transformed program for the {\NDFPP} is
\begin{subequations}\label{Prog:NDFPP-firststage}
 \begin{align}
    \min \mbox{: } \label{constr:objctive} & \mathbb{E}_{\boldsymbol{\vartheta}}[Q_T(x,\boldsymbol{\vartheta})] &&\\
    \mbox{s.t. } \label{constr:facility-protec} & \sum_{p \in \mathcal{P}} x_f^p = 1 && \quad \forall\ f \in \mathcal{N_F} \\
    \label{constr:budget} & \sum_{f \in \mathcal{N_F}}\sum_{p \in \mathcal{P}}c_f^p x_f^p + \sum_{e \in \mathcal{E}} c_e z_e \leq C && \\
    \label{constr:x-domain}& x_f^p \in \{0,1\}&& \quad \forall\ f \in \mathcal{N_F},\quad p \in \mathcal{P} \\
    \label{constr:z-domain}& z_e \in \{0,1\}&& \quad \forall\ e \in \mathcal{E}, 
\end{align}
\end{subequations}
where the transformed second-stage program $Q_T(x,\vartheta_s)$ for realization $\vartheta_s$ of $\boldsymbol{\vartheta}$ is
\begin{subequations}\label{Prog:NDFPP-secondstage}
 \begin{align}
    \label{constr:objctive-secondstage} \min \mbox{: } &\sum_{(i,j) \in \mathcal{A}}q_{ij}y_{ij}^{s}  && \\
    \label{constr:flow-conservation} & \sum_{j \in \mathcal{N_T}} y_{jn}^s - \sum_{j \in \mathcal{N_T}} y_{nj}^s = b_n && \quad \forall\ n \in \mathcal{N_C} \\
    \label{constr:facility-capacity} & \sum_{j \in \mathcal{N_T}} y_{fj}^s - \sum_{j \in \mathcal{N_T}} y_{jf}^s \leq \xi_x^{fs} && \quad \forall\ f \in \mathcal{N_F} \\
    \label{constr:arc-capacity} & y_{ij}^s \leq B z_{e(i,j)}  && \quad \forall\ (i,j) \in \mathcal{A} \backslash \mathcal{A}_m \\
    \label{constr:facility-transformation} & \xi_x^{fs}=t(x,\vartheta^s) && \quad \forall\ f \in \mathcal{N_F} \\
    \label{constr:y-domain} & y_{ij}^s \geq 0 && \quad \forall\ (i,j) \in \mathcal{A}.
\end{align}
\end{subequations}
In the first stage, the objective function~\eqref{constr:objctive} minimizes the expected value of the transformed second-stage program. We also have Constraints~\eqref{constr:facility-protec} and \eqref{constr:budget}, which specify the selection of one level of protection for each facility and ensure that the cost budget for investing in facility protection and edge installation, respectively, is respected. The objective function \eqref{constr:objctive-secondstage} of the second stage minimizes the average cost of transporting goods from facilities to clients. In addition, the second stage contains Constraints~\eqref{constr:flow-conservation}--\eqref{constr:facility-transformation}. These ensure flow conservation at client nodes (Constraints~\eqref{constr:flow-conservation}), the use only of available capacity (Constraints~\eqref{constr:facility-capacity}), the existence of flow only in installed arcs (Constraints \eqref{constr:arc-capacity}), and a suitable definition of the endogenous realizations $\xi_x^{fs}$ (Constraints~\eqref{constr:facility-transformation}). We note that, following the notation in Program~\eqref{Prog:basic}, only $h(\omega)$ is uncertain in this problem. We show the transformed {\SAA} Program~\eqref{Prog:saa} of the {\NDFPP} in Appendix~\ref{app:transformed-saa}.

In the following sections, we describe the various problem variants based on the decision-dependent distribution of $\xi_x^{fs}$. Specifically, we detail Constraints~\eqref{constr:facility-transformation} and define $\vartheta^s$. For discrete endogenous distributions, we assume that the post-disruption capacity of facilities can take on $|\mathcal{W}|$ possible values, where $\mathcal{W}=\{0,\ldots, W\}$ is an ordered index set. Here, level 0 represents no capacity, and level $W$ represents maximum capacity. The parameter $\nu_w$ for $w \in \mathcal{W}$ denotes the capacity at level $w$, with $\nu_{w+1} = \nu_{w} + \overline{\nu}$, where $\overline{\nu}$ is the increase in capacity per unit increase in capacity level.
 
\subsubsection{Discrete Selection of Distributions.}\label{subsec:discretechoice-problem}

First, we consider the discrete selection of distributions introduced by \citet{bhuiyan2020stochastic}. If a protection level $p \in \mathcal{P}$ is installed at facility $f \in \mathcal{N_F}$, then its capacity level is considered to follow $Bin(W,\overline{\phi}_{f}^{dp})$ when event $d \in \mathcal{D}$ is realized. The parameter $\overline{\phi}_{f}^{dp}$ is the probability of successfully increasing the capacity level by one. We follow the enumration in Section~\ref{subsec:enumeration}. For facility $f \in \mathcal{N_F}$, protection $p \in \mathcal{P}$, and scenario $s \in \mathcal{S}$, we define $\vartheta^s_{fp}$ as the value of the realization of the capacity level following the $Bin(W,\overline{\phi}_{f}^{dp})$. In this case, no extra variables are needed. We have $\xi_{x}^{fs} = \sum_{p \in \mathcal{P}} \vartheta^s_{fp} x_f^p$ for every facility $f \in \mathcal{N_F}$, which model Constraints~\eqref{constr:facility-transformation}. Finally, this variant of the considered problem is denoted {\NDFPP}-Selection.

\subsubsection{Binomial Distribution.} \label{subsec:binomial-problem}

We assume that the network-wide protection allocation affects the capacity level distribution of all facilities. Thus, we aim to represent situations in which facilities collaborate by sharing resources, situations in which facilities intercommunicate with transmission lines (e.g., power systems), or situations where failures propagate in the network (e.g., the spread of diseases). We consider that the capacity level of a facility $f \in \mathcal{N_F}$ follows $Bin(W,{\phi}_{f}^d)$ when event $d \in \mathcal{D}$ occurs, where ${\phi}_{f}^d \in [0,1]$ is a variable computed as:
\begin{equation}\label{eq:binomial-probability}
    {\phi}_{f}^d = \frac{1}{\rho} \sum_{i \in \mathcal{N_F}} \sum_{p \in \mathcal{P}} \widetilde{\phi}_{fi}^{dp} x_{i}^p \ \quad \forall\ d \in \mathcal{D}, \quad f \in \mathcal{N_F},
\end{equation}
where the parameter $\widetilde{\phi}_{fi}^{dp}$ is the impact on the capacity level at facility $f \in \mathcal{N_F}$ of installing protection level $p \in \mathcal{P}$ at facility  $i \in \mathcal{N_F}$ when event $d$ occurs, and $\rho$ is a normalization parameter to guarantee that ${\phi}_{f}^d \in [0,1]$. Hence, we employ the convolution technique for the Binomial distribution in Section~\ref{subsec:convolution}; that is, $\vartheta_{fw}^s$ is the realization following the $U(0,1)$ for scenario $s \in \mathcal{S}$, facility $f \in \mathcal{N_F}$, and capacity level $w \in \{1,\ldots,W\}$. The variables $\upsilon^s_{fw} \in \{0,1\}$ are equal to one if $\vartheta_{fw}^s < {\phi}^{\hat{d}(s)}_f$ and zero otherwise for $f \in \mathcal{N_F}$ and $w \in \{1,\ldots,W\}$. Finally, we adapt Constraints~\eqref{const:bernoulli} to guarantee the adequate definition of variables $\upsilon^s_{fw}$ (see Appendix~\ref{app:linear-reformu-stdnorma}). Thus, we determine the capacity value in every scenario as $\xi_x^{fs} = \overline{\nu} \left( \sum_{w=1}^{W} \upsilon_{fw}^s  \right)$ for every facility $f \in \mathcal{N_F}$, which are equivalent to Constraints~\eqref{constr:facility-transformation}.  
This variant of the considered problem is denoted {\NDFPP}-Binomial. Note that the {\DEOP} of this variant is nonlinear nonconvex (see Appendix~\ref{app:deop-nonconvex}).

\subsubsection{Discrete Distribution.} \label{subsec:stdnorma-problem}

We continue to assume that network-wide protection influences the capacity levels of all facilities. We consider that the random variables follow a discrete distribution inspired by the multinomial logit models in discrete choice theory~\citep{Bierlaire1998}. Therefore, the probability $\hat{p}_{df}^{w}(u)$ of a facility $f \in \mathcal{N_F}$ having capacity level $w \in \mathcal{W}$ after the occurrence of an event $d \in \mathcal{D}$ is defined as
\begin{equation}\label{prob_std_norma}
    \hat{p}_{df}^{w}(u) = \frac{ u_{df}^w } { \sum_{j \in \mathcal{W}} u_{df}^j } \ \quad \forall\ d \in \mathcal{D},\quad f \in \mathcal{N_F},\quad w \in \mathcal{W},
\end{equation} 
where $u_{df}^w \geq 0$ are auxiliary variables expressing the likelihood of facility $f \in \mathcal{N_F}$ having capacity level $w \in W$ when event $d \in \mathcal{D}$ occurs. These variables are 
\begin{equation}\label{eq:stdnorma-utility}
    u_{df}^w = \sum_{i \in \mathcal{N_F}} \sum_{p \in \mathcal{P}} \widetilde{u}_{fi}^{dpw} x_{i}^p \ \quad \forall\ d \in \mathcal{D},\quad f \in \mathcal{N_F},\quad w \in \mathcal{W}, 
\end{equation}
where the parameter $\widetilde{u}_{fi}^{dpw}$ is the impact of installing protection $p \in \mathcal{P}$ at facility $i \in \mathcal{N}_\mathcal{F}$ on the likelihood of capacity level $w \in \mathcal{W}$ at facility $f \in \mathcal{N_F}$ when event $d \in \mathcal{D}$ occurs. Thus, we adapt the inverse transformation in Section~\ref{subsubsec:discrete-rv}. For each facility $f \in \mathcal{N_F}$ and scenario $s \in \mathcal{S}$, the realization $\vartheta_f^s$ follows the $U(0,1)$. The variables $\theta_{fw}^{s} \in \{0,1\}$ equal one if the capacity level $w \in \mathcal{W}$ is selected for facility $f \in \mathcal{N_F}$ in scenario $s \in \mathcal{S}$ and zero otherwise. We also define the variables $\pi_{fw}^s \in \{0,1\}$, which are equal to $\sum_{j=0}^{w} \theta_{fj}^s$, with $\pi_{fW}^s$ set to one. We adapt the Constraints~\eqref{general_drv_1}--\eqref{general_drv_4} for the {\NDFPP} to ensure the correct inverse transformation of $\vartheta_f^s$ according to the probability distribution defined by Equation~\eqref{prob_std_norma}. The obtained constraints initially contain bilinear terms involving binary variables, and thus, to satisfy Assumption~\ref{assump:convexity}, we apply the exact McCormick envelope. This adaptation and linearization procedures are detailed in Appendix~\ref{app:linear-reformu-stdnorma}. As a result, we have $\xi_x^{fs} = \sum_{w \in \mathcal{W}} \nu_w \theta_{fw}^s$ for every facility $f \in \mathcal{N_F}$, which model Constraints~\eqref{constr:facility-transformation}. This variant is denoted {\NDFPP}-Discrete. Finally, the {\DEOP} in this case is also nonlinear nonconvex (see Appendix~\ref{app:deop-nonconvex}).

\subsubsection{Normal Distribution.}\label{subsec:normal-problem}

We also assume that network-wide protection affects the facilities' capacity level. We consider that the capacity level of facility $f \in \mathcal{N_F}$ follows $N(\mu_f^d,{(\widetilde{\sigma}_f^d)}^2)$ when event $d \in \mathcal{D}$ occurs. The standard deviation $\widetilde{\sigma}_f^d \geq 0$ is a parameter and the mean $\mu_f^d \geq 0$ is a decision variable such that
\begin{equation}
    \label{constr:normal-mean} \mu_f^d = \frac{1}{\rho}  \sum_{i \in \mathcal{N_F}} \sum_{p \in \mathcal{P}} \widetilde{\mu}_{fi}^{dp}  x_{i}^p \ \quad \forall\ d \in \mathcal{D},\quad f \in \mathcal{N_F}.
\end{equation}
The parameter $\widetilde{\mu}_{fi}^{dp}$ is the impact of installing protection $p \in \mathcal{P}$ at facility $i \in \mathcal{N}_\mathcal{F}$ on the mean of the capacity level of facility $f \in \mathcal{N}_\mathcal{F}$ when event $d \in \mathcal{D}$ occurs, and $\rho$ is a normalization value. We ensure that the capacity is always positive, i.e., $\mu_f^d \geq 4 \widetilde{\sigma}_f^d$ for all feasible decisions. We use the standardization method in Section~\ref{subsec:standardization}, so the realization $\vartheta_f^s$ follows $N(0,1)$ for each facility $f \in \mathcal{N_F}$ and scenario $s \in \mathcal{S}$. Constraints~\eqref{constr:facility-transformation} are modeled as $\xi^{fs}_x = \mu_f^{\hat{d}(s)}  + \vartheta_f^s  \sigma_f^{\hat{d}(s)}$ for every $f \in \mathcal{N_F}$. This variant is denoted {\NDFPP}-Normal, and its {\DEOP} is nonlinear nonconvex (see Appendix~\ref{app:deop-nonconvex}).

\subsection{Experimental Results}\label{sec:experimental-results}

This section validates the application of our methodology for addressing {\SMIPD}. First, Section~\ref{subsec:instances-generation} summarizes the instance generation. Section~\ref{subsec:results-comparison-baseline} then compares the {\SAA} applied to the transformed program ({\SAATP}) and the {\DEOP} for {\NDFPP}-Selection, specifically, the {\MILP} reformulation obtained by auxiliary linearizations in~\citet{bhuiyan2020stochastic}. The {\DEOP} for the {\NDFPP}-Binomial, {\NDFPP}-Discrete, and {\NDFPP}-Normal is nonlinear nonconvex (see Section~\ref{sec:case-study}). Thus, Section~\ref{subsec:results-estimators} only explores the performance of the estimators obtained by the {\SAATP}. Finally, we investigate the value of the stochastic solution ({\VSS}) of the {\NDFPP} in Section~\ref{subsec:results-EEV}. Note that different methods can be used to solve our obtained two-stage program with exogenous uncertainty. We apply the {\SAA}, given its simplicity and desirable properties (see Section~\ref{subsec:saa}). 

The instances of the transformed {\SAA} program, the {\DEOP}, the {\EV} problem, and the expected result of the {\EV} solution ({\EEV}) for computing the {\VSS} were coded in C++ and solved directly with Gurobi 10.0.3. All experiments were processed on the Beluga cluster of the Compute Canada network on an Intel Gold 6148 Skylake @ 2.4 GHz with a maximum of 100~GB. We limited Gurobi to using one thread. The details on the application of the {\SAA} are given in Appendix~\ref{app:saa-estimators}. All instances, results, and code are publicly available on GitHub~(\url{https://github.com/mariabazotte/saa_endog}\label{fn:github}).

\subsubsection{Instance Generation.}\label{subsec:instances-generation}

We follow a procedure similar to that of~\citet{bhuiyan2020stochastic} for network, transportation costs, and client demand generation, which is detailed in Appendix~\ref{app:instances} and available on the same GitHub project. We consider five networks, with $|\mathcal{N}_\mathcal{T}^0|=15$, 25, 30, 39, and 48, and assume the values 4 and 5 for $| \mathcal{N}_{\mathcal{F}}|$. The protection costs $c_{f}^{p}$, which are randomly generated, and the edge-opening costs $c_e$ are defined similarly as in~\citet{bhuiyan2020stochastic} (see Appendix~\ref{app:instances}). The budget $C$ is a percentage of the total cost with maximum protection, $C= 0.5 ( \sum_{n \in \mathcal{N_F}} c_n^{P} + \sum_{e \in \mathcal{E}} c_e )$, making it adjustable to different instance sizes. We consider four protection levels ($P=4$), and define the maximum facility capacity as $\nu_{\text{max}}=| { \sum_{n \in \mathcal{N_C}} b_n }/{0.9 |\mathcal{N_F}| } |$. Thus, the capacity increase $\overline{\nu}$ is ${\nu_{\text{max}}}/{W}$ with capacity values $\nu_w=w \overline{\nu}$ for $w \in \mathcal{W}$. We use $W = 2, 3,$ and $4$. As in~\citet{bhuiyan2020stochastic}, we consider three disruption events with distinct centers of occurrence ($|\mathcal{D}_0|=3$), and we assume three disruption intensity levels ($L=3$). We use the procedure in~\citet{bhuiyan2020stochastic} to define the intensity level $\hat{l}(f,d)$ in $\mathcal{L}$ at which facility $f \in \mathcal{N_F}$ is impacted by disruption $d \in \mathcal{D}_0$ (see Appendix~\ref{app:instances}).

Analogously to~\citet{bhuiyan2020stochastic}, we define the probability of successfully increasing the capacity level by one, $\overline{\phi}_{f}^{dp}$, as $({0.95 \times c_f^p}/{c_f^P} )^{{\hat{l}(f,d)}/{L}}$, with a 5\% risk of protection malfunction for each capacity level increase. For the $nd$ event, $\overline{\phi}_{f}^{nd,p}$ is defined as $\eta ( {0.95 \times c_f^p}/{c_f^P} ) + (1-\eta)$, where $\eta=0.7$ reflects the influence of protection decisions on capacity level uncertainty on this event. Unlike \citet{bhuiyan2020stochastic}, we assume that facilities are affected by the $nd$ event and low-intensity disruption events, which results in instances with a larger number of scenarios. For the {\NDFPP}-Binomial and {\NDFPP}-Discrete, we set $\widetilde{\phi}_{fi}^{dp}$ and $\widetilde{u}_{fi}^{dpw}$ to $\overline{\phi}_{i}^{dp}$ if $f = i$, and to $0.3 \times \overline{\phi}_{i}^{dp}$ otherwise. For the {\NDFPP}-Binomial and {\NDFPP}-Normal, $\rho$ is defined as $1 + 0.3 \times (|\mathcal{N_F}| - 1)$. For the {\NDFPP}-Normal, $\widetilde{\mu}_{fi}^{dp}$ is $\overline{\mu}_f^{dp}$ if $f = i$, and $0.3 \times \overline{\mu}_f^{dp}$ otherwise. The parameter $\overline{\mu}_f^{dp}$ is given by $\nu_{\text{max}} \left[1 - \frac{(P + 1 - p)(L + \hat{l}(f,d) + 2)}{2(P+1)(L+1)}\right]$ for each protection level $p \in \mathcal{P}$, with $\hat{l}(f,d)$ returning zero for the $nd$ event, meaning that the mean increases with protection and decreases with intensity. The standard deviation $\widetilde{\sigma}^d_f$ is $\frac{1}{\rho} \left[\overline{\sigma}_f^d + 0.3 \sum_{i \in \mathcal{N}_\mathcal{F} \backslash \{f\}} \overline{\sigma}_i^d\right]$, where $\overline{\sigma}_f^d$ is $\nu_{\text{max}}\left[\frac{L + \hat{l}(f,d) + 2}{8(P+1)(L+1)}\right]$. This standard deviation increases with intensity levels and is independent of protection, ensuring $4\overline{\sigma}_f^d \leq \overline{\mu}_{f}^{d1} \leq \cdots \leq \overline{\mu}_{f}^{dP}$.

Given that the costs $c_f^p$ are randomly generated, we consider five different seeds per instance configuration, (network, number of facilities, and capacity level). Thus, for {\NDFPP}-Selection, {\NDFPP}-Binomial, and {\NDFPP}-Discrete, we have $5 \times 2 \times 3=30$ instances per seed (150 in total). For {\NDFPP}-Normal, we have $5 \times 2=10$ instances per seed (50 in total).

\subsubsection{Comparative results for {\NDFPP}-Selection.}\label{subsec:results-comparison-baseline}

For the case of discrete selection of distributions, the literature on endogenous uncertainty has mainly linearized the probability function~\eqref{probability} with different approaches. Here, we consider the {\DEOP} formulated in~\citet{bhuiyan2020stochastic} as an {\MILP}. To solve this {\MILP}, we implemented the same methodology as \citet{bhuiyan2020stochastic}:  an {\sf L}-shaped decomposition accelerated with nondominated cuts and valid inequalities. We compare this approach with the {\SAATP}. Note that we cannot apply the {\SAA} to the {\DEOP} as it has endogenous uncertainty.
\begin{table}[]
\TABLE
{Comparison of the {\SAATP} and the {\DEOP} for {\NDFPP}-Selection instances. \label{tab:linea-vs-saa}}
{\resizebox{\textwidth}{!}{%
\begin{tabular}{rr|rr|rrr|rr|rr|rrr|rr}
\hline \up\down
&  & \multicolumn{7}{c|}{Four-facility instances}                                                                                                                          & \multicolumn{7}{c}{Five-facility instances}                                                                                                                             \\ \cline{3-16} \up\down
                              &                       & \multicolumn{2}{c|}{{\DEOP}}           & \multicolumn{3}{c|}{{\SAATP}}                                       & \multicolumn{2}{c|}{Comparison}                         & \multicolumn{2}{c|}{{\DEOP}}            & \multicolumn{3}{c|}{{\SAATP}}                                        & \multicolumn{2}{c}{Comparison}                          \\ \cline{3-16} \up\down
$|\mathcal{N}_\mathcal{T}^0|$ & $W$                   & $|K|$ & \multicolumn{1}{r|}{$GAP(\%)$}  & $|K_t|$ & $\overline{gap}(\overline{x})(\%)$ & \multicolumn{1}{r|}{Time(s)}        & $GAP_{UB}$(\%) & $GAP_{LB}$(\%) & $|K|$ & \multicolumn{1}{r|}{$GAP$(\%)}  & $|K_t|$ & $\overline{gap}(\overline{x})(\%)$ & \multicolumn{1}{r|}{Time (s)}         & $GAP_{UB}$(\%) & $GAP_{LB}$(\%) \\ \hline \up
15                            & 2   & 324   & \multicolumn{1}{r|}{$23.5\pm0.3$} & $4\cdot3^{16}$ & $0.4\pm0.2$                & \multicolumn{1}{r|}{$1086.4\pm947.5$} & $0.2\pm0.2$                       & $-29.9\pm0.5$           & 972   & \multicolumn{1}{r|}{$16.3\pm0.5$} & $4\cdot 3^{20}$  & $0.2\pm0.1$                & \multicolumn{1}{r|}{$1540.8\pm279.6$} & $0.5\pm0.4$                       & $-18.8\pm0.4$           \\
                              & 3   & 1024  & \multicolumn{1}{r|}{$13.1\pm0.3$} & $4\cdot 4^{16}$ & $0.3\pm0.1$                & \multicolumn{1}{r|}{$1135.6\pm937.6$} & $0.0\pm0.0$                       & $-14.8\pm0.3$           & 4096  & \multicolumn{1}{r|}{$14.4\pm0.2$}  & $4\cdot 4^{20}$ & $0.2\pm0.1$                & \multicolumn{1}{r|}{$982.8\pm98.7$}   & $0.3\pm0.2$                       & $-16.3\pm0.3$           \\ \down
                              & 4   & 2500  & \multicolumn{1}{r|}{$10.8\pm0.4$} & $4\cdot 5^{16}$ & $0.1\pm0.1$                & \multicolumn{1}{r|}{$952.8\pm503.9$}  & $0.0\pm0.1$                       & $-11.9\pm0.4$           & 12500 & \multicolumn{1}{r|}{$10.9\pm7.3$} & $4\cdot 5^{20}$  & $0.6\pm0.1$                & \multicolumn{1}{r|}{$861.8\pm107.5$}  & $3.9\pm7.9$                       & $-7.3\pm0.1$            \\ \hline \up
25                            & 2   & 324   & \multicolumn{1}{r|}{$23.0\pm3.4$} & $4\cdot 3^{16}$ &$ 0.6\pm0.2                $& \multicolumn{1}{r|}{$300.4\pm93.6$}   &$ 0.1\pm0.1                       $&$ -29.1\pm5.5           $& 972   &\multicolumn{1}{r|}{$17.8\pm0.1$}  & $4\cdot 3^{20}$ & $ 0.3\pm0.1                $& \multicolumn{1}{r|}{$683.4\pm60.5$}   &$ 0.4\pm0.1                       $&$ -20.9\pm0.1  $         \\
                              & 3   & 1024  & \multicolumn{1}{r|}{$13.1\pm1.6$} & $4\cdot 4^{16}$ & $0.2\pm0.1$                & \multicolumn{1}{r|}{$151.0\pm11.1$}   & $0.1\pm0.0$                       & $-14.8\pm2.2$           & 4096  & \multicolumn{1}{r|}{$15.4\pm0.3$}  & $4\cdot 4^{20}$ & $0.3\pm0.2$                & \multicolumn{1}{r|}{$603.4\pm58.7$}   & $0.7\pm0.3$                       & $-17.0\pm0.2$           \\ \down
                              & 4   & 2500  & \multicolumn{1}{r|}{$11.5\pm0.0$} & $4\cdot 5^{16}$ & $0.2\pm0.1$                & \multicolumn{1}{r|}{$168.4\pm15.7$}   & $0.1\pm0.0$                       & $-12.6\pm0.1$           & 12500 & \multicolumn{1}{r|}{$100.0\pm0.0$} & $4\cdot 5^{20}$ & $0.2\pm0.1$                & \multicolumn{1}{r|}{$298.8\pm92.1$}   & $100.0\pm0.0$                     & $-7.3\pm0.1$            \\ \hline \up
30                            & 2   & 324   & \multicolumn{1}{r|}{$21.2\pm3.0$} & $4\cdot 3^{16}$ & $0.6\pm0.2$                & \multicolumn{1}{r|}{$402.6\pm45.8$}   & $0.2\pm0.1$                       & $-26.0\pm4.8$           & 972   & \multicolumn{1}{r|}{$18.6\pm0.2$} & $4\cdot 3^{20}$ & $0.2\pm0.1$                & \multicolumn{1}{r|}{$1067.0\pm133.0$} & $0.8\pm0.3$                       & $-21.5\pm0.1$           \\
                              & 3   & 1024  & \multicolumn{1}{r|}{$13.3\pm0.1$} & $4\cdot 4^{16}$ & $0.2\pm0.1$                & \multicolumn{1}{r|}{$579.8\pm64.2$}   & $0.2\pm0.1$                       & $-15.0\pm0.1$           & 4096  & \multicolumn{1}{r|}{$16.1\pm0.3$} & $4\cdot 4^{20}$  & $0.2\pm0.2$                & \multicolumn{1}{r|}{$918.0\pm78.1$}   & $1.1\pm0.3$                       & $-17.7\pm0.2$           \\ \down
                              & 4   & 2500  & \multicolumn{1}{r|}{$10.9\pm0.1$} & $4\cdot 5^{16}$ & $0.3\pm0.1$                & \multicolumn{1}{r|}{$204.0\pm12.9$}   & $0.2\pm0.1$                       & $-11.6\pm0.1$           & 12500 & \multicolumn{1}{r|}{$100.0\pm0.0$} & $4\cdot 5^{20}$ & $0.4\pm0.1$                & \multicolumn{1}{r|}{$406.0\pm46.8$}   & $100.0\pm0.0$                     & $-7.9\pm0.1$            \\ \hline \up
39                            & 2   & 324   & \multicolumn{1}{r|}{$13.0\pm0.1$} & $4\cdot 3^{16}$ & $0.3\pm0.1$                & \multicolumn{1}{r|}{$406.6\pm78.5$}   & $0.2\pm0.1$                       & $-14.4\pm0.1$           & 972   & \multicolumn{1}{r|}{$12.0\pm0.2$} & $4\cdot 3^{20}$  & $0.2\pm0.1$                & \multicolumn{1}{r|}{$195.4\pm16.5$}   & $0.3\pm0.2$                       & $-13.0\pm0.2$           \\
                              & 3   & 1024  & \multicolumn{1}{r|}{$8.5\pm0.1$} & $4\cdot 4^{16}$ & $0.1\pm0.1$                & \multicolumn{1}{r|}{$706.2\pm70.9$}   & $0.3\pm0.1$                       & $-8.8\pm0.1$            & 4096  & \multicolumn{1}{r|}{$11.2\pm7.8$} & $4\cdot 4^{20}$  & $0.2\pm0.1$                & \multicolumn{1}{r|}{$252.8\pm10.0$}   & $3.4\pm8.5$                       & $-8.6\pm0.1$            \\ \down
                              & 4   & 2500  & \multicolumn{1}{r|}{$5.2\pm0.0$} & $4\cdot 5^{16}$ & $0.1\pm0.1$                & \multicolumn{1}{r|}{$297.2\pm29.6$}   & $0.3\pm0.0$                       & $-5.1\pm0.1$            & 12500 & \multicolumn{1}{r|}{$100.0\pm0.0$} & $4\cdot 5^{20}$ & $0.1\pm0.0$                & \multicolumn{1}{r|}{$365.8\pm13.7$}   & $100.0\pm0.0$                     & $-4.1\pm0.0$            \\ \hline \up
48                            & 2   & 324   & \multicolumn{1}{r|}{$13.7\pm0.2$} & $4\cdot 3^{16}$ & $0.3\pm0.1$                & \multicolumn{1}{r|}{$5835.0\pm707.1$} & $0.4\pm0.2$                       & $-15.1\pm0.1$           & 972   & \multicolumn{1}{r|}{$11.8\pm0.4$} & $4\cdot 3^{20}$  & $0.2\pm0.1$                & \multicolumn{1}{r|}{$2320.8\pm935.7$} & $0.7\pm0.5$                       & $-12.3\pm0.1$           \\
                              & 3   & 1024  & \multicolumn{1}{r|}{$8.5\pm0.1$} & $4\cdot 4^{16}$ & $0.1\pm0.1$                & \multicolumn{1}{r|}{$2388.8\pm152.3$} & $0.4\pm0.1$                       & $-8.7\pm0.1$            & 4096  & \multicolumn{1}{r|}{$9.2\pm0.5$} & $4\cdot 4^{20}$  & $0.1\pm0.1$                & \multicolumn{1}{r|}{$1090.2\pm163.4$} & $1.5\pm0.5$                       & $-8.3\pm0.1$            \\ \down
                              & 4   & 2500  & \multicolumn{1}{r|}{$5.9\pm0.2$} & $4\cdot 5^{16}$ & $0.1\pm0.1$                & \multicolumn{1}{r|}{$2186.0\pm708.0$} & $0.7\pm0.2$                       & $-5.4\pm0.1$            & 12500 & \multicolumn{1}{r|}{$100.0\pm0.0$} & $4\cdot 5^{20}$ & $0.2\pm0.0$                & \multicolumn{1}{r|}{$1212.6\pm615.8$} & $100.0\pm0.0$                     & $-4.2\pm0.0$            \\ \hline
\end{tabular}}}
{\centering The values are averages over the five instances generated per instance configuration, with a 95\% confidence interval.}
\end{table}
\begin{figure}[]
\centering
\FIGURE
{\centering {\resizebox{0.35\textwidth}{!}{\begin{tikzpicture}
\begin{axis}[
    xlabel={Instances (\%)},
    xmin=0, xmax=100,
    ymin=-5, ymax=120,
    xtick={0,10,20,30,40,50,60,70,80,90,100},
    ytick={0,10,20,30,40,50,60,70,80,90,100},
    legend pos=north west,
    legend style={nodes={scale=0.8, transform shape}},
]
\addplot[
    color=black,
    ]
    table{saa_gap.dat};
    \addlegendentry{\small {\SAATP} $\left(\overline{gap}(\overline{x})(\%)\right)$}
\addplot[
    color=black,
    style=dashed,
    ]
    table{dep_gap.dat};
    \addlegendentry{\small {\DEOP}  $\left(GAP(\%)\right)$}
\end{axis}
\end{tikzpicture}}}
}
{\centering Comparison of the optimality gap obtained by the {\SAATP} and the {\DEOP}.\label{fig:opt-gap}}
{}
\end{figure} 
Table~\ref{tab:linea-vs-saa} shows the experimental results of the 4- and 5-facility instances. In the table, $|K|$ and $|K_t|$ are the set of scenarios in the {\DEOP} (Program~\eqref{Prog:discrete}) and in the deterministic equivalent of our transformed Program~\eqref{Prog:transformed}, respectively. Moreover, $GAP(\%)$ is the optimality gap for the {\DEOP}. For the {\DEOP}, all runs reached the time limit of two hours. We followed the procedure explained in Section~\ref{subsec:saa} for the {\SAATP} to obtain the feasible solution $\overline{x}$ and its estimators (see Appendix~\ref{app:saa-estimators} for details). Even though the number of scenarios $|K_t|$ of the transformed program is extremely large (i.e.,$|\mathcal{D}|{|\mathcal{W}|}^{|\mathcal{N_F}||\mathcal{P}|}$), the number $|K|$ of scenarios of the {\DEOP} is small (i.e., $|\mathcal{D}|{|\mathcal{W}|}^{|\mathcal{N_F}|}$) for the considered instances. Thus, we can compute the real stochastic value of the feasible solution $\overline{x}$, i.e., $v(\overline{x})=c^{\top}\overline{x} + g(\overline{x})$. For the {\SAATP}, we solve $M=50$ problems with $N=750$ sample scenarios because these values sufficed to obtain estimators with low variance and feasible, high-quality solutions. We present the over-estimator of the relative optimality gap, i.e., $\overline{gap}(\overline{x})(\%)={100\% \times [v(\overline{x})-\overline{v}_M^N]}/{v(\overline{x})}$. Because the lower-bound estimator $\overline{v}_M^N$ is a random variable, this over-estimator of the optimality gap may assume negative values. The computational time of the {\SAATP} is in seconds. For all runs, the {\SAATP} finished before the time limit of two hours. For comparison, we display the gap between the upper-bound values obtained by the {\SAATP} and the {\DEOP}, i.e., $GAP_{UB}(\%)={100\%\times [UB-v(\overline{x})]}/{UB}$, and the gap between their lower bounds, i.e., $GAP_{LB}(\%)={100\% \times [LB-\overline{v}_M^N]}/{LB}$, where $UB$ and $LB$ are the upper and lower bounds for the {\DEOP}, respectively. Next, we summarize our main findings.
\begin{result}\label{result:discretechoice-saa}
      The {\SAATP} produces a feasible solution $\overline{x}$ with a considerably low over-estimator of the optimality gap, $\overline{gap}(\overline{x})(\%)$ is between $0\%$ and $1\%$ for all instances. The lower bound $\overline{v}_{M}^N$ and the optimality gap estimators have low variability, as shown by the standard deviation  $\hat{\sigma}_{N,M}$ of the lower-bound estimator, which is inferior to $0.5\%$ of $\overline{v}_{M}^N$ for all instances (see Appendix~\ref{app:additional-results-comp-baseline} for details). For the considered instances, the computational time required by the {\SAATP} was inferior to 2000s in 87\% of the cases, with an average of around 990s.
\end{result}
\begin{result}\label{result:linearization-saa}
    The {\SAATP}'s over-estimator of the optimality gap is significantly smaller than the {\DEOP}'s optimality gap (the difference is around $24\%$ on average). The {\SAATP} obtained these optimality gaps in significantly less time than the {\DEOP}, and, within the time limit, consistently found feasible solutions with slightly better values than those from the {\DEOP}.
\end{result}

Result~\ref{result:discretechoice-saa} highlights the good performance of the {\SAATP}. Result~\ref{result:linearization-saa} shows that it outperformed the {\DEOP}. Figure~\ref{fig:opt-gap}, which shows the optimality gap for the {\DEOP} and the {\SAATP}, corroborates this analysis. The {\SAATP} obtained optimality gaps close to $0\%$ for all instances, whereas the {\DEOP} obtained optimality gaps inferior to $10\%$ only for $22\%$ of the instances. The {\DEOP} obtained optimality gaps of $10\%$--$20\%$ for around $56\%$ of the instances, of $20\%$--$30\%$ for around $8\%$ of the instances, and of $100\%$ for around $14\%$ of the instances. Even though the transformed Program~\eqref{Prog:transformed} has more scenarios than the {\DEOP}, the {\SAATP} solves different problem instances with fewer scenarios. Therefore, although solving the {\DEOP} is constrained to a relatively small number of scenarios, the transformed program can handle larger instances as we can use sampling techniques to address it. The {\SAATP} obtained high-quality solutions by sampling only $N=750$ scenarios for each of the $M=50$ {\SAA} problems, which,  in most cases, is less than the number $|K|$ of scenarios of the {\DEOP}. Despite the number of scenarios, the {\SAATP} solves mathematical programs smaller than the {\DEOP} (in terms of number of variables and constraints). The {\DEOP}'s {\MILP} formulation contains numerous big-M constraints due to the linearization of the probability function~\eqref{probability}, which leads to weak linear relaxations. In contrast, the transformation applied to obtain the {\NDFPP}-Selection formulation (see Section~\ref{subsec:discretechoice-problem}) includes no extra variables or constraints (see Appendix~\ref{app:vars-constrs} for a comparison of the number of variables and constraints in the {\DEOP} and the {\SAATP}). Finally, although the {\SAA} performed well in this case study, other sampling methods can be employed based on the problem's specific requirements.

\subsubsection{Analysis for {\NDFPP}-Binomial, {\NDFPP}-Discrete, and {\NDFPP}-Normal.} \label{subsec:results-estimators}

In this section, we analyze the statistical estimators and feasible solution $\overline{x}$ of the {\SAATP} for {\NDFPP}-Binomial, {\NDFPP}-Discrete, and {\NDFPP}-Normal, which are {\MILP}s (see Sections~\ref{subsec:binomial-problem}-\ref{subsec:normal-problem}). In these cases, the {\DEOP} (Program~\eqref{Prog:discrete} for discrete distributions) is nonlinear nonconvex, and the linearization of~\citet{bhuiyan2020stochastic} cannot be applied. Moreover, note that our transformed Program~\eqref{Prog:transformed} for these variants has continuous random variables, making the explicit enumeration of scenarios impossible and requiring the implementation of sampling techniques. The {\DEOP} has a small number $|K|$ of scenarios for the {\NDFPP}-Binomial and {\NDFPP}-Discrete instances. Thus, we compute $v(\overline{x})$ for these cases. The number of scenarios of the {\DEOP} is infinite for {\NDFPP}-Normal because the endogenous random variables are continuous. Consequently, we compute the estimator $\hat{v}_{N{'}}(\overline{x})$ of $v(\overline{x})$. Tables~\ref{tab:saa-results} and~\ref{tab:saa-results-1} lists the experimental results of the 4- and 5-facility instances. The over-estimator of the optimality gap $\overline{gap}(\overline{x})(\%)$ is defined as in the previous section for {\NDFPP}-Binomial and {\NDFPP}-Discrete. For  {\NDFPP}-Normal, $\overline{gap}(\overline{x})(\%)={100\% \times [\hat{v}_{N{'}}(\overline{x})-\overline{v}_M^N]}/{\hat{v}_{N{'}}(\overline{x})}$.
We set $M=50$ and $N=750$ for the same reason as in the previous section. For {\NDFPP}-Normal, we use $N{'}=5\times10^4$. The {\SAATP} finished for all instances before reaching 3 hours. The result is as follows:
\begin{result}\label{res:saa-estimators}
    The {\SAATP} obtained a feasible solution $\overline{x}$ with a significantly low statistical over-estimator of the optimality gap $\overline{gap}(\overline{x})(\%)$ between $-1\%$ and $1\%$ for all instances of the problem versions. The {\SAA}'s estimators have high confidence, as shown by the standard deviation $\hat{\sigma}_{N,M}$  of the lower-bound estimator, which is inferior to $0.5\%$ of $\overline{v}_M^N$ for all cases, and the standard deviation  $\hat{\sigma}_{N{'}}(\overline{x})$ of the upper bound estimator, which is inferior to $0.2\%$ of $\hat{v}_{N{'}}(\overline{x})$ for all {\NDFPP}-Normal instances (see Appendix~\ref{app:additional-statistical-estimators} for details). The computational time was, on average, around 1500s for {\NDFPP}-Binomial, 3000s for {\NDFPP}-Discrete, and 1200s for {\NDFPP}-Normal.
\end{result}

\begin{table}[]
\centering
\TABLE
{Experimental results for {\NDFPP}-Binomial and {\NDFPP}-Discrete of the {\SAATP}. \label{tab:saa-results}}
{\begin{adjustbox}{width=0.85\textwidth,center}%
\begin{tabular}{rr|rrrrrr|rrrrrr}
\hline 
\multicolumn{2}{c|}{\multirow{2}{*}{Network}} & \multicolumn{6}{c|}{{\NDFPP}-Binomial}                                                                                              & \multicolumn{6}{c}{{\NDFPP}-Discrete}                                                                                               \\ \cline{3-14} 
\multicolumn{2}{c|}{}                         & \multicolumn{3}{c|}{Four-facility instances}                                        & \multicolumn{3}{c|}{Five-facility instances}                   & \multicolumn{3}{c|}{Four-facility instances}                                        & \multicolumn{3}{c}{Five-facility instances}                    \\ \hline \down \up
$|\mathcal{N}_\mathcal{T}^0|$      & $W$      & $|K|$ & $\overline{gap}(\overline{x})(\%)$ & \multicolumn{1}{r|}{Time (s)}          & $|K|$ & $\overline{gap}(\overline{x})(\%)$ & Time (s)          & $|K|$ & $\overline{gap}(\overline{x})(\%)$ & \multicolumn{1}{r|}{Time (s)}          & $|K|$ & $\overline{gap}(\overline{x})(\%)$ & Time (s)          \\ \hline \up
15                                 & 2        & 324   & $0.5\pm0.1$                        & \multicolumn{1}{r|}{$1929.0\pm1362.8$} & 972   & $0.5\pm0.1$                        & $1749.4\pm167.4$  & 324   & $0.4\pm0.1$                        & \multicolumn{1}{r|}{$2373.6\pm871.8$}  & 972   & $0.3\pm0.1$                        & $2291.2\pm121.0$  \\
                                   & 3        & 1024  & $0.5\pm0.1$                        & \multicolumn{1}{r|}{$1890.2\pm973.9$}  & 4096  & $0.3\pm0.1$                        & $1200.0\pm180.4$  & 1024  & $0.4\pm0.1$                        & \multicolumn{1}{r|}{$3520.4\pm957.4$}  & 4096  & $0.2\pm0.1$                        & $2642.6\pm310.2$  \\ \down
                                   & 4        & 2500  & $0.3\pm0.1$                        & \multicolumn{1}{r|}{$1507.6\pm583.5$}  & 12500 & $0.6\pm0.1$                        & $1375.0\pm73.9$   & 2500  & $0.1\pm0.1$                        & \multicolumn{1}{r|}{$3733.6\pm799.9$}  & 12500 & $0.6\pm0.0$                        & $3898.0\pm692.2$  \\ \hline \up 
25                                 & 2        & 324   & $0.7\pm0.1$                        & \multicolumn{1}{r|}{$458.6\pm115.6$}   & 972   & $0.6\pm0.1$                        & $1015.6\pm169.8$  & 324   & $0.6\pm0.1$                        & \multicolumn{1}{r|}{$1287.2\pm160.0$}  & 972   & $0.4\pm0.1$                        & $2211.6\pm348.2$  \\
                                   & 3        & 1024  & $0.4\pm0.1$                        & \multicolumn{1}{r|}{$317.4\pm15.3$}    & 4096  & $0.5\pm0.1$                        & $1092.2\pm81.5$   & 1024  & $0.3\pm0.1$                        & \multicolumn{1}{r|}{$1834.2\pm205.3$}  & 4096  & $0.3\pm0.1$                        & $2885.6\pm153.5$  \\ \down
                                   & 4        & 2500  & $0.4\pm0.1$                        & \multicolumn{1}{r|}{$390.6\pm24.7$}    & 12500 & $0.2\pm0.1$                        & $601.0\pm97.9$    & 2500  & $0.2\pm0.1$                        & \multicolumn{1}{r|}{$2441.0\pm173.7$}  & 12500 & $0.2\pm0.0$                        & $3470.0\pm460.0$  \\ \hline \up 
30                                 & 2        & 324   & $0.6\pm0.1$                        & \multicolumn{1}{r|}{$683.2\pm32.9$}    & 972   & $0.5\pm0.1$                        & $1507.6\pm166.9$  & 324   & $0.5\pm0.1$                        & \multicolumn{1}{r|}{$1502.4\pm146.2$}  & 972   & $0.3\pm0.1$                        & $2535.4\pm177.8$  \\
                                   & 3        & 1024  & $0.4\pm0.1$                        & \multicolumn{1}{r|}{$957.8\pm100.6$}   & 4096  & $0.3\pm0.1$                        & $1489.4\pm160.8$  & 1024  & $0.2\pm0.1$                        & \multicolumn{1}{r|}{$2525.8\pm55.9$}   & 4096  & $0.2\pm0.1$                        & $3209.4\pm179.1$  \\ \down
                                   & 4        & 2500  & $0.5\pm0.1$                        & \multicolumn{1}{r|}{$498.2\pm38.4$}    & 12500 & $0.5\pm0.1$                        & $720.4\pm98.1$    & 2500  & $0.3\pm0.1$                        & \multicolumn{1}{r|}{$2578.0\pm149.4$}  & 12500 & $0.4\pm0.0$                        & $3218.6\pm294.8$  \\ \hline \up  
39                                 & 2        & 324   & $0.4\pm0.0$                        & \multicolumn{1}{r|}{$681.2\pm115.0$}   & 972   & $0.3\pm0.1$                        & $391.2\pm35.5$    & 324   & $0.2\pm0.1$                        & \multicolumn{1}{r|}{$1457.2\pm81.6$}   & 972   & $0.2\pm0.1$                        & $1758.2\pm89.1$   \\
                                   & 3        & 1024  & $0.3\pm0.1$                        & \multicolumn{1}{r|}{$1133.4\pm69.4$}   & 4096  & $0.2\pm0.1$                        & $543.8\pm37.2$    & 1024  & $0.2\pm0.1$                        & \multicolumn{1}{r|}{$2350.6\pm206.3$}  & 4096  & $0.2\pm0.1$                        & $2469.6\pm243.0$  \\ \down
                                   & 4        & 2500  & $0.1\pm0.0$                        & \multicolumn{1}{r|}{$566.8\pm57.4$}    & 12500 & $0.1\pm0.1$                        & $757.2\pm59.0$    & 2500  & $0.1\pm0.1$                        & \multicolumn{1}{r|}{$2509.8\pm268.1$}  & 12500 & $0.1\pm0.0$                        & $4135.6\pm444.9$  \\ \hline \up 
48                                 & 2        & 324   & $0.4\pm0.0$                        & \multicolumn{1}{r|}{$7520.0\pm967.4$}  & 972   & $0.3\pm0.1$                        & $2774.0\pm1377.4$ & 324   & $0.2\pm0.1$                        & \multicolumn{1}{r|}{$7730.4\pm1258.2$} & 972   & $0.2\pm0.1$                        & $4639.0\pm1769.7$ \\
                                   & 3        & 1024  & $0.3\pm0.1$                        & \multicolumn{1}{r|}{$3717.6\pm449.1$}  & 4096  & $0.2\pm0.1$                        & $1585.0\pm281.4$  & 1024  & $0.2\pm0.1$                        & \multicolumn{1}{r|}{$4984.2\pm521.0$}  & 4096  & $0.1\pm0.1$                        & $3166.4\pm189.6$  \\ \down
                                   & 4        & 2500  & $0.2\pm0.0$                        & \multicolumn{1}{r|}{$2767.4\pm523.7$}  & 12500 & $0.2\pm0.0$                        & $1526.2\pm574.6$  & 2500  & $0.2\pm0.1$                        & \multicolumn{1}{r|}{$4989.4\pm288.8$}  & 12500 & $0.2\pm0.0$                        & $5558.6\pm348.3$  \\ \hline 
\end{tabular}%
\end{adjustbox}}
{\centering The values are averaged over the five instances generated per instance configuration, with a 95\% confidence interval.}
\end{table}

\begin{table}[]
\centering
\TABLE
{Experimental results for {\NDFPP}-Normal of the {\SAATP}. \label{tab:saa-results-1}}
{\begin{adjustbox}{width=0.35\textwidth,center}%
\begin{tabular}{r|rr|rr}
\hline
\multicolumn{1}{c|}{\multirow{2}{*}{$|\mathcal{N}_\mathcal{T}^0|$}} & \multicolumn{2}{c|}{Four-facility instances}                                           & \multicolumn{2}{c}{Five-facility instances}                                           \\ \cline{2-5} 
\multicolumn{1}{c|}{}                                               & \multicolumn{1}{c}{$\overline{gap}(\overline{x})(\%)$} & \multicolumn{1}{c|}{Time (s)} & \multicolumn{1}{c}{$\overline{gap}(\overline{x})(\%)$} & \multicolumn{1}{c}{Time (s)} \\ \hline \up
15                                                                  & $-0.2\pm0.2$                                           & $1478.0\pm256.3$              & $-0.3\pm0.3$                                           & $878.8\pm51.2$               \\
25                                                                  & $-0.2\pm0.2$                                           & $851.4\pm63.0$                & $-0.3\pm0.3$                                           & $842.6\pm39.2$               \\
30                                                                  & $-0.2\pm0.2$                                           & $1127.0\pm145.5$              & $-0.3\pm0.3$                                           & $1058.6\pm104.9$             \\
39                                                                  & $-0.2\pm0.2$                                           & $1436.0\pm76.9$               & $-0.3\pm0.3$                                           & $1352.4\pm126.3$             \\ \down
48                                                                  & $-0.2\pm0.2$                                           & $1837.2\pm125.6$              & $-0.2\pm0.3$                                           & $1746.6\pm72.1$              \\ \hline
\end{tabular}
\end{adjustbox}}
{\centering The values are averaged over the five instances generated per instance configuration, with a 95\% confidence interval.}
\end{table}

The {\SAATP} found high-quality solutions and estimators for the considered problems. The computational time of this approach for the {\NDFPP}-Binomial, {\NDFPP}-Discrete, and {\NDFPP}-Normal exceeded that for the {\NDFPP}-Selection. Solving the {\NDFPP}-Discrete required the longest computational time. The transformation used in this case (i.e., inversion for discrete random variables; see Section~\ref{subsubsec:discrete-rv}) includes numerous binary variables and big-M constraints in the formulation. The initial formulation of {\NDFPP}-Discrete contained bilinear terms, which were linearized by adding more big-M constraints in the formulation. For the same reason, the computational time of the {\SAATP} for {\NDFPP}-Binomial exceeded that for {\NDFPP}-Selection. The computational time of the {\SAATP} for solving {\NDFPP}-Discrete increases with the number of capacity levels (possible realizations of the random variables). Therefore, a trade-off exists between how sample size affects the performance and confidence of the {\SAA}'s statistical estimators and the associated computational time. For {\NDFPP}-Normal, the number of evaluation scenarios $N{'}$ is considerably large, which increases the computational time of the {\SAATP}. Moreover, the standard deviation of the optimality gap is slightly larger than for the remaining versions because it includes the variability of $\hat{v}_{N{'}}(\overline{x})$.

\subsubsection{Value of the Stochastic Solution.} \label{subsec:results-EEV}

This section analyzes the relative {\VSS} of all problem variants. The relative {\VSS} is defined as ${(\textrm{{\EEV}} - v^{*} )}/{\textrm{{\EEV}}}$, where $v^{*}$ is the optimal value of the recourse problem, and {\EEV} is ${v(\hat{x}^{*})= c^{\top}\hat{x}^{*} + g(\hat{x}^{*})}$, with $\hat{x}^{*}$ the optimal solution of the {\EV}. In Appendix~\ref{app:EV}, we detail the {\EV}s (i.e., the deterministic problems with the average case scenario) according to the problem variants. We use the feasible solution $\overline{x}$ obtained by the {\SAATP} to compute under estimators of the {\VSS}, and we detail this computation in Appendix~\ref{app:vss}. We resume our main findings in Result~\ref{res:vss}.
\begin{result}\label{res:vss}
    The {\SAATP} obtained better solutions than the {\EV} in all cases. Thus, solving the recourse problem instead of the {\EV} is also advantageous in the context of endogenous uncertainty, where the distribution defined by  $\hat{x}^{*}$ is not the same as that defined by $\overline{x}$ or $x^{*}$. Different instances {of \NDFPP}-Selection showed a high {\VSS}, around $30\%$ or more. Some instances {of \NDFPP}-Binomial and {\NDFPP}-Discrete also presented high {\VSS}s, whereas others had VSSs of  $1\%$--$4\%$. Although the VSS for {\NDFPP}-Normal approached $0\%$ for the 5-facility instances, it hovered around $3\%$--$6\%$ for the 4-facility instances. Although the {\VSS} is primarily instance-dependent, it is also affected by the problem (e.g., {\NDFPP}-Selection and {\NDFPP}-Normal for 15-node instances have considerably different {\VSS}s). Finally, Appendix~\ref{app:vss} contains the complete results.
\end{result}

\section{Conclusions}\label{conclusion}

We proposed a general method for addressing {\SMIPD} that uses random variable transformation. Our methodology
\begin{enumerate*}[label=(\roman*)]
    \item works for a general class of programs with endogenous uncertainty;
    \item obtains a stochastic program with exogenous uncertainty for which solution procedures are well-studied;
    \item yields {\MILP}s or {\MICP}s with exogenous uncertainty in many cases;
    \item works for both discrete and continuous endogenous random variables; and
    \item is applicable to problems with independent and correlated random variables.
\end{enumerate*}
We investigated techniques to define appropriate transformations of random variables to treat cases where the decision-dependent random variables follow classical probability distributions. We showed the existence of a simple transformation for the case of the discrete selection of distributions, which is the most commonly studied in previous research. We also demonstrated the proposed methodology by applying it to the {\NDFPP}, a broadly studied problem. In the computational experiments, the {\SAATP} produced feasible, high-quality solutions and high-confidence estimators in a short computational time.

Numerous research perspectives stem from this work. Future work could focus on 
\begin{enumerate*}[label=(\roman*)]
\item exploring additional techniques to define transformation functions;
\item investigating approaches to deal with the current limitations of the method when $q(\omega)$, $W(\omega)$ and/or $H(\omega)$ are random, such as exploring valid inequalities for the discrete random variables inverse transformation;
\item analyzing how endogenous uncertainty affects the definition of the VSS;
\item examining stochastic programs with endogenous random variables for which the probability distribution is not known \emph{a priori}; and
\item extending the {\NDFPP} to account for assumptions investigated in previous research, such as post-disruption edge failure, protection (see, e.g.,~\citep{du2014stochastic}), and facility location~\citep{li2021stochastic,li2022locating}.
\end{enumerate*}

\ACKNOWLEDGMENT{This research was funded by the SCALE-AI Chair in Data Driven Supply Chains, the FRQ-IVADO Research Chair in Data Science for Combinatorial Game Theory, and  NSERC Grant No. 2024-04051. It was also enabled by support provided by Calcul Québec (www.calculquebec.ca) and the Digital Research Alliance of Canada (https://alliancecan.ca/).}

\bibliographystyle{class_paper}
\bibliography{myrefs}

\begin{APPENDICES}
\section{Description of the {\SAA} Method and its Statistical Estimators}\label{app:saa-estimators}

We detail the {\SAA} method in Algorithm~\ref{alg:saa}. We define $M$ different instances of the transformed {\SAA} Program~\eqref{Prog:saa}, each with $N$ independent scenarios. Moreover, the scenarios from distinct instances are also independently generated. The $\hat{v}_N^j$ is the value of the optimal solution $\hat{x}_N^j$ of the transformed {\SAA} instance $j \in \mathcal{M}=\{1,\ldots,M\}$. For the sake of clarity, we simplified the {\SAA} method in Algorithm~\ref{alg:saa}. In practice, $M$, $N$, and/or $N{'}$ should be increased and Algorithm~\ref{alg:saa} repeated as long as the obtained variances are large. Additionally, the $M$ instances of the transformed {\SAA} program and the evaluation problem to select the solution $\overline{x}$ (Steps~1 and~2 of Algorithm~\ref{alg:saa}, respectively) are solved sequentially. The conduction of the validation analyses in Step~3 of the algorithm is applicable when the transformed Program~\eqref{Prog:transformed} has relatively complete recourse. Thus, next, we detail the computation of the statistical estimators in Step~3.
\begin{algorithm}
    \small
    \caption{{\SAA} method} \label{alg:saa} 
    \Input $( c^{\top}, A, b)$, $\Omega$, $\boldsymbol{\xi}_x$, $\boldsymbol{\vartheta}$, $\mathcal{X}$, $\mathcal{Y}$, $M$, $\mathcal{M}$, $N$, $N{'}$, and $\alpha$;\\
    \Output $\overline{x}$, $\overline{v}_{N}^M$, $\hat{v}_{N{'}}(\overline{x})$, $\overline{gap}(\overline{x})$, $\hat{\sigma}^2_{N,M}$, $\hat{\sigma}^2_{N{'}}(\overline{x})$, $\hat{\sigma}^2_{gap}(\overline{x})$, $v_{lb}$, $v_{ub}(\overline{x})$, and $sgap(\overline{x})$; \\ 
    \Stepone Solve $M$ {\SAA} Programs~\eqref{Prog:saa} with $N$ scenarios to obtain $\hat{x}_N^j$ and $\hat{v}_N^j$, $j \in \mathcal{M}$; \\
    \Steptwo  Set $k = \arg \min_{j \in \mathcal{M}}\{c^{\top}\hat{x}_N^j + \hat{g}_{N{'}}(\hat{x}_N^j)\}$ and $\overline{x} = \hat{x}_N^k$;  \\
    \Stepthree Compute $\overline{v}_{N}^M$, $\hat{v}_{N{'}}(\overline{x})$, $\overline{gap}(\overline{x})$, $\hat{\sigma}^2_{N,M}$, $\hat{\sigma}^2_{N{'}}(\overline{x})$, $v_{lb}$, $v_{ub}(\overline{x})$, and $sgap(\overline{x})$;\\
    \Return $\overline{x}$, $\overline{v}_{N}^M$, $\hat{v}_{N{'}}(\overline{x})$, $\overline{gap}(\overline{x})$, $\hat{\sigma}^2_{N,M}$, $\hat{\sigma}^2_{N{'}}(\overline{x})$, $\hat{\sigma}^2_{gap}(\overline{x})$, $v_{lb}$, $v_{ub}(\overline{x})$, and $sgap(\overline{x})$.
\end{algorithm}
The unbiased estimator of $\mathbb{E}[v_N]$ is defined as
\begin{equation}
    \overline{v}_N^M = \frac{1}{M}\sum_{j \in \mathcal{M}} \hat{v}_N^j. 
\end{equation}
Therefore, a statistical lower bound of $v^{*}$ with $100(1-\alpha)\%$ confidence is  
\begin{equation}
    v_{lb} = \overline{v}_N^M - \mathfrak{t}_{\alpha,M-1} \hat{\sigma}_{N,M}, 
\end{equation}
where $\mathfrak{t}_{\alpha,M-1}$ is the $\alpha$ critical value of the $\mathfrak{t}$-student distribution with $M-1$ degrees of freedom, and $\hat{\sigma}^2_{N,M}$ is the unbiased estimator of the variance of $\overline{v}_N^M$:
\begin{equation}
    \hat{\sigma}^2_{N,M} = \frac{1}{M \left(M-1\right)} \sum_{j=1}^{M} \left( \hat{v}_N^j - \overline{v}_N^M \right)^2. 
\end{equation}    
To obtain a statistical upper bound, and considering a feasible first-stage solution such as $\overline{x} \in \mathcal{X}$ (see  Algorithm~\ref{alg:saa}), we generate a sufficiently larger number of scenarios ($N{'} \gg N$) and obtain the value of the unbiased estimator of $g(\overline{x})$, which is $\hat{g}_{N{'}}(\overline{x})$. The estimator of the stochastic value of this solution is $\hat{v}_{N{'}}(\overline{x})= c^{\top}\overline{x} + \hat{g}_{N{'}}(\overline{x})$. Thus, the statistical upper bound with $100(1-\alpha)\%$ confidence is
\begin{equation}
    v_{ub}(\overline{x}) = \hat{v}_{N{'}}(\overline{x}) + \mathfrak{z}_{\alpha}\hat{\sigma}_{N{'}}(\overline{x}), 
\end{equation}
where $\mathfrak{z}_{\alpha}$ is the $\alpha$ critical value of the standard normal distribution, and $\hat{\sigma}_{N{'}}^2(\overline{x})$ is the unbiased estimator of the variance of $\hat{v}_{N{'}}(\overline{x})$, which is calculated as
\begin{equation}
    \hat{\sigma}_{N{'}}^2(\overline{x}) = \frac{1}{N \left(N{'}-1\right)} \sum_{i=1}^{N{'}} \left[ Q(\overline{x},\vartheta^i) - \hat{g}_{N{'}}(\overline{x}) \right]^2. 
\end{equation}
Finally, $v_N$ is a downward-biased estimator of $v^{*}$ (i.e., $v^{*} = v_{t(x,\boldsymbol{\vartheta})}^{*} \geq \mathbb{E}[v_N]$)~\citep[Chap. 5]{shapiro2021lectures}. Thus, following the bias of $v_N$, the over-estimator of $c^{\top}\overline{x} + g(\overline{x}) - v_{t(x,\boldsymbol{\vartheta})}^{*}$, which is equal to the optimality gap $\overline{gap}(\overline{x})=c^{\top}\overline{x} + g(\overline{x}) - v^{*}$, is the difference $\overline{gap}(\overline{x})= \hat{v}_{N{'}}(\overline{x})-\overline{v}^M_N$. The variance of this over-estimator is $\hat{\sigma}^2_{gap}(\overline{x})=\hat{\sigma}_{N{'}}^2(\overline{x}) + \hat{\sigma}^2_{N,M}$. Thus, a statistical gap estimator with the same confidence as the previous estimators is
\begin{equation}
   sgap(\overline{x}) = \hat{v}_{N{'}}(\overline{x}) - \overline{v}_M^N + \mathfrak{z}_{\alpha} \sqrt{\left(\hat{\sigma}_{N{'}}^2(\overline{x}) + \hat{\sigma}^2_{N,M}\right)}, 
\end{equation}
where $\mathfrak{z}_{\alpha}$ is again the $\alpha$ critical value of the standard normal distribution. For a detailed analysis of these statistical estimators, please refer to~\citet[Ch.\ 5]{shapiro2021lectures}, Section 5.6.

\section{Proofs of Theorems~\ref{theorem_general_discrete} and~\ref{theorem:continuous-dist-conv}} \label{app:proofs}

\proof{Proof of Theorem~\ref{theorem_general_discrete}.}
If the functions $\hat{p}^u_{r}(x)$ are concave over $\left[0,1\right]^{n_1}$, the sum $\sum_{j=1}^{r} \hat{p}_{j}^u(x)$ has the same property and the continuous relaxation of Constraints~\eqref{general_drv_1} and~\eqref{general_drv_2} are convex. If the functions $\hat{p}_{r}^u(x)$ are convex over $\left[0,1\right]^{n_1}$, we can reformulate Constraints~\eqref{general_drv_1} and~\eqref{general_drv_2} using 
\begin{equation}
    \label{equi1}\sum_{j=1}^{r} \hat{p}_{j}^u(x) = 1 - \sum_{j=r+1}^{R_u} \hat{p}_{j}^u(x), 
\end{equation}
and we obtain constraints with a convex linear relaxation. \Halmos
\endproof

\proof{Proof of Theorem~\ref{theorem:continuous-dist-conv}.}
From Criterion~\ref{item1}, we have $t_1\boldsymbol{(} x, \widetilde{q}(\omega)^{\top} \boldsymbol{)} = q^{\top}$, $t_3\boldsymbol{(} x, \widetilde{T}(\omega) \boldsymbol{)} = T$ and $t_4\boldsymbol{(} x, \widetilde{W}(\omega) \boldsymbol{)} = W$, so the objective function and the left-hand side of the recourse constraint of Program~\eqref{Prog:2stagetransformed} are linear. Additionally, Criterion~\ref{item2} ensures that $t(x,\boldsymbol{\vartheta_u}) = t(x,\boldsymbol{h}^u) = \mathbb{F}^{-1}_{\boldsymbol{\xi}^u}(x,r)(\cdot)$ is the common inverse function, which has a convex continuous relaxation on $x \in \mathcal{X}$ by Criterion~\ref{item3}. Therefore, the recourse constraint in Program~\eqref{Prog:2stagetransformed} is convex. \Halmos
\endproof 

\section{Examples of Transformation Techniques}\label{app:examples}

This section illustrates the application of the inversion (Example~\ref{ex:inverse}),  enumeration (Example~\ref{ex:discrete-selection}), and convolution (Example~\ref{ex:correlated}) techniques for obtaining appropriate transformation functions.

\begin{example}\label{ex:inverse}
    Consider two independent (correlated) random variables $\boldsymbol{\xi}_x^A$ and $\boldsymbol{\xi}_x^B$ following the marginal cumulative distribution $\mathbb{F}_{\boldsymbol{\xi}^A}(x,\beta)(\cdot)$ and $\mathbb{F}_{\boldsymbol{\xi}^B}(x,\beta)(\cdot)$, respectively. We have two independent (correlated) random variables $\boldsymbol{\vartheta}_A$ and $\boldsymbol{\vartheta}_B$ with marginal distribution $U(0,1)$. For the independent case, the $\boldsymbol{\vartheta}_A$ and $\boldsymbol{\vartheta}_B$ follow the independence bi-variate copula. For the correlated case, $\boldsymbol{\vartheta}_A$ and $\boldsymbol{\vartheta}_B$ follow a bi-variate copula modeling the correlation of $\boldsymbol{\xi}_x^A$ and $\boldsymbol{\xi}_x^B$, e.g., the Gaussian copula with the appropriate correlation. A scenario is composed of the realizations of $\boldsymbol{\vartheta}_A$ and $\boldsymbol{\vartheta}_B$, which are transformed into $\boldsymbol{\xi}_x^A$ and $\boldsymbol{\xi}_x^B$ by $t(x,\boldsymbol{\vartheta}_u) =\mathbb{F}^{-1}_{\boldsymbol{\xi}^u}(x,\beta)$ for $u \in \{A,B\}$.
\end{example}

\begin{example}\label{ex:discrete-selection}
Consider two independent (correlated) random elements, $\boldsymbol{\xi}_x^{A}$ and $\boldsymbol{\xi}_x^{B}$, with realization sets $\mathcal{R}_{A} = \{r_A^1,r_A^2,r_A^3\}$ and $\mathcal{R}_{B} = \{r_B^1,r_B^2,r_B^3\}$, respectively, resulting in nine scenarios for the original program. These variables have alternative marginal distributions $\mathcal{D}_{A} = \{d_A^1,d_A^2\}$ and $\mathcal{D}_{B} = \{d_B^1,d_B^2\}$. In distribution $d_A^1$, only $r_A^1$ and $r_A^2$ have non-zero probabilities, whereas $d_A^2$ assigns probabilities to all three realizations. The same applies to $d_B^1$ and $d_B^2$. If $d_A^1$ and $d_B^1$ are selected, only four scenarios have non-zero probabilities; if $d_A^2$ and $d_B^2$ are chosen, all nine scenarios do. We have two independent (correlated) exogenous random vectors: $\boldsymbol{\vartheta}_{A} = \{\boldsymbol{\vartheta}_{A,d_A^1}, \boldsymbol{\vartheta}_{A,d_A^2}\}$ and $\boldsymbol{\vartheta}_{B} = \{\boldsymbol{\vartheta}_{B,d_B^1},\boldsymbol{\vartheta}_{B,d_B^2}\}$. To obtain $ \boldsymbol{\vartheta}_{A} $, we consider a random variable $\boldsymbol{\varepsilon}_A$ with a marginal distribution $U(0,1)$, and use the inverse transformation method to derive $\boldsymbol{\vartheta}_{A,d_A^1}$ and $\boldsymbol{\vartheta}_{A,d_A^2}$ according to the distributions $d_A^1$ and $ d_A^2$, respectively. The same approach applies to $\boldsymbol{\vartheta}_B$. As in Example~\ref{ex:inverse}, for the independent case, $\boldsymbol{\varepsilon}_A$ and $\boldsymbol{\varepsilon}_B$ follow an independent bi-variate copula. For the correlated case, $\boldsymbol{\varepsilon}_A$ and $\boldsymbol{\varepsilon}_B$ are modeled using a bi-variate copula that captures the correlation between $\boldsymbol{\xi}_x^A$ and $\boldsymbol{\xi}_x^B$. The random variable $\boldsymbol{\vartheta}_{A,d_A^1}$ follows distribution $d_A^1$ and has two realizations ($r_A^1$ and $r_A^2$), while $\boldsymbol{\vartheta}_{A,d_A^2}$ follows distribution $d_A^2$ and has three realizations ($r_A^1$, $r_A^2$, and $r_A^3$). The same applies to $\boldsymbol{\vartheta}_{B,d_B^1}$ and $\boldsymbol{\vartheta}_{B,d_B^2}$. Thus, the transformed program has 36 scenarios. If distributions $d_A^1$ and $d_B^1$ are chosen, then $t(x,\boldsymbol{\vartheta}_{A})$ and $t(x,\boldsymbol{\vartheta}_{B})$ return random variables $\boldsymbol{\vartheta}_{A,d_A^1}$ and $\boldsymbol{\vartheta}_{B,d_B^1}$, respectively. For the scenario with realizations $r_A^1$, $r_A^3$, $r_B^1$, and $r_B^3$ for random variables $\boldsymbol{\vartheta}_{A,d^1_A}$, $\boldsymbol{\vartheta}_{A,d^2_A}$, $\boldsymbol{\vartheta}_{B,d^1_B}$, and $\boldsymbol{\vartheta}_{B,d^2_B}$, respectively, $t(x,\boldsymbol{\vartheta}_{A})$ and $t(x,\boldsymbol{\vartheta}_{B})$ return realizations $r_A^1$ and $r_B^1$.
\end{example}

\begin{example}\label{ex:discrete-selection-all}
 If an enumeration of the alternative distributions $\mathcal{D}$ of the entire random vector $\boldsymbol{\xi}$ is viable, then we can apply the same transformation considering all the elements of $\boldsymbol{\xi}$. That is, ${t(x,\boldsymbol{\vartheta})=\sum_{d = 1}^{D}\boldsymbol{\vartheta}_{d}x^d}$, where $\boldsymbol{\vartheta} = \{\boldsymbol{\vartheta}_d\}_{d \in \mathcal{D}}$, $\boldsymbol{\vartheta}_d$ follows distribution $d \in \mathcal{D}$, and $x^d$ is defined as in Section~\ref{subsec:background}.
\end{example}

\begin{example}\label{ex:correlated}
    Consider two correlated random elements $\boldsymbol{\xi}_x^A$ and $\boldsymbol{\xi}_x^B$, such that $\boldsymbol{\xi}_x^A = \boldsymbol{\zeta}_x^{A1} + \boldsymbol{\zeta}_x^{AB2}$ and  $\boldsymbol{\xi}_x^B = \boldsymbol{\zeta}_x^{B1} + \boldsymbol{\zeta}_x^{AB2}$. In this case, the correlation is a result of the dependency on $\boldsymbol{\zeta}_x^{AB2}$. Therefore, if $\boldsymbol{\zeta}_x^{AB2} = \hat{t}_{AB2}(x,\boldsymbol{\vartheta}_{AB2})$, $\boldsymbol{\zeta}_x^{A1} = \hat{t}_{A1}(x,\boldsymbol{\vartheta}_{A1})$, and $\boldsymbol{\zeta}_x^{B1} = \hat{t}_{B1}(x,\boldsymbol{\vartheta}_{B1})$, then $\boldsymbol{\xi}_x^A = \hat{t}(x,\boldsymbol{\vartheta}_{A1}) + \hat{t}(x,\boldsymbol{\vartheta}_{AB2})$ and $\boldsymbol{\xi}_x^B = \hat{t}(x,\boldsymbol{\vartheta}_{B1}) + \hat{t}(x,\boldsymbol{\vartheta}_{AB2})$, with $\boldsymbol{\vartheta_A}=\{\boldsymbol{\vartheta}_{A1},\boldsymbol{\vartheta}_{AB2}\}$ and $\boldsymbol{\vartheta_B}=\{\boldsymbol{\vartheta}_{B1},\boldsymbol{\vartheta}_{AB2}\}$.
\end{example}

\section{Comparison of Models of the Inversion Method for Discrete Random Variables}\label{app:comparison-tranf-discrete}

We prove that Constraints~\eqref{general_drv_1} and~\eqref{general_drv_2} lead to a tighter linear relaxation feasibility set than the modeling of the inversion method for discrete random variables proposed by~\citet{holzmann2021shortest}. Toward this end, we generalize the constraints proposed by \citet{holzmann2021shortest} to consider the general functions $\hat{p}_{r}^u(x)$. We obtain 
\begin{align}
    \label{app_drv_1}& \theta_{u}^{r} \leq 1 + \sum_{j=1}^{r} \hat{p}_{j}^u(x) - \vartheta_u - \epsilon && \forall\ u \in \mathcal{U},\quad r \in \mathcal{R}_u \\
    \label{app_drv_2}& \theta_{u}^r \leq \vartheta_u + \sum_{j=r}^{R_u} \hat{p}_{j}^u(x)  && \forall\ u \in \mathcal{U},\quad r \in \mathcal{R}_u .
\end{align}
Constraints~\eqref{app_drv_1} serve the same purpose as Constraints~\eqref{general_drv_1}, being active when $\vartheta_u \geq \sum_{j=1}^{r} \hat{p}_{j}^u(x)$. Similarly, Constraints~\eqref{app_drv_2} correspond to Constraints~\eqref{general_drv_2} and are active when $\vartheta_u + \sum_{j=1}^{r} \hat{p}_{j}^u(x) \leq 1$. Therefore, we first prove that the right-hand side of Constraints~\eqref{general_drv_1} is less than that of Constraints~\eqref{app_drv_1} when these constraints are active; that is,
\begin{equation}
     1 + \sum_{j=1}^{r} \hat{p}_{j}^u(x) - \vartheta_u - \epsilon \geq \frac{\sum_{j=1}^{r} \hat{p}_{j}^u(x) - \epsilon}{\vartheta_u} \quad \forall\ u \in \mathcal{U}, \quad r \in \mathcal{R}_u. \nonumber
\end{equation}
We can rewrite these inequalities as:
\begin{equation}
     \label{constr:right-side-1} (\vartheta_u - 1) \left[\sum_{j=1}^{r} \hat{p}_{j}^u(x) - \vartheta_u - \epsilon\right] \geq 0 \quad \forall \ u \in \mathcal{U},\quad r \in \mathcal{R}_u. \nonumber
\end{equation}
These inequalities hold because $(\vartheta_u - 1) \geq 0$ by the definition of $\vartheta_u$, and $\left[\sum_{j=1}^{r} \hat{p}_{j}^u(x) - \vartheta_u - \epsilon\right] \leq 0$ when Constraints~\eqref{general_drv_1} and~\eqref{app_drv_1} are active. Additionally, we demonstrate that the right-hand side of Constraints~\eqref{general_drv_2} is less than that of Constraints~\eqref{app_drv_2} when they are active, as shown in the following relation:
\begin{equation}
     \vartheta_u + \sum_{j=r}^{R_u} \hat{p}_{j}^u(x) \geq  \frac{\sum_{j=r}^{R_u} \hat{p}_{j}^u(x)}{1 - \vartheta_u} \quad \forall\  u \in \mathcal{U},\quad r \in \mathcal{R}_u. \nonumber
\end{equation}
We rewrite these inequalities as:
\begin{equation}
     \vartheta_u \left[ 1 - \vartheta_u - \sum_{j=r}^{R_u} \hat{p}_{j}^u(x) \right] \geq 0 \quad \forall\ u \in \mathcal{U},\quad r \in \mathcal{R}_u. \nonumber
\end{equation}
Thus, these inequalities hold when Constraints~\eqref{general_drv_2} and~\eqref{app_drv_2} are active because $\vartheta_u$ is defined as being greater than or equal to zero, and $\left[1 - \vartheta_u - \sum_{j=r}^{R_u} \hat{p}_{j}^u(x)\right] \geq 0$ when the constraints are active. Finally, the left-hand side of Constraints~\eqref{general_drv_1} and~\eqref{general_drv_2} uses the terms $\pi_{u}^r = \sum_{j=1}^{r}\theta_{u}^{j}$ and $(1 - \pi_{iu}^r) = \sum_{j=r+1}^{R_u}\theta_{u}^{j}$, respectively, instead of the term $\theta_{u}^r$ as in Constraints~\eqref{app_drv_1} and~\eqref{app_drv_2}. This further restricts the feasibility region of the continuous relaxation.

\section{Transformed {\SAA} program of the {\NDFPP}}\label{app:transformed-saa}

We solve the transformed program of the {\NDFPP} with the {\SAA} method. Thus, considering the set of {\SAA} scenarios as $\mathcal{S}_\mathcal{N} = \{1,\ldots,N\}$, the transformed {\SAA} Program~\eqref{Prog:saa} of the {\NDFPP} is defined as
\begin{subequations}
\begin{align}
    \min \mbox{: } \label{appconstr:objctive} &\frac{1}{N}\left( \sum_{s \in \mathcal{S_N}}\sum_{(i,j) \in \mathcal{A}}q_{ij}y_{ij}^{s} \right) &&\\
    \mbox{s.t. } & \eqref{constr:facility-protec}-\eqref{constr:z-domain} && \nonumber \\
    &\eqref{constr:flow-conservation}-\eqref{constr:y-domain} && \forall\ s \in \mathcal{S_N}. \nonumber
\end{align}
\end{subequations}

\section{Adaptation of the Transformations to {\NDFPP}-Binomial and {\NDFPP}-Discrete} \label{app:linear-reformu-stdnorma}

As introduced in Section~\ref{subsec:binomial-problem}, we adapt the convolution for binomial random variables to the {\NDFPP}-Binomial. Thus, the adaptation of Constraints~\eqref{const:bernoulli} to the {\NDFPP}-Binomial leads to the following constraints:
\begin{equation}
    \{\vartheta_{fw}^s  \upsilon^s_{fw} \leq {\phi}_f^{\hat{d}(s)} - \epsilon \ ; \ \upsilon^s_{fw} \geq {\phi}_f^{\hat{d}(s)} - \vartheta_{fw}^s \} \quad\ \forall\ f \in \mathcal{N_F},\quad w \in \{1,\ldots,W\}.
\end{equation}
As stated in Section~\ref{subsec:stdnorma-problem}, for the {\NDFPP}-Discrete, we use the inversion for discrete random variables. Therefore, the adaptation of Constraints~\eqref{general_drv_1}--\eqref{general_drv_4} for the {\NDFPP}-Discrete results in the following constraints:
\begin{align}
    \label{normalization_2}&  \pi_{fw}^s \sum_{w \in \mathcal{W}}u_{\hat{d}(s)f}^w \leq \frac{1}{\vartheta_f^s} \sum_{j=0}^{w} u_{\hat{d}(s)f}^j - \epsilon && \forall\ f \in \mathcal{N_F},\quad w \in \{0,\ldots,W-1\} \\
    \label{normalization_3}&  (1 - \pi_{f,w-1}^s) \sum_{w \in \mathcal{W}}u_{\hat{d}(s)f}^w \leq \frac{1}{(1 - \vartheta_f^s)} \sum_{j=w}^{W} u_{\hat{d}(s)f}^j && \forall\ f \in \mathcal{N_F},\quad w \in \{1,\ldots,W\} \\
    \label{normalization_4}& \pi_{fw}^s = \pi_{f,w+1}^{s} - \theta_{f,w+1}^s && \forall\ f \in \mathcal{N_F},\quad w \in \{0,\ldots,W-1\} \\
    \label{normalization_5}& \sum_{w \in \mathcal{W}} \theta_{fw}^s = 1 && \forall\ f \in \mathcal{N_F}. 
\end{align}
Constraints~\eqref{normalization_2} and~\eqref{normalization_3} contain bilinear terms on binary variables. Hence, we apply the McCormick envelope to linearize these constraints. For this, we include auxiliary variables ${u_T}_f^d \geq 0$ for every facility $f \in \mathcal{N_F}$ and event $d \in \mathcal{D}$. We also define the continuous variables $\tau_{fw}^s \geq 0$ to represent the products $\pi_{fw}^s {u_T}^{\hat{d}(s)}_{f}$ for every $f \in \mathcal{N_F}$, $s \in \mathcal{S}$, and $w \in \{1,\ldots,W-1\}$. We obtain the following inequalities:
\begin{align}
    \label{normalization_1}& {u_T}_f^d = \sum_{w \in \mathcal{W}} u_{df}^w && \forall\ d \in \mathcal{D},  \quad f \in \mathcal{N_F}\\
    \label{mccormick_1}&\tau_{fw}^{s} \leq \hat{u}^{\hat{d}(s)}_f \pi_{fw}^s && \forall\ f \in \mathcal{N_F}, \quad s \in \mathcal{S},\quad w \in \{1,\ldots, W-1\} \\
    \label{mccormick_2}&\tau_{fw}^s \leq {u_T}_f^{\hat{d}(s)} && \forall\  f\in \mathcal{N_F},\quad s \in \mathcal{S},\quad w \in \{1,\ldots, W-1\}\\
    \label{mccormick_3}&\tau_{fw}^s \geq {u_T}_f^{\hat{d}(s)} - \hat{u}^{\hat{d}(s)}_f(1 - \pi_{fw}^s)&& \forall\ f \in \mathcal{N_F},\quad  s \in \mathcal{S},\quad w \in \{1,\ldots, W-1\}, 
\end{align}
where  $\hat{u}^d_f$ is an upper bound of the sum $\sum_{w \in W} u_{df}^w$ and is defined as
\begin{equation}
    \hat{u}^d_f = \sum_{w \in \mathcal{W}} \sum_{i \in \mathcal{N_F}}  \max_{p \in \mathcal{P}} \left( \widetilde{u}_{fi}^{dpw} \right) \quad \forall\ d \in \mathcal{D},\quad f \in \mathcal{N_F}.
\end{equation}
Constraints~\eqref{normalization_1} specify the auxiliary variables $
{u_T}_f^d$ to be the sum of the capacity levels' utility for each node $f \in \mathcal{N_F}$ and event $d \in \mathcal{D}$. Constraints~\eqref{mccormick_1}--\eqref{mccormick_3} are the exact McCormick envelope, which guarantees the correct definition of variables $\tau_{fw}^s$. Finally, Constraints~\eqref{normalization_2} and~\eqref{normalization_3} are replaced by the following inequalities:
\begin{align}
    & \tau_{fw}^s \leq \frac{1}{\vartheta_f^s} \sum_{j=0}^{h} u_{\hat{d}(s)f}^j - \epsilon && \forall\ f \in \mathcal{N_F},\quad s \in \mathcal{S},\quad w \in \{0,\ldots,W-1\} \\
    & {u_T}_f^{\hat{d}(s)} - \tau_{f,w-1}^s \leq \frac{1}{(1 - \vartheta_f^s)} \sum_{j=w}^{W} u_{\hat{d}(s)f}^j && \forall\ f \in \mathcal{N_F},\quad s \in \mathcal{S},\quad w \in \{1,\ldots,W\},
\end{align}
and we obtain a linear reformulation.

\section{{\DEOP} nonconvex formulations of the {\NDFPP}} \label{app:deop-nonconvex}

We follow the notation in Program~\eqref{Prog:discrete} and define the set $K$ of scenarios. Each scenario is composed of a realization of an event and a capacity level realization for each facility $\mathcal{N_F}$. Therefore, the objective function of the {\DEOP} in the extensive form of the {\NDFPP}-Binomial is
\begin{equation}
    \min \sum_{k \in K} \left( \prod_{f \in \mathcal{N_F}}\left( 
    \begin{smallmatrix}
        W \\
        \hat{w}(f,k) 
    \end{smallmatrix} \right) 
    \left({\phi_f^{\hat{d}(k)}}\right)^{\hat{w}(f,k)} \left(1-{\phi_f^{\hat{d}(k)}} \right)^{W - \hat{w} (f,k)} \right) \left( \sum_{(i,j) \in \mathcal{A}} q_{i,j} y_{ij}^k \right),
\end{equation}
for the {\NDFPP}-Discrete is
\begin{equation}
    \min \sum_{k \in K} \left( \prod_{f \in \mathcal{F}} \hat{p}_{\hat{d}(k),f}^{\hat{w}(f,k)}(u) \right) \left( \sum_{(i,j) \in \mathcal{A}} q_{i,j} y_{ij}^k \right) ,
\end{equation}
and for the {\NDFPP}-Normal is
\begin{equation}
    \min \int_{k \in K} \left( \prod_{f \in \mathcal{F}} \frac{1}{\widetilde{\sigma}_f^{\hat{d}(k)}\sqrt{2\pi}} \cdot e^{\frac{-\left(\hat{w}(f,k)-\mu_f^{\hat{d}(k)}\right)}{\left(\widetilde{\sigma}_f^{\hat{d}(k)}\right)^2}} \right) \left( \sum_{(i,j) \in \mathcal{A}} q_{i,j} y_{ij}^k \right) dk.
\end{equation}
Here, $\hat{d}(k)$ and $\hat{w}(f,k)$ are the functions returning the event in $\mathcal{D}$ and the capacity level in $\mathcal{W}$ (discrete distributions) or in $[0,W]$ (continuous distribution) of facility $f \in \mathcal{N_F}$ in scenario $k \in \mathcal{K}$, respectively. In addition, the variables $\phi_f^d$, the functions $\hat{p}_{d,f}^{w}(u)$, and the variables $\mu_f^d$ are defined as in Equations~\eqref{eq:binomial-probability},~\eqref{prob_std_norma}, and~\eqref{constr:normal-mean}, respectively. The first stage variables and constraints are defined as in the transformed program (i.e., we have Constrainsts~\eqref{constr:facility-protec}--\eqref{constr:z-domain}). For the second stage, the variables $y_{ij}^k$ have the same definition as in the transformed program, however, they are indexed by the set $K$ instead of the set $\mathcal{S}$ of scenarios. Thus, we have Constraints~\eqref{constr:flow-conservation}--\eqref{constr:arc-capacity} and~\eqref{constr:y-domain} for every scenario in $K$. 

\section{Expected Value Problems}\label{app:EV}

This section presents the {\EV} for the {\NDFPP}. For this, we define the parameters ${prob}_d$, which are the probability of occurrence of  events $d \in \mathcal{D}$. Table~\ref{tab:disruptions} shows the probability of each disruption, and the probability of no-natural disruption is $
{prob}_{nd} = 1 - \sum_{d \in \mathcal{D}_0} {prob}_d$. We specify a general EV formulation for the {\NDFPP} as:
\begin{align}
    \min \mbox{: } & \sum_{(i,j) \in \mathcal{A}} q_{ij} y_{ij} && \\
    \mbox{s.t. } & \eqref{constr:facility-protec},\eqref{constr:budget},\eqref{constr:x-domain},\eqref{constr:y-domain} && \nonumber \\
    & \sum_{j \in \mathcal{N_T}} y_{jn} - \sum_{j \in \mathcal{N_T}} y_{nj} = b_n  && \forall \n \in \mathcal{N_C} \\
    & \sum_{j \in \mathcal{N_T}} y_{fj} - \sum_{j \in \mathcal{N_T}} y_{jf} \leq \overline{\xi}_x^{f} && \forall\ f \in \mathcal{N_F} \\
    & y_{ij} \leq M z_{e(i,j)} && \forall\ (i,j) \in \mathcal{A}\backslash\mathcal{A}_m \\
    & y_{ij} \geq 0 && \forall\ (i,j) \in \mathcal{A},
\end{align}
where $\overline{\xi}_x^{f}$ is the mean of the endogenous capacity of facility $f \in \mathcal{N_F}$. Variables $x_f^p$ and $z_e$ have the same definition as before, and variables $y_{ij}$ are the flow of goods in arcs $(i,j) \in \mathcal{A}$ in the average case.

For each considered endogenous distribution,  $\overline{\xi}_x^{f}$ is defined differently. In the case of {\NDFPP}-Selection,  the parameter $\zeta_f^p$ is the average capacity value of node $f \in \mathcal{N_F}$ when protection $p \in \mathcal{P}$ is installed. This average value considers the events $d \in \mathcal{D}$. Thus, we have 
\begin{equation}
    \label{constr:ev-transf-discrete-selection}\overline{\xi}_x^{f} = \sum_{p \in \mathcal{P}}\zeta_f^p x_f^p  \ \quad \forall\ f \in \mathcal{N_F}.
\end{equation}
Taking into account  {\NDFPP}-Binomial, we include,  the Constraints~\eqref{eq:binomial-probability} and the variables $
{\phi}^d_f \in [0,1]$ for every $d \in \mathcal{D}$ and $f \in \mathcal{N_F}$, which are introduced in Section~\ref{subsec:binomial-problem}. We specify  $\overline{\xi}_x^{f}$ as follows:
\begin{equation}
    \overline{\xi}_x^{f} = \sum_{d \in \mathcal{D}} {prob}_d  \nu_W {\phi}^d_f \  \quad \forall\ f \in \mathcal{N_F}.
\end{equation}
For  {\NDFPP}-Discrete, we do not use the auxiliary variables $u_{df}^w$,\ but include their definition (see Section~\ref{subsec:stdnorma-problem}) directly in the constraints of the EV problem. We define an auxiliary variable $\varphi^d_f \geq 0$ as the mean of the endogenous distribution considering the event $d \in \mathcal{D}$ and the facility $f \in \mathcal{N_F}$. The $\overline{\xi}_x^{f}$ is then specified as follows:
\begin{align}
    \label{const:def-mean-per-event}&\varphi_f^d = \frac{\sum_{w \in \mathcal{W}} \nu_w \left( \frac{1}{\rho} \sum_{i \in \mathcal{N_F}} \sum_{p \in \mathcal{P}} \widetilde{u}_{fi}^{dpw} x_{i}^p \right)}{\sum_{w \in \mathcal{W}} \left( \frac{1}{\rho} \sum_{i \in \mathcal{N_F}} \sum_{p \in \mathcal{P}} \widetilde{u}_{fi}^{dpw} x_{i}^p \right)} && \forall\ f \in \mathcal{N_F} \\
    &\overline{\xi}_x^{f} = \sum_{d \in \mathcal{D}} {prob}_d  \varphi_f^d   && \forall\ f \in \mathcal{N_F}.
\end{align}
Because the variables $x_f^p$  are binary, we linearize the Constraints~\eqref{const:def-mean-per-event} with the McCormick linearization.

For  {\NDFPP}-Normal, we include the variables $\mu_f^d \geq 0$ for every $d \in \mathcal{D}$ and $f \in \mathcal{N_F}$, with the same definition as in Section~\ref{subsec:normal-problem}, and  Constraints~\eqref{constr:normal-mean}. Thus, we define the $\overline{\xi}_x^{f}$ as follows:
\begin{equation}
    \label{constr:ev-transf-normal}\overline{\xi}_x^{f} = \overline{\nu} \sum_{d \in \mathcal{D}} {prob}_d \mu_f^d  \ \quad \forall\ f \in \mathcal{N_F}.
\end{equation}

\section{Case Study: Summary of Notation} \label{app:detailed-notation}
 Table~\ref{tab:notation} summarizes the notation used in Section~\ref{sec:case-study} for the {\NDFPP}.
\begin{table}[ht]
\TABLE
{Notation of  {\NDFPP}.\label{tab:notation}}
{\begin{adjustbox}{width=.88\textwidth,center}%
\begin{tabular}{ll}
\hline
Sets                   & Definition                                                                                                                                                               \\ \hline
$\mathcal{N_T}$        & Set of all nodes in the network, including a dummy facility (index $j$)                                                                                  \\
$\mathcal{N}_\mathcal{T}^0$        & Set of all nodes in the network without the dummy facility (index $j$)                                                                                  \\
$\mathcal{N_F}$        & Set of facility nodes (index $f$)                                                                                                                                                  \\
$\mathcal{N_C}$        & Set of client nodes (index $n$)                                                                                                                                                   \\
$\mathcal{E}$          & Set of (undirected) edges or links in the network (index $e$)                                                                                                                      \\
$\mathcal{A}$          & Set of (directed) arcs in the network [index $(i,j)$]                                                                                                                                  \\
$\mathcal{A}_m$        & Set of (directed) arcs leaving the dummy facility and going to the other nodes                                                                                        \\
$\mathcal{D}_0$        & Set of disruption types (index $d$)                                                                                                                                               \\
$\mathcal{D}$          & Set of events including the different disruptions and the case of no-disruption (index $d$)                                                                                         \\
$\mathcal{L}$          & Set of disruption intensity levels (index $l$)                                                                                                                                      \\
$\mathcal{P}$          & Set of node protection levels (index $p$)                                                                                                                                       \\
$\mathcal{W}$          & Set of facility capacity levels (index $w$)                                                                                                                                         \\
$\mathcal{S}$          & Set of scenarios (index $s$)                                                                                                                                           \\
$\mathcal{S_N}$          & Set of scenarios for the SAA (index $s$)                                                                                                                                           \\\hline
Parameters             & Definition                                                                                                                                                               \\ \hline
$c_f^p$                & Cost of installing protection level $p \in \mathcal{P}$ at facility $f \in \mathcal{N_F}$                                                                                     \\
$c_e$                & Cost of opening edge $e \in \mathcal{E}$                                                                                                                      \\
$b_n$                  & Demand of client node $n \in \mathcal{N_C}$                                                                                                                                                                      \\
$q_{ij}$               & Per unit transportation cost for flow in arc $(i,j) \in \mathcal{A}$                                                                                                     \\
$m$                    & Dummy facility                                                                                                                                                           \\
$a$                    & Penalty for each unit of unsatisfied demand                                                                                                                              \\
$nd$                   & No-disruption event                                                                                                                                                      \\
$\nu_w$                & Value of facility's capacity at level $w \in \mathcal{W}$                                                                                                                \\
$\overline{\nu}$       & Capacity increase by the increment of one capacity level                                                                                                                 \\
$B$                    & $\sum_{n \in \mathcal{N_C}} b_n$, a big-M value for edge capacity                                                                                                        \\
$C$                    & Maximum cost-budget for edge opening and facility protection allocation                                                                                                      \\ \hline
First-stage variables  & Definition                                                                                                                                                               \\ \hline
$x_f^p$                & Binary variable, equals  1 if protection level $p \in \mathcal{P}$ is installed at facility $f \in \mathcal{N_F}$, and 0 otherwise\\
$z_e^p$                & Binary variable, equals  1 if edge $e \in \mathcal{E}$ is installed with, and 0 otherwise                              \\ \hline
Second-stage variables & Definition                                                                                                                                                                                                \\ \hline
$y_{ij}^s$             & Amount of flow passing through arc $(i,j) \in \mathcal{A}$ in scenario $s \in \mathcal{S}$                                                                                                              \\ \hline
\end{tabular}
\end{adjustbox}
}
{}
\end{table}

\section{Detailed Instance Generation}\label{app:instances}

We follow a similar procedure to that of~\citet{bhuiyan2020stochastic} for network generation. We use the latitude, longitude, and population data of cities in the southeastern United States from the basic US Cities dataset, which is available at \url{https://simplemaps.com/data/us-cities}. First, we select cities with a population exceeding a predefined threshold as the nodes of the network. Then, departing from~\citet{bhuiyan2020stochastic}, we use Delaunay triangulation to define the edges of the network. Thus, we do not restrict the edge lengths and the number of neighbors for nodes. We investigate five networks, with $|\mathcal{N}_\mathcal{T}^0|=15$, 25, 30, 39, and 48. The thresholds and number of edges for these networks are presented in Table~\ref{tab:networks}. The transportation cost $q_{ij}$ of arcs $(i,j) \in \mathcal{A} \backslash \mathcal{A}_m$ is the geometric distance $\varrho_{e(i,j)}$ between cities $i$ and $j$. The demands $b_n$ are set to the node population divided by $10^4$. We set the dummy facility cost multiplier $a$ to $10$. Following~\citet{bhuiyan2020stochastic}, we solve the capacitated facility location and network design problem from~\citet{melkote2001capacitated} to define the set of facilities. For this, we use the median home value from \url{www.zillow.com}, 2023, and we set the node capacity to ${ \sum_{n \in \mathcal{N}_\mathcal{T}^0} b_n }/{0.9 |\mathcal{N}_\mathcal{T}^0|}$. Moreover, similarly as~\citet{bhuiyan2020stochastic}, we generate a maximum cost $c^{\text{max}}_{f}$ for protecting facility $f \in \mathcal{N_F}$ from the continuous uniform distribution $U(7500,15\,000)$. The cost $c_{f}^{p}$ of installing protection $p \in \mathcal{P}$ is ${c^{\text{max}}_{f} p}/{P}$. The cost $c_e$ for opening edge $e \in \mathcal{E}$ is $\varrho_e \overline{b}_{\text{max}}$, where $\overline{b}_{\text{max}}=10$ is the edge-construction cost per unit length.
\begin{table}[ht]
\TABLE
{Details of network instances. \label{tab:networks}}
{\begin{adjustbox}{width=.48\textwidth,center}
\begin{tabular}{rrr||rrr||rrr}
\hline
\up \down $|\mathcal{N}_\mathcal{T}^0|$ &$|\mathcal{E}|$ &Threshold & $|\mathcal{N}_\mathcal{T}^0|$ &$|\mathcal{E}|$ &Threshold & $|\mathcal{N}_\mathcal{T}^0|$ &$|\mathcal{E}|$ &Threshold \\ \hline
\up 15 & 36  & 650\,000  &
25 & 66  & 403\,500  &
30 & 80  & 350\,000  \\
39 & 106 & 300\,000   & 
48 & 131 & 450\,000   \\ \hline
\end{tabular}
\end{adjustbox}}
{}
\end{table}

We follow the same procedure as in~\citet{bhuiyan2020stochastic} to define the intensity level at which facilities are impacted by the disruptions. For each disruption event, facilities situated within a first radius of the disruption's center of occurrence have a (high) disruption intensity level of 1. Facilities beyond this first radius but within a second radius have a (medium) intensity level of 2. The remaining facilities have a (low) intensity level of 3. The center of occurrence, impact radius, and probability of disruptions are shown in Table~\ref{tab:disruptions}. The no-disruption event has $75\%$ probability.
\begin{table}[ht]
\TABLE
{Set of disruption events, adapted from~\citet{bhuiyan2020stochastic}.\label{tab:disruptions}} 
{\begin{adjustbox}{width=.48\textwidth,center} \begin{tabular}{lllll}
\hline
\up \down Disruption & Center         & First radius (km) & Second radius (km) & Probability \\ \hline
\up Hurricane  & Tampa, FL      & 250               & 800                & $10\%$        \\
Snowstorm & Raleigh, NC    & 160               & 800                & $5\%$         \\
\down Tornado    & Huntsville, AL & 80                & 400                & $10\%$        \\ \hline
\end{tabular}
\end{adjustbox}}
{}
\end{table}

\section{Additional Comparative Results for {\NDFPP}-Selection}\label{app:additional-results-comp-baseline}

Table~\ref{tab:linea-saa-seed0} shows the results obtained for the 4- and 5-facility instances generated with a specific seed. We show the number $|K|$ of scenarios, and the optimality gap $GAP$ for the {\DEOP}. All runs of the {\DEOP} reached the time limit of two hours. For the {\SAATP}, we present the estimator of the lower bound $\overline{v}_{M}^N$, its standard deviation $ \hat{\sigma}_{N,M}$, the upper bound $v(\overline{x})$, the over-estimator of the relative optimality gap $\overline{gap}(\overline{x})(\%) = {100\% \times \overline{gap}(\overline{x})}/{v(\overline{x})}$, and the computational time in seconds. In all the runs, the {\SAATP} finished before the time limit of two hours. We also display the gap between the upper-bound values obtained by the {\SAATP} and by the {\DEOP} [$GAP_{UB}(\%)={100\% \times [ UB-v(\overline{x}) ]}/{UB}$], and the gap of the lower bounds [$GAP_{LB}(\%)={100\% \times (LB-\overline{v}_M^N ) }/{LB}$], where $UB$ and $LB$ are the upper and lower bounds, respectively, obtained by the {\DEOP} within the time limit.
\begin{table}[ht]
\TABLE
{Detailed comparison of the {\SAATP} with the {\DEOP} for {\NDFPP}-Selection instances.\label{tab:linea-saa-seed0}}
{\begin{adjustbox}{width=.9\textwidth,center}
\begin{tabular}{rr|rrrrrrrrr|rrrrrrrrr}
\hline \up\down
                              &     & \multicolumn{9}{c|}{4-facility instances}                                                                                                                                                                                          & \multicolumn{9}{c}{5-facility instances}                                                                                                                                                                                                               \\ \cline{3-20} \up\down
                              &     & \multicolumn{2}{c|}{{\DEOP}}       & \multicolumn{5}{c|}{{\SAATP}}                                                                                                & \multicolumn{2}{c|}{Comparison}                                 & \multicolumn{2}{c|}{{\DEOP}}        & \multicolumn{5}{c|}{{\SAATP}}                                                                                                                    & \multicolumn{2}{c}{Comparison}                                  \\ \hline \up\down
$|\mathcal{N}_\mathcal{T}^0|$ & $W$ & $|K|$ & \multicolumn{1}{r|}{$GAP(\%)$} & $\overline{v}_M^N$ & $\hat{\sigma}_{N,M}$ & $v(\overline{x})$ & $\overline{gap}(\overline{x})(\%)$ & \multicolumn{1}{r|}{Time(s)} & $GAP_{UB}(\%)$ & $GAP_{LB}(\%)$ & $|K|$ & \multicolumn{1}{r|}{$GAP(\%)$} & $\overline{v}_M^N$ & $\hat{\sigma}_{N,N}$ & \multicolumn{1}{r|}{$v(\overline{x})$} & $\overline{gap}(\overline{x})(\%)$ & \multicolumn{1}{r|}{Time(s)} & $GAP_{UB}(\%)$ & $GAP_{LB}(\%)$ \\ \hline \up
15                            & 2   & 324   & \multicolumn{1}{r|}{23.5}      & 5099.2             & 14.8                 & 5129.7            & 0.6                     & \multicolumn{1}{r|}{869}    & 0.3                                 & $-29.4$                     & 972   & \multicolumn{1}{r|}{16.0}      & 3474.7             & 12.2                 & \multicolumn{1}{r|}{3479.4}            & 0.1                     & \multicolumn{1}{r|}{1708}  & 0.2                                 & $-18.7$                     \\
                              & 3   & 1024  & \multicolumn{1}{r|}{13.1}      & 4464.6             & 10.6                 & 4476.3            & 0.3                     & \multicolumn{1}{r|}{884}    & 0.0                                 & $-14.8$                     & 4096  & \multicolumn{1}{r|}{14.1}      & 3359.5             & 8.7                  & \multicolumn{1}{r|}{3369.1}            & 0.3                     & \multicolumn{1}{r|}{1096}  & 0.0                                 & $-16.1$                     \\ \down
                              & 4   & 2500  & \multicolumn{1}{r|}{11.1}      & 4321.2             & 7.7                  & 4331.2            & 0.2                     & \multicolumn{1}{r|}{1026}   & 0.0                                 & $-12.2$                     & 12500 & \multicolumn{1}{r|}{21.4}      & 3098.1             & 7.0                  & \multicolumn{1}{r|}{3116.3}            & 0.6                     & \multicolumn{1}{r|}{716}   & 15.2                                & $-7.3$                      \\ \hline \up
25                            & 2   & 324   & \multicolumn{1}{r|}{25.0}      & 7346.9             & 21.3                 & 7403.9            & 0.8                     & \multicolumn{1}{r|}{280}    & 0.0                                 & $-32.3$                     & 972   & \multicolumn{1}{r|}{17.8}      & 4966.1             & 18.7                 & \multicolumn{1}{r|}{4979.3}            & 0.3                     & \multicolumn{1}{r|}{644}   & 0.3                                 & $-20.9$                     \\
                              & 3   & 1024  & \multicolumn{1}{r|}{14.0}      & 6442.3             & 15.4                 & 6452.1            & 0.2                     & \multicolumn{1}{r|}{137}    & 0.1                                 & $-16.0$                     & 4096  & \multicolumn{1}{r|}{15.4}      & 4798.3             & 13.0                 & \multicolumn{1}{r|}{4819.8}            & 0.4                     & \multicolumn{1}{r|}{588}   & 0.7                                 & $-16.8$                     \\ \down
                              & 4   & 2500  & \multicolumn{1}{r|}{11.5}      & 6244.6             & 11.3                 & 6265.1            & 0.3                     & \multicolumn{1}{r|}{168}    & 0.1                                 & $-12.5$                     & 12500 & \multicolumn{1}{r|}{100.0}     & 4405.9             & 10.6                 & \multicolumn{1}{r|}{4415.0}            & 0.2                     & \multicolumn{1}{r|}{210}   & 100.0                               & $-7.3$                      \\ \hline \up
30                            & 2   & 324   & \multicolumn{1}{r|}{23.8}      & 8437.2             & 23.7                 & 8498.9            & 0.7                     & \multicolumn{1}{r|}{452}    & 0.1                                 & $-30.2$                     & 972   & \multicolumn{1}{r|}{18.6}      & 5598.4             & 21.1                 & \multicolumn{1}{r|}{5608.2}            & 0.2                     & \multicolumn{1}{r|}{968}   & 0.8                                 & $-21.6$                     \\
                              & 3   & 1024  & \multicolumn{1}{r|}{13.2}      & 7449.2             & 17.1                 & 7460.3            & 0.1                     & \multicolumn{1}{r|}{658}    & 0.1                                 & $-15.0$                     & 4096  & \multicolumn{1}{r|}{15.9}      & 5413.0             & 14.7                 & \multicolumn{1}{r|}{5429.4}            & 0.3                     & \multicolumn{1}{r|}{879}   & 0.8                                 & $-17.5$                     \\ \down
                              & 4   & 2500  & \multicolumn{1}{r|}{10.9}      & 7224.1             & 12.5                 & 7255.2            & 0.4                     & \multicolumn{1}{r|}{205}    & 0.3                                 & $-11.5$                     & 12500 & \multicolumn{1}{r|}{100.0}     & 4969.9             & 11.9                 & \multicolumn{1}{r|}{4991.0}            & 0.4                     & \multicolumn{1}{r|}{361}   & 100.0                               & $-7.9$                      \\ \hline \up
39                            & 2   & 324   & \multicolumn{1}{r|}{13.0}      & 11954.6            & 24.5                 & 11997.9           & 0.4                     & \multicolumn{1}{r|}{359}    & 0.1                                 & $-14.3$                     & 972   & \multicolumn{1}{r|}{11.8}      & 10237.6            & 24.4                 & \multicolumn{1}{r|}{10241.9}           & 0.0                     & \multicolumn{1}{r|}{193}   & 0.1                                 & $-13.2$                     \\
                              & 3   & 1024  & \multicolumn{1}{r|}{8.5}       & 11379.1            & 16.8                 & 11394.8           & 0.1                     & \multicolumn{1}{r|}{734}    & 0.3                                 & $-8.8$                      & 4096  & \multicolumn{1}{r|}{8.6}       & 9820.8             & 15.4                 & \multicolumn{1}{r|}{9843.7}            & 0.2                     & \multicolumn{1}{r|}{249}   & 0.5                                 & $-8.6$                      \\ \down
                              & 4   & 2500  & \multicolumn{1}{r|}{5.2}       & 10981.8            & 11.2                 & 10999.0           & 0.2                     & \multicolumn{1}{r|}{296}    & 0.3                                 & $-5.0$                      & 12500 & \multicolumn{1}{r|}{100.0}     & 9415.3             & 12.6                 & \multicolumn{1}{r|}{9426.8}            & 0.1                     & \multicolumn{1}{r|}{375}   & 100.0                               & $-4.1$                      \\ \hline \up
48                            & 2   & 324   & \multicolumn{1}{r|}{13.7}      & 13040.6            & 26.6                 & 13089.2           & 0.4                     & \multicolumn{1}{r|}{5850}   & 0.4                                 & $-15.0$                     & 972   & \multicolumn{1}{r|}{11.9}      & 10400.1            & 24.7                 & \multicolumn{1}{r|}{10413}           & 0.1                     & \multicolumn{1}{r|}{1608.0}  & 0.9                                 & $-12.4$                     \\
                              & 3   & 1024  & \multicolumn{1}{r|}{8.4}       & 12331.3            & 18.5                 & 12348.0           & 0.1                     & \multicolumn{1}{r|}{2383}   & 0.2                                 & $-8.8$                      & 4096  & \multicolumn{1}{r|}{8.9}       & 10018.0            & 15.4                 & \multicolumn{1}{r|}{10037.0}           & 0.2                     & \multicolumn{1}{r|}{935}   & 1.2                                 & $-8.3$                      \\ \down
                              & 4   & 2500  & \multicolumn{1}{r|}{5.7}       & 11937.5            & 12.5                 & 11965.2           & 0.2                     & \multicolumn{1}{r|}{2060}   & 0.4                                 & $-5.3$                      & 12500 & \multicolumn{1}{r|}{100.0}     & 9641.9             & 12.2                 & \multicolumn{1}{r|}{9660.4}            & 0.2                     & \multicolumn{1}{r|}{889}   & 100.0                               & $-4.2$                      \\ \hline 
\end{tabular}
\end{adjustbox}
}
{\centering The values correspond to statistics of specific instances generated at the same seed.}
\end{table}

\section{Case Study: Comparison of Formulation Sizes ({\NDFPP}-Selection)}\label{app:vars-constrs}

The nonlinear nonconvex version of the {\DEOP} for the {\NDFPP}-Selection involves $|\mathcal{N_F}| \cdot |\mathcal{P}| + |\mathcal{E}| + |\mathcal{A}| \cdot |K|$ variables, with $|\mathcal{N_F}| \cdot |\mathcal{P}| + |\mathcal{E}|$ binary and $|\mathcal{A}| \cdot |K|$ continuous nonnegative. It also includes $|\mathcal{N_F}| + 1 + |K| \cdot (|\mathcal{N_T}| + |\mathcal{A} \backslash \mathcal{A}_m|)$ constraints, excluding variable domain constraints. Here, $|K|$ is the number of endogenous scenarios. The {\MILP} version of the {\DEOP} proposed by~\citet{bhuiyan2020stochastic}, and discussed in Section~\ref{subsec:results-comparison-baseline},  introduces $|\mathcal{N_F}| \cdot |\mathcal{P}| \cdot (|\mathcal{A}| + 1) \cdot |K|$ additional continuous variables and $|\mathcal{P}| + (|\mathcal{N_F}| - 1) + |\mathcal{N_F}| \cdot |\mathcal{P}| \cdot (3 \cdot |\mathcal{A}| + 1)$ extra constraints due to the linearization of the nonlinear objective function. Moreover, each {\SAATP} problem contains $|\mathcal{N_F}| \cdot |\mathcal{P}| + |\mathcal{E}| + |\mathcal{A}| \cdot |\mathcal{S_N}|$ variables, with the same number of binary variables as the {\DEOP}, and $|\mathcal{N_F}| + 1 + |\mathcal{S_N}| \cdot (|\mathcal{N_T}| + |\mathcal{A} \backslash \mathcal{A}_m|)$ constraints, where $|\mathcal{S_N}|$ refers to the number of {\SAA} scenarios (750 in Section~\ref{subsec:results-comparison-baseline}). To keep the paper self-contained, Table~\ref{tab:numvars} only provides the variable counts for the four- and five-facility instances of the {\NDFPP}-Selection across the three formulations.
\begin{table}[ht]
\TABLE
{Number of variables of the nonlinear {\DEOP}, {\MILP} {\DEOP} and of each problem of the {\SAATP} for the {\NDFPP}-Selection instances.\label{tab:numvars}}
{\begin{adjustbox}{width=.7\textwidth,center} \begin{tabular}{ll|rrr|rrr}
\hline
                  &     & \multicolumn{3}{c|}{Four-facility instances}                                                                                                                & \multicolumn{3}{c}{Five-facility instances}                                                                                                                \\ \hline
{$|\mathcal{N_F}|$} & {$W$} & \multicolumn{1}{l}{{{\DEOP}-Nonlinear}} & \multicolumn{1}{l}{{{\DEOP}-{\MILP}}} & \multicolumn{1}{c|}{{{\SAATP}}} & \multicolumn{1}{l}{{{\DEOP}-Nonlinear}} & \multicolumn{1}{l}{{{\DEOP}-{\MILP}}} & \multicolumn{1}{c}{{{\SAATP}}} \\ \hline \up
{15}                & {2}   & {26940}                                                  & {353532}                                            & {62298}                                          & {26619}                                                  & {429999}                                            & {61551}                                         \\
                  & {3}   & {85044}                                                  & {1461300}                                           & {62302}                                          & {84024}                                                  & {1783864}                                           & {61556}                                         \\
                  & {4}   & {207556}                                                 & {4407556}                                           & {62306}                                          & {205061}                                                 & {5392561}                                           & {61561}                                         \\ \hline
{25}                & {2}   & {49650}                                                  & {648402}                                            & {114828}                                         & {49329}                                                  & {792909}                                            & {114081}                                        \\
                  & {3}   & {156754}                                                & {2679890}                                           & {114832}                                         & {155734}                                                 & {3289174}                                           & {114086}                                        \\
                  & {4}   & {382586}                                                 & {8082586}                                           & {114836}                                         & {380091}                                                 & {9942591}                                           & {114091}                                        \\ \hline
{30}                & {2}   & {60356}                                                  & {787412}                                            & {139592}                                         & {60035}                                                  & {963995}                                            & {138845}                                        \\
                  & {3}   & {190560}                                                 & {3254368}                                           & {139596}                                         & {189540}                                                 & {3998820}                                           & {138850}                                        \\
                  & {4}   & {465100}                                                 & {9815100}                                           & {139600}                                         & {462605}                                                 & {12087605}                                          & {138855}                                        \\ \hline \up
{39}                & {2}   & {80146}                                                  & {1044370}                                           & {185368}                                         & {79825}                                                  & {1280245}                                           & {184621}                                        \\
                  & {3}   & {253050}                                                 & {4316282}                                           & {185372}                                         & {252030}                                                 & {5310590}                                           & {184626}                                        \\ \down
                  & {4}   & {617626}                                                 & {13017626}                                          & {185376}                                         & {615131}                                                 & {16052631}                                          & {184631}                                        \\ \hline \up
{48}                & {2}   & {99287}                                                  & {1292903}                                           & {229643}                                         & {98966}                                                  & {1586126}                                           & {228896}                                        \\
                  & {3}   & {313491}                                                 & {5343379}                                           & {229647}                                         & {312471}                                                 & {6579351}                                           & {228901}                                        \\ \down
                  & {4}   & {765151}                                                 & {16115151}                                          & {229651}                                         & {762656}                                                 & {19887656}                                          & {228906}                                        \\ \hline
\end{tabular}%
\end{adjustbox}
}
{}
\end{table}
While the {\MILP} {\DEOP} has a considerably larger number of variables and constraints compared to the {\SAATP}, the latter has fewer variables and constraints than the nonlinear {\DEOP} for most instances. Additionally, with 750 sampled scenarios, the {\SAATP} achieves extremely low variance (see Result~\ref{result:discretechoice-saa}), and using fewer scenarios may also be efficient depending on the allowable variance, thus reducing the number of variables and constraints of the problems in the {\SAATP}.

\section{Additional Results for {\NDFPP}-Binomial, {\NDFPP}-Discrete, and {\NDFPP}-Normal}\label{app:additional-statistical-estimators}

Table~\ref{tab:saa-results-seed0} shows the results for the 4- and 5 -facility instances generated at a given seed. We present the number $|K|$ of scenarios of the {\DEOP}, the lower-bound estimator $\overline{v}_{M}^N$, its standard deviation $\hat{\sigma}_{N,M}$, the upper-bound estimator $v_{N{'}}(\overline{x})$ for {\NDFPP}-Normal, the upper bound $v(\overline{x})$ for  {\NDFPP}-Binomial and {\NDFPP}-Discrete, the standard deviation $\hat{\sigma}_{N{'}}(\overline{x})$ of the upper-bound estimator for {\NDFPP}-Normal, the over-estimator $\overline{gap}(\overline{x})(\%)$ of the relative optimality gap using either $v_{N{'}}(\overline{x})$ or $v(\overline{x})$, which is ${100 \%\times \overline{gap}(\overline{x})}/{\hat{v}_{N{'}}(\overline{x})} \mbox{ or }{100\% \times \overline{gap}(\overline{x})}/{v(\overline{x})}$, and the computational time in seconds of the {\SAATP}. We set $M=50$ and $N=750$. For  {\NDFPP}-Normal, we use $N{'}=5\times10^4$. In all runs, the {\SAATP} finished before reaching 3 hours.
\begin{table}[ht]
\TABLE
{Detailed results of the {\SAATP}.\label{tab:saa-results-seed0}}
{\begin{adjustbox}{width=.8\textwidth,center}
\begin{tabular}{rrr|rrrrrrr|rrrrrrr}
\hline \up\down
                              &                               &     & \multicolumn{7}{c|}{4-facility instances}                                                                                                                                      & \multicolumn{7}{c}{5-facility instances}                                                                                                                                       \\ \hline \up\down
{\NDFPP}                  & $|\mathcal{N}_\mathcal{T}^0|$ & $W$ & $|K|$ & $\overline{v}_M^N$ & $\hat{\sigma}_{N,M}$ & $v(\overline{x})$ or $\hat{v}_{N{'}}(\overline{x})$ & $\hat{\sigma}_{N{'}}(\overline{x})$ & $\overline{gap}(\overline{x})(\%)$ & Time (s) & $|K|$ & $\overline{v}_M^N$ & $\hat{\sigma}_{N,M}$ & $v(\overline{x})$ or $\hat{v}_{N{'}}(\overline{x})$ & $\hat{\sigma}_{N{'}}(\overline{x})$ & $\overline{gap}(\overline{x})(\%)$ & Time (s) \\ \hline \up
\multirow{15}{*}{Binomial}    & 15                            & 2   & 324   & 5101.6                                 & 21.4                                     & 5130.1                                                                   & -                              & 0.6                                         & 1563                          & 972   & 3458.9                                 & 12.9                                     & 3480.5                                                                   & -                              & 0.6                                         & 1934                         \\
                              &                               & 3   & 1024  & 4459.9                                 & 12.7                                     & 4476.5                                                                   & -                              & 0.4                                         & 1487                          & 4096  & 3358.2                                 & 5.9                                      & 3369.4                                                                   & -                              & 0.3                                         & 1038                         \\ \down
                              &                               & 4   & 2500  & 4318.4                                 & 8.7                                      & 4331.2                                                                   & -                              & 0.3                                         & 1482                          & 12500 & 3095.8                                 & 5.4                                      & 3116.6                                                                   & -                              & 0.7                                         & 1296                         \\ \cline{2-17} \up
                              & 25                            & 2   & 324   & 7350.0                                 & 31.0                                     & 7404.3                                                                   & -                              & 0.7                                         & 481                           & 972   & 4944.5                                 & 19.7                                     & 4981.0                                                                   & -                              & 0.7                                         & 895                          \\
                              &                               & 3   & 1024  & 6435.2                                 & 18.4                                     & 6452.3                                                                   & -                              & 0.3                                         & 308                           & 4096  & 4796.6                                 & 9.3                                      & 4820.7                                                                   & -                              & 0.5                                         & 1018                         \\ \down
                              &                               & 4   & 2500  & 6240.4                                 & 12.6                                     & 6264.9                                                                   & -                              & 0.4                                         & 386                           & 12500 & 4403.0                                 & 8.1                                      & 4415.5                                                                   & -                              & 0.3                                         & 517                          \\ \cline{2-17} \up
                              & 30                            & 2   & 324   & 8440.1                                 & 34.5                                     & 8499.0                                                                   & -                              & 0.7                                         & 703                           & 972   & 5574.5                                 & 22.3                                     & 5610.1                                                                   & -                              & 0.6                                         & 1367                         \\
                              &                               & 3   & 1024  & 7441.5                                 & 20.5                                     & 7460.3                                                                   & -                              & 0.3                                         & 971                           & 4096  & 5411.1                                 & 10.5                                     & 5430.4                                                                   & -                              & 0.4                                         & 1290                         \\ \down
                              &                               & 4   & 2500  & 7219.2                                 & 14.0                                     & 7254.8                                                                   & -                              & 0.5                                         & 530                           & 12500 & 4966.7                                 & 9.1                                      & 4991.7                                                                   & -                              & 0.5                                         & 716                          \\ \cline{2-17} \up
                              & 39                            & 2   & 324   & 11950.5                                & 34.5                                     & 11998.9                                                                  & -                              & 0.4                                         & 671                           & 972   & 10204.4                                & 27.7                                     & 10242.8                                                                  & -                              & 0.4                                         & 412                          \\
                              &                               & 3   & 1024  & 11373.2                                & 20.7                                     & 11395.4                                                                  & -                              & 0.2                                         & 1108                          & 4096  & 9814.0                                 & 12.3                                     & 9843.4                                                                   & -                              & 0.3                                         & 524                          \\ \down
                              &                               & 4   & 2500  & 10981.5                                & 13.2                                     & 10999.0                                                                  & -                              & 0.2                                         & 598                           & 12500 & 9411.9                                 & 9.9                                      & 9426.1                                                                   & -                              & 0.2                                         & 754                          \\ \cline{2-17}  \up
                              & 48                            & 2   & 324   & 13037.8                                & 37.8                                     & 13090.5                                                                  & -                              & 0.4                                         & 8426                          & 972   & 10368.2                                & 26.3                                     & 10415.6                                                                  & -                              & 0.5                                         & 1730                         \\
                              &                               & 3   & 1024  & 12324.7                                & 22.7                                     & 12348.7                                                                  & -                              & 0.2                                         & 3872                          & 4096  & 10020.0                                & 11.9                                     & 10038.2                                                                  & -                              & 0.2                                         & 1344                         \\ \down
                              &                               & 4   & 2500  & 11936.6                                & 14.6                                     & 11965.3                                                                  & -                              & 0.2                                         & 2966                          & 12500 & 9640.8                                 & 9.5                                      & 9661.0                                                                   & -                              & 0.2                                         & 1379                         \\ \hline \up
\multirow{15}{*}{Discrete} & 15                            & 2   & 324   & 5106.9                                 & 15.1                                     & 5131.0                                                                   & -                              & 0.5                                         & 2191                          & 972   & 3475.4                                 & 11.7                                     & 3481.6                                                                   & -                              & 0.2                                         & 2280                         \\
                              &                               & 3   & 1024  & 4460.8                                 & 10.8                                     & 4478.1                                                                   & -                              & 0.4                                         & 3069                          & 4096  & 3360.7                                 & 9.2                                      & 3370.7                                                                   & -                              & 0.3                                         & 2882                         \\ \down
                              &                               & 4   & 2500  & 4321.0                                 & 7.6                                      & 4332.9                                                                   & -                              & 0.3                                         & 3736                          & 12500 & 3100.4                                 & 7.0                                      & 3117.9                                                                   & -                              & 0.6                                         & 3864                         \\ \cline{2-17} \up
                              & 25                            & 2   & 324   & 7357.6                                 & 21.8                                     & 7405.5                                                                   & -                              & 0.6                                         & 1261                          & 972   & 4967.9                                 & 17.8                                     & 4982.6                                                                   & -                              & 0.3                                         & 2038                         \\
                              &                               & 3   & 1024  & 6436.6                                 & 15.6                                     & 6454.6                                                                   & -                              & 0.3                                         & 1653                          & 4096  & 4801.2                                 & 13.5                                     & 4822.7                                                                   & -                              & 0.4                                         & 2813                         \\ \down
                              &                               & 4   & 2500  & 6244.1                                 & 11.1                                     & 6267.2                                                                   & -                              & 0.4                                         & 2302                          & 12500 & 4408.7                                 & 10.8                                     & 4417.5                                                                   & -                              & 0.2                                         & 3047                         \\ \cline{2-17} \up
                              & 30                            & 2   & 324   & 8449.1                                 & 24.2                                     & 8500.4                                                                   & -                              & 0.6                                         & 1539                          & 972   & 5600.3                                 & 20.1                                     & 5611.9                                                                   & -                              & 0.2                                         & 2376                         \\
                              &                               & 3   & 1024  & 7442.6                                 & 17.3                                     & 7462.8                                                                   & -                              & 0.3                                         & 2490                          & 4096  & 5416.2                                 & 15.3                                     & 5432.6                                                                   & -                              & 0.3                                         & 3276                         \\ \down
                              &                               & 4   & 2500  & 7223.4                                 & 12.3                                     & 7257.3                                                                   & -                              & 0.5                                         & 2451                          & 12500 & 4973.2                                 & 12.2                                     & 4993.9                                                                   & -                              & 0.4                                         & 3081                         \\ \cline{2-17} \up
                              & 39                            & 2   & 324   & 11965.7                                & 25.2                                     & 12000.7                                                                  & -                              & 0.3                                         & 1514                          & 972   & 10234.9                                & 23.6                                     & 10245.0                                                                  & -                              & 0.1                                         & 1827                         \\
                              &                               & 3   & 1024  & 11373.2                                & 17.0                                     & 11397.9                                                                  & -                              & 0.2                                         & 2232                          & 4096  & 9820.6                                 & 16.7                                     & 9845.8                                                                   & -                              & 0.3                                         & 2309                         \\ \down
                              &                               & 4   & 2500  & 10981.8                                & 11.1                                     & 11001.4                                                                  & -                              & 0.2                                         & 2496                          & 12500 & 9417.3                                 & 13.1                                     & 9428.3                                                                   & -                              & 0.1                                         & 3948                         \\ \cline{2-17} \up
                              & 48                            & 2   & 324   & 13053.4                                & 27.4                                     & 13092.4                                                                  & -                              & 0.3                                         & 6927                          & 972   & 10401.0                                & 23.2                                     & 10417.6                                                                  & -                              & 0.2                                         & 3736                         \\
                              &                               & 3   & 1024  & 12324.9                                & 18.7                                     & 12351.5                                                                  & -                              & 0.2                                         & 4526                          & 4096  & 10020.5                                & 15.6                                     & 10040.3                                                                  & -                              & 0.2                                         & 3141                         \\ \down
                              &                               & 4   & 2500  & 11937.5                                & 12.4                                     & 11967.9                                                                  & -                              & 0.3                                         & 5103                          & 12500 & 9645.3                                 & 12.5                                     & 9663.1                                                                   & -                              & 0.2                                         & 5151                         \\ \hline \up
\multirow{5}{*}{Normal}       & 15                            & -   & -     & 8604.8                                 & 16.2                                     & 8569.1                                                                   & 14.1       & -0.4                                        & 1637                          & -     & 7133.2                                 & 16.0                                     & 7098.6                                                                   & 12.1       & -0.5                                        & 896                          \\
                              & 25                            & -   & -     & 12564.0                                & 23.3                                     & 12513.1                                                                  & 20.2       & -0.4                                        & 802                           & -     & 10481.0                                & 23.9                                     & 10432.8                                                                  & 18.2       & -0.5                                        & 817                          \\
                              & 30                            & -   & -     & 14194.2                                & 25.8                                     & 14137.7                                                                  & 22.3       & -0.4                                        & 1046                          & -     & 11889.7                                & 26.8                                     & 11835.6                                                                  & 20.3       & -0.5                                        & 1083                         \\
                              & 39                            & -   & -     & 16502.5                                & 28.7                                     & 16444.6                                                                  & 24.6       & -0.4                                        & 1536                          & -     & 15599.8                                & 34.3                                     & 15531.9                                                                  & 25.7       & -0.4                                        & 1339                         \\ \down
                              & 48                            & -   & -     & 18305.6                                & 30.8                                     & 18241.6                                                                  & 26.7       & -0.4                                        & 1704                          & -     & 15758.0                                & 32.5                                     & 15693.6                                                                  & 24.3       & -0.4                                        & 1741                         \\ \hline
\end{tabular}%
\end{adjustbox}
}{\centering The values correspond to statistics of specific instances generated at the same seed.}
\end{table}

\section{Value of the Stochastic Solution}\label{app:vss}

We describe the computation of the relative {\VSS} under estimators and present the detailed results. We have the under estimator $\textrm{{\VSS}}_1={100\% \times [\textrm{{\EEV}} - v(\overline{x})]}/{\textrm{{\EEV}}}$ for {\NDFPP}-Selection, {\NDFPP}-Binomial, and {\NDFPP}-Discrete, for which we can enumerate the {\DEOP} scenarios. For {\NDFPP}-Normal, we compute $\textrm{{\VSS}}_2= 100\% \times {[\overline{\textrm{{\EEV}}} - \hat{v}_{N{'}}(\overline{x})]}/{\overline{\textrm{{\EEV}}}}$, where $\overline{\textrm{{\EEV}}}= \hat{v}_{N{'}}(\hat{x}^{*})$ and $N{'}=5\times10^4$. Table~\ref{tab:vss} summarizes the underestimators of the {\VSS} for the 4- and 5-facility instances.
\begin{table}[ht]
\TABLE
{$\textrm{{\VSS}}_1$ or $\textrm{{\VSS}}_2$ for the different problem variants. \label{tab:vss}}
{\begin{adjustbox}{width=.8\textwidth,center}
\begin{tabular}{r|r|rrrrrr|rrrrrr}
\hline \up\down
\multirow{2}{*}{{\NDFPP}} & \multicolumn{1}{c|}{\multirow{2}{*}{$W$}} & \multicolumn{6}{c|}{4-facility instances}                                                                                                                    & \multicolumn{6}{c}{5-facility instances}                                                                                                                    \\ \cline{3-14} \up\down
                              & \multicolumn{1}{c|}{}                     & $|K|$ & \multicolumn{1}{r}{15-node} & \multicolumn{1}{r}{25-node} & \multicolumn{1}{r}{30-node} & \multicolumn{1}{r}{39-node} & \multicolumn{1}{r|}{48-node} & $|K|$ & \multicolumn{1}{r}{15-node} & \multicolumn{1}{r}{25-node} & \multicolumn{1}{r}{30-node} & \multicolumn{1}{r}{39-node} & \multicolumn{1}{r}{48-node} \\ \hline \up
Selection               & 2                                         & 324   & $30.0\pm0.2$                  & $2.3\pm0.0$                   & $1.6\pm0.0$                   & $2.1\pm0.0$                   & $2.7\pm0.0$                    & 972   & $36.4\pm0.0$                  & $3.0\pm0.0$                   & $52.7\pm0.0$                  & $42.4\pm0.0$                  & $36.4\pm0.0 $                 \\
                              & 3                                         & 1024  & $29.4\pm0.2$                  & $2.6\pm0.0$                   & $1.5\pm0.0$                   & $1.7\pm0.0$                   & $2.4\pm0.0$                    & 4096  & $30.8\pm0.0$                  & $2.2\pm0.0$                   & $51.4\pm0.0$                  & $32.1\pm0.0$                  & $34.0\pm0.0$                  \\ \down
                              & 4                                         & 2500  & $27.2\pm0.1$                  & $2.3\pm0.0$                   & $1.3\pm0.0$                   & $1.5\pm0.0$                   & $2.1\pm0.0$                    & 12500 & $28.2\pm0.0$                  & $2.2\pm0.0$                   & $51.0\pm0.0$                  & $26.7\pm0.0$                  & $31.4\pm0.0$                  \\ \hline \up
Binomial                      & 2                                         & 324   & $1.8\pm0.2$                   & $2.3\pm0.0$                   & $1.5\pm0.0$                   & $2.0\pm0.0$                   & $2.6\pm0.0$                    & 972   & $36.3\pm0.0$                  & $3.0\pm0.0$                   & $3.7\pm0.0$                   & $24.4\pm0.0$                  & $2.4\pm0.0$                   \\
                              & 3                                         & 1024  & $2.1\pm0.2$                   & $2.6\pm0.0$                   & $1.5\pm0.0$                   & $1.7\pm0.0$                   & $2.3\pm0.0$                    & 4096  & $30.7\pm0.0$                  & $2.3\pm0.0$                   & $3.2\pm0.0$                   & $15.5\pm0.0$                  & $2.0\pm0.0$                   \\ \down
                              & 4                                         & 2500  & $1.5\pm0.1$                   & $2.3\pm0.0$                   & $1.3\pm0.0$                   & $1.5\pm0.0$                   & $2.1\pm0.0$                    & 12500 & $28.1\pm0.0$                  & $2.2\pm0.0$                   & $3.4\pm0.0$                   & $11.7\pm0.0$                  & $1.9\pm0.0 $                  \\ \hline \up
Discrete                   & 2                                         & 324   & $1.8\pm0.2$                   & $2.3\pm0.0$                   & $0.7\pm0.0$                   & $1.0\pm0.0$                   & $1.6\pm0.0$                    & 972   & $36.3\pm0.0$                  & $3.0\pm0.0$                   & $3.7\pm0.0$                   & $24.1\pm0.0$                  & $1.2\pm0.0$                   \\
                              & 3                                         & 1024  & $2.1\pm0.2$                   & $2.6\pm0.0$                   & $0.8\pm0.0$                   & $0.7\pm0.0$                   & $1.3\pm0.0$                    & 4096  & $30.7\pm0.0$                  & $2.3\pm0.0$                   & $3.1\pm0.0$                   & $15.0\pm0.0$                  & $1.0\pm0.0$                   \\ \down
                              & 4                                         & 2500  & $1.5\pm0.1$                   & $2.3\pm0.0$                   & $0.6\pm0.0$                   & $0.6\pm0.0$                   & $1.1\pm0.0$                    & 12500 & $28.1\pm0.0$                  & $2.2\pm0.0$                   & $3.1\pm0.8$                   & $11.1\pm0.0$                  & $0.9\pm0.0$                   \\ \hline \up \down
Normal                        & -                                         & -     & $2.9\pm0.1$                   & $6.3\pm0.1$                   & $3.4\pm0.1$                   & $3.1\pm0.1$                   & $3.3\pm0.1$                    & -     & $0.1\pm0.2$                   & $0.2\pm0.2$                   & $0.3\pm0.2$                   & $0.1\pm0.2$                   & $0.5\pm0.2$                   \\ \hline
\end{tabular}%
\end{adjustbox}
}
{\centering The values are averages over the five instances generated per instance configuration, with a 95\% confidence interval.}
\end{table}

\end{APPENDICES}

\end{document}